\documentclass[10pt,journal,compsoc]{IEEEtran}
\usepackage[nocompress]{cite}
\usepackage{amsmath,amssymb}
\usepackage{algorithm}
\usepackage{algorithmic}
\usepackage{subfigure}
\usepackage{graphicx}
\usepackage{textcomp}
\usepackage{epstopdf}
\usepackage{indentfirst}
\usepackage[justification=centering]{caption}
\usepackage{multirow}
\usepackage{multicol}
\usepackage{caption}
\usepackage{amsthm}
\usepackage{color}
\usepackage{tabularx}
\usepackage{booktabs} 
\usepackage{diagbox}
\usepackage{fancyhdr}
\newtheorem{remark}{Remark}
\newtheorem{assumption}{Assumption}
\newtheorem{lemma}{Lemma}

\newtheorem{theorem}{Theorem}
\newtheorem{proposition}{Proposition}
\ifCLASSINFOpdf

\ifCLASSOPTIONcompsoc
  \usepackage[nocompress]{cite}
\else
  \usepackage{cite}
\fi
\ifCLASSINFOpdf
\else
\fi
\begin{document}
%
\title{(Corrected Version) Push-LSVRG-UP: Distributed Stochastic Optimization over Unbalanced Directed Networks with Uncoordinated Triggered Probabilities}
%
%
%
%

\author{Jinhui~Hu,
        Guo~Chen,
        Huaqing~Li,
        Zixiang~Shen,
        and Weidong Zhang
\IEEEcompsocitemizethanks{\IEEEcompsocthanksitem J. Hu, G. Chen, and Z. Shen are with the Department
of Automation, Central South University, Changsha 410083, P.R. China. E-mail: jinhuihu\_csu@163.com; guochen@ieee.org \protect\\
\IEEEcompsocthanksitem H. Li is with Chongqing Key Laboratory of Nonlinear Circuits and Intelligent Information Processing, College of Electronic and Information Engineering, Southwest University, Chongqing 400715, P.R. China. E-mail: huaqingli@swu.edu.cn\protect\\
\IEEEcompsocthanksitem W. Zhang is with School of Information and Communication Engineering, Hainan University, Haikou 570228, Hainan, P.R. China, and also with Department of Automation, Shanghai Jiaotong University, Shanghai 200240, P.R. China. E-mail: wdzhang@sjtu.edu.cn\protect\\
}
\thanks{This work is supported by the National Natural Science Foundation of China under Grant 62073344. (Corresponding author: Guo Chen.)}
}

%
%

%

\IEEEtitleabstractindextext{%
\begin{abstract}
Distributed stochastic optimization, arising in the crossing and integration of traditional stochastic optimization, distributed computing and storage, and network science, has advantages of high efficiency and a low per-iteration computational complexity in resolving large-scale optimization problems. This paper concentrates on resolving a large-scale convex finite-sum optimization problem in a multi-agent system over unbalanced directed networks.
To tackle this problem in an efficient way, a distributed consensus optimization algorithm, adopting the "push-sum" technique and a distributed loopless stochastic variance-reduced gradient (\textit{LSVRG}) method with uncoordinated triggered probabilities, is developed and named \textit{Push-LSVRG-UP}. Each agent under this algorithmic framework performs only local computation and communicates only with its neighbors without leaking their private information.
The convergence analysis of \textit{Push-LSVRG-UP} is relied on analyzing the contraction relationships between four error terms associated with the multi-agent system. Theoretical results provide an explicit feasible range of the constant step-size, a linear convergence rate, and an iteration complexity of \textit{Push-LSVRG-UP} when achieving the globally optimal solution. It is shown that
\textit{Push-LSVRG-UP} achieves the superior characteristics of accelerated linear convergence, fewer storage costs, and a lower per-iteration computational complexity than most existing works. Meanwhile, the introduction of an uncoordinated probabilistic triggered mechanism allows \textit{Push-LSVRG-UP} to facilitate the independence and flexibility of agents in computing local batch gradients. In simulations, the practicability and improved performance of \textit{Push-LSVRG-UP} are manifested via resolving two distributed learning problems based on real-world datasets.
\end{abstract}

}

\maketitle
\thispagestyle{fancy}
\lfoot{}
\cfoot{Citation information: DOI: 10.1109/TNSE.2022.3225229, IEEE Transactions on Network Science and Engineering. Copyright © 20XX IEEE. Permission to use this material for any other purposes must be obtained from the IEEE by sending an email to pubs-permissions@ieee.org. }
\renewcommand{\headrulewidth}{0mm}

\IEEEdisplaynontitleabstractindextext

%
\IEEEpeerreviewmaketitle

\IEEEraisesectionheading{\section{Introduction}\label{Section 1}}
\IEEEPARstart{D}{istributed} optimization has found extensive applications in various fields such as machine learning \cite{Boyd2010b, Nedic2020}, deep learning \cite{Assran2019}, power systems \cite{Li2021,Yang2019}, signal processing \cite{Li2019i}, resource allocation \cite{Lu2020e}, and distributed model predictive control \cite{Shi2018d,Camponogara2011} thanks to its advantages of alleviating the computational burden for the agents, high efficiency for the multi-agent system, and guaranteed privacy for each agent in a peer-to-peer network. However, when facing a category of large-scale optimization problems, distributed batch gradient methods still suffer from a high per-iteration computational complexity result from the local batch gradient computation at each iteration. A way of avoiding such issue is to design stochastic gradient methods. Therefore, this paper aims at studying the following generic finite-sum optimization problem
\begin{equation}\label{E1-1}
\mathop {\min }\limits_{\tilde z \in {\mathbb{R}^n}} \tilde f\left( {\tilde z} \right) := \frac{1}{m}\sum\limits_{i = 1}^m {{f_i}\left( {\tilde z} \right)}, \quad{f_i}\left( {\tilde z} \right) = \frac{1}{{{q_i}}}\sum\limits_{j = 1}^{{q_i}} {{f_{i,j}}\left( {\tilde z} \right)},
\end{equation}
where ${f_i}:{\mathbb{R}^n} \to \mathbb{R}$ is the local objective function and can be further decomposed as $q_i$ component functions $f_{i,j}$ in many machine learning or deep learning problems \cite{Xin2020g,Nedic2020,Assran2019,Ye2020,Xin2020f,Pu2021a}.
The decision variable is $\tilde z$ and the mutual goal of all agents is to seek the optimal solution ${\tilde z^*}$ to problem (\ref{E1-1}) through exchanging information with its neighbors.

\subsection{Literature review}\label{Section 1-1}
\begin{figure}[htp]
\centering
\subfigure[A directed structure.]{\includegraphics[width=1.62in,height=1.52in]{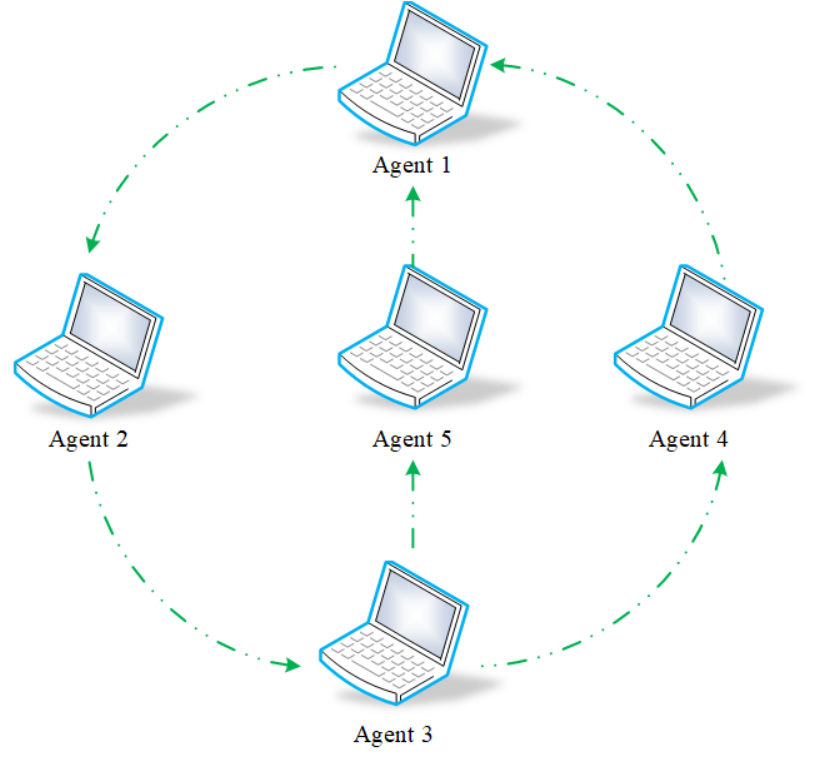}
\label{fig. 0-1}}
\hfil
\subfigure[An undirected structure.]{\includegraphics[width=1.62in,height=1.52in]{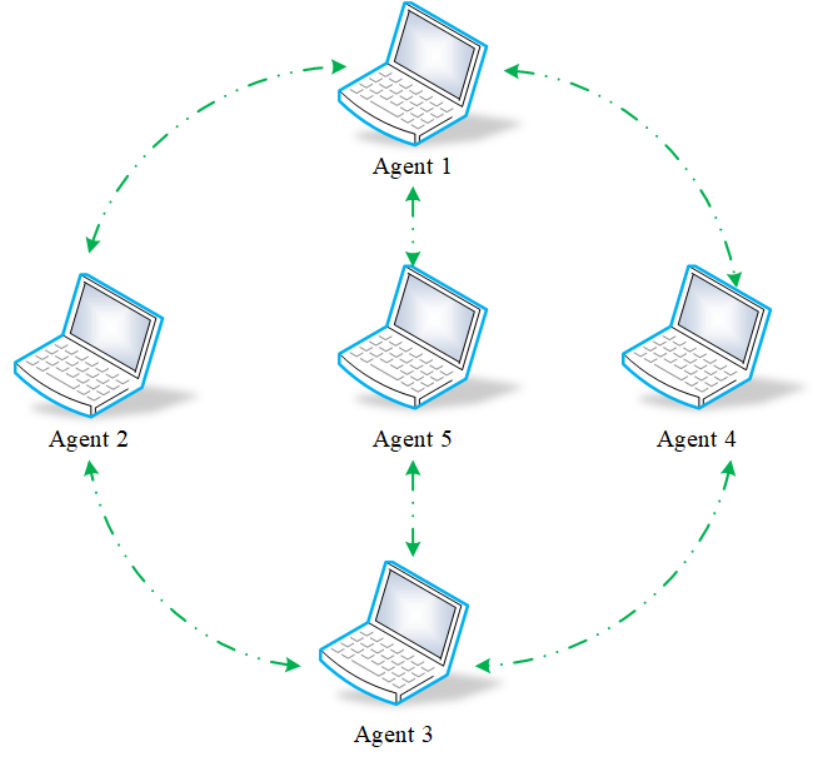}
\label{fig. 0-2}}
\caption{Directed network (a) vs undirected network (b).}
\label{fig. 0}
\end{figure}

Distributed first-order optimization methods can be divided into two categories from the perspective of gradient computation, one of which is the distributed batch gradient methods. Early distributed batch gradient methods include the distributed gradient descent (\textit{DGD}) algorithm \cite{Nedicr2009} and the distributed dual averaging algorithm \cite{J.C.DuchiA.Agarwal}, both of which achieve the globally optimal solution at sub-linear exact convergence rates. Then, \textit{EXTRA} \cite{Shi2015a} adopting a constant step-size achieves linear exact convergence when the local objective functions are strongly convex and have Lipschitz gradients via considering two consecutive gradients of the local objective function. To further facilitate the convergence of \textit{EXTRA}, \textit{DIGing} \cite{Nedic2017a} is designed via adopting the gradient-tracking (GT) technique \cite{Xu2015b}, which is a combine-then-adapt variant of \textit{Aug-DGM} \cite{Xu2015b}. Both \textit{Aug-DGM} \cite{Xu2015b} and \textit{DIGing} \cite{Nedic2017a} are basically two variants of \textit{GT-DGD} methods. Therefore, \cite{Jakovetic2019a} unifies \cite{Shi2015a,Xu2015b,Nedic2017a} into a general primal-dual framework. Nevertheless, the above mentioned distributed algorithms can only work in undirected networks due to the employment of doubly-stochastic weight matrices. In common broadcast-based communication protocols, agents in the system may broadcast at diverse power levels, which indicates the communication capability in one direction while not in the other \cite{Xi2018e}. Here, a simple example that declares the difference between directed communication and undirected communication is exemplified in Fig. \ref{fig. 0}. It is clear from Fig. \ref{fig. 0-1} that each agent in the system may focus on transmitting its information in one direction, while Fig. \ref{fig. 0-2} indicates the necessity of bidirectional information exchange.
Therefore, some outstanding distributed optimization algorithms over unbalanced directed networks are developed.
An earlier work  \cite{Xi2016} incorporates a surplus-based technique into \textit{DGD} \cite{Nedicr2009} to realize a sub-linear exact convergence rate. \textit{DEXTRA} \cite{Xi2017h} combines a "push-sum" method \cite{Nedic2015f} with \textit{EXTRA} \cite{Shi2015a}  to achieve linear exact convergence under the standard strong convexity assumption with the step-size lying in some non-trivial interval. Based on the "push-sum" method and GT technique, \textit{Push-DIGing} \cite{Nedic2017a} and \textit{ADD-OPT} \cite{Xi2018d} employ a column-stochastic weight matrix to achieve accelerated linear exact convergence. \cite{Xi2018e} combines the "push-sum" method and GT technique designs a row-stochastic weight matrix based distributed optimization algorithm. Follow-up papers \cite{Saadatniaki2020,Hu2021c,Pu2021} employ both row- and column-stochastic weight matrices to further explore the generality and novelty in both communication networks and algorithm structures. Although the above mentioned distributed batch gradient methods enjoy fast linear exact convergence under some necessary assumptions, they require each agent in the system computing local batch gradients at each iteration. This inevitably incurs an expensive per-iteration computational complexity for each agent and may also increase the burden of the whole multi-agent system.

Therefore, in another category, distributed stochastic gradient algorithms inspired by centralized stochastic gradient methods receive extensive attention. Based on a decentralized SAGA \cite{Defazio2014c} method, \textit{DSA} \cite{Ribeiro2015} is the first distributed stochastic optimization algorithm, which replaces the batch gradients in \textit{EXTRA} \cite{Shi2015a} with stochastic variance-reduced gradients. Then, \cite{Wang2019} gives an edge-based variant of \textit{DSA}. \textit{S-DIGing} \cite{Li2022} and \textit{GT-SAGA} \cite{Xin2020f} combine the decentralized SAGA method with \textit{GT-DGD} \cite{Xu2015b,Nedic2017a} to achieve accelerated convergence.
Recent work \textit{PMGT-SAGA} \cite{Ye2020} combines a proximal gradient method with a Fast-Mix multi-consensus \cite{Ye2020a} technique to extend \textit{GT-SAGA} for convex composite optimization problems considering a possibly non-smooth term. However, when facing optimization problem (\ref{E1-1}), \textit{SAGA}-based algorithms suffer from an expensive storage cost of $\mathcal{O}\left( {nQ} \right)$ at each iteration. Moreover, when the optimization problem becomes large-scale or high dimensions. That is to say, when $q_i$ or $n$ become larger, the storage cost and the whole multi-agent system can be unbearable. Therefore, \textit{GT-SVRG} \cite{Xin2020f} incorporates a decentralized \textit{SVRG} \cite{Johnson2013c} method into \textit{GT-DGD} \cite{Xu2015b,Nedic2017a}, which is indeed a double-loop distributed stochastic optimization algorithm and the decentralized \textit{SVRG} method enjoys both lower storage costs (almost storage-free) and reduced variance. However, the
double-loop \textit{SVRG} \cite{Johnson2013c} method brings in a global parameter known as the inner-loop iteration number into the distributed algorithm, which requires the multi-agent system paying additional communication to coordinate it in practical unmanned applications. Furthermore, some analytical and estimated issues \cite{Kovalev2019}, especially for more complex networks, are also incurred. Therefore, \cite{Kovalev2019} provides a loopless variant of \textit{SVRG} \cite{Johnson2013c}, which removes the inner-loop iteration in \textit{GT-SVRG} via introducing a probabilistic mechanism. The probabilistic mechanism can trigger the local batch gradient computation in a predefined probability. \textit{LSVRG} is shown in \cite{Kovalev2019} to have the same convergence rate with \textit{SVRG} without requiring any additional assumptions. Then, \textit{PMGT-LSVRG} \cite{Ye2020} extends \textit{LSVRG}\cite{Kovalev2019} to a decentralized setting.

Nevertheless, the above decentralized \textit{SVRG} \cite{Johnson2013c} or \textit{LSVRG} \cite{Kovalev2019} based stochastic algorithms can only work over undirected networks due to the employment of doubly-stochastic weight matrices. As explained before, this may restrict \textit{GT-SVRG} \cite{Xin2020f}  and \textit{PMGT-LSVRG} \cite{Ye2020} in some practical applications with communication capability in one direction while not in the other. Thus to address the issue, this paper devises a distributed stochastic optimization algorithm named \textit{Push-LSVRG-UP}, which employs the "push-sum" technique \cite{Nedic2015f} to cancel the imbalance incurred by unbalanced information exchange, thereby available to handle large-scale optimization problems over unbalanced directed networks. Considering that all agents in the system are restrictive to one common triggered probability, we further introduce the uncoordinated probabilistic triggered mechanism into \textit{Push-LSVRG-UP}
to improve the independence and flexibility of agents. The main contributions of this paper are summarized in the following four aspects.
\subsection{Statement of Contributions}\label{Section 1-2}
\begin{enumerate}
	\item[1)]This chapter designs the first \textit{LSVRG}-based distributed stochastic gradient algorithm \textit{Push-LSVRG-UP} for distributed multi-agent optimization over unbalanced directed networks. Owing to employment of the ``push-sum" technique, \textit{Push-LSVRG-UP} is available to distributed stochastic optimization over unbalanced directed networks, which is more practical than the existing distributed stochastic gradient algorithms \textit{DSA} \cite{Ribeiro2015}, \textit{Diffusion-AVRG}\cite{Yuan2019e}, \textit{S-DIGing} \cite{Li2022}, \textit{GT-SAGA}/\textit{GT-SVRG} \cite{Xin2020f}, \textit{PMGT-SAGA}/\textit{PMGT-LSVRG} \cite{Ye2020}, and \cite{Wang2019} that only available to undirected networks. In theoretical aspects, \textit{Push-LSVRG-UP} is proved to converge linearly to the globally optimal solution under some common assumptions. In simulations, the improved performance of \textit{Push-LSVRG-UP} is shown through making comparisons with existing well-known algorithms when resolving two machine learning problems based on real-world data sets.
	\item[2)]Compared with the expensive storage-required \textit{SAGA-based} distributed stochastic algorithms \textit{DSA} \cite{Ribeiro2015}, \textit{S-DIGing} \cite{Li2022}, \textit{GT-SAGA} \cite{Xin2020f}, \textit{Push-SAGA} \cite{Qureshi2020a}, and \cite{Wang2019}, \textit{Push-LSVRG-UP} reduces a storage cost of $\mathcal{O}\left( {nQ} \right)$, which will be expansive when the dimension of training sample or decision variable is high and the number of training samples is large. Furthermore, from a perspective of theoretical aspects, the main challenge is how to design the contraction relationships of the \textit{LSVRG}-based method over unbalanced directed networks. To overcome this challenge, we present a unified sketch of the proof in Section \ref{Section 4.2}.
	\item[3)]Although there are a large amount of notable distributed batch gradient methods \cite{Nedic2015f,Xi2016,Xi2017h,Nedic2017a,Xi2018d,Xi2018e,Hu2021a,Lu2020e,Pu2021} over unbalanced directed networks, they all suffer from a high per-iteration computational complexity in computing local batch gradients. Especially, when the number of training samples is significantly large, the computational cost for agent $i$ may be unbearable. When compared with \textit{DSGD} \cite{Assran2019,Nedic2016}, \textit{DSGT} \cite{Pu2021a}, $\mathcal{S}$-$\mathcal{A}\mathcal{B}$ \cite{Xin2019c}, and \textit{S-ADDOPT} \cite{Qureshi2021}, \textit{Push-LSVRG-UP} shows its superiority in achieving linear exact convergence. Since these noise-based algorithms \cite{Assran2019,Nedic2016,Xin2019c,Qureshi2021,Pu2021a} do not employ any variance-reduced techniques, they can only converge to the globally optimal solution with sub-linear convergence rates, or achieve linear inexact convergence to an error ball around the globally optimal solution.
	\item[4)] Different with all agents in the system are subject one global triggering probability \cite{Xin2020f, Ye2020}, \textit{Push-LSVRG-UP} adopts an uncoordinated probabilistic triggering mechanism to improve the independence and flexibility of agents in the system. This improvement has certain engineering significance since agents are not necessary to adhere the same coordinated probability to trigger the computation of local batch gradients.
\end{enumerate}

\subsection{Organization}\label{Section 1-3}
The remainder of this paper is organized here.
Some preliminaries including the basic notations, the communication network model, and the problem reformulation are presented in Section \ref{Section 2}. Section \ref{Section 3} develops \textit{Push-LSVRG-UP} and discusses its superior characteristics in contrast to existing distributed optimization algorithms.
The linear convergence rate and the iteration complexity of \textit{Push-LSVRG-UP} are analyzed in Section \ref{Section 4}. Section \ref{Section 5} compares \textit{Push-LSVRG-UP} with existing well-known algorithms based on two different distributed learning problems. We draw a conclusion and state our future work
in Section \ref{Section 6}. Some detailed derivations for the main results are placed in Section \ref{Appendix 1}.
\section{Preliminaries}\label{Section 2}
\subsection{Basic Notations}\label{Section 2-1}
In this section, we give some essential notations that are frequently used in this paper. Note that all vectors are recognized as column vectors if no otherwise stated. some specific definitions are presented in Table \ref{Table 1}.
Notice that nonnegative vectors or matrices indicate all elements of the vectors or matrices are nonnegative.
\begin{table}[!htp]
\centering
\caption{Basic notations.}
\begin{tabularx}{8.8cm}{lX}  
\hline                      
\bf{Symbols}  & \bf{Definitions}  \\
\hline
${\mathbb R}$, ${{\mathbb R}^n}$, ${{\mathbb R}^{m \times n}}$ & the set of real numbers, $n$-dimensional column real vectors, $m \times n$ real matrices, respectively\\
$\mathbb{E}\left[ {{s_k}|{\mathcal{F}_k}} \right]$ & the expectation of a random variable $s_k$ conditioned on a filter $\mathcal{F}_k$\\
$:=$ & the definition symbol\\
${I_n}$   & the $n \times n$ identity matrix \\
$1_m$   & an $m$-dimensional column vector of all ones\\
${x^{\top}}$ & transpose of vector $x$\\
${A^{\top}}$ & transpose of matrix $A$\\
${\rm{diag}}\left\{ x \right\}$ & a diagonal matrix with all the elements of vector $x$ laying on its main diagonal\\
$X \le Y$ & each element in $Y - X$  is nonnegative, where $X$ and $Y$ are two vectors or matrices with same dimensions\\
$X \otimes Y$ & the Kronecker product of matrices $X$ and $Y$\\
$\rho ( X )$ & the spectral radius for matrix $X$\\
$\left\|  \cdot  \right\|_2$ & the standard 2-norm for vectors and matrices\\
${\left\| \cdot \right\|_\pi }$ & a weighted-norm such that $\forall b \in {\mathbb{R}^m}$, ${\left\| b \right\|_\pi } := {\left\| {{{\left( {{\text{diag}}\left\{ {\sqrt \pi  } \right\}} \right)}^{ - 1}}b} \right\|_2}$ or $\forall B\in {\mathbb{R}^{m\times m}}$,
${\left\| B \right\|_\pi } := {\left\| {{{\left( {{\text{diag}}\left\{ {\sqrt \pi  } \right\}} \right)}^{ - 1}}B\left( {{\text{diag}}\left\{ {\sqrt \pi  } \right\}} \right)} \right\|_2}$  \\
\hline
\end{tabularx}
\label{Table 1}
\end{table}
\subsection{Communication Network Model}\label{Section 2-2}
Consider an unbalanced directed network ${\mathcal{G}}{\text{ = }}\left( {\mathcal{V},\mathcal{E}} \right)$, where $\mathcal{V}=\left\{ {1, \ldots ,m} \right\}$ is the set of agents and $\mathcal{E} \subseteq \mathcal{V} \times \mathcal{V}$ is the collected ordered pairs.
Moreover, if $\left( {j,i} \right) \in \mathcal{E}$, there exists ${a_{ji}} > 0$ and ${a_{ji}} = 0$ otherwise. Specifically, for arbitrary two agents, $i,j \in \mathcal{V}$, if ${a_{ji}} > 0$, then agent $i$ can send information to agent $j$ and ${a_{ji}} = 0$ otherwise. The in-neighbors of agent $i$ is denoted as $\mathcal{N}_i^{\text{in}}$, i.e., the set of agents sending information to agent $i$. Similarly, the out-neighbors of agent $i$ is denoted as $\mathcal{N}_i^{\text{out}}$, i.e., the set of agents receiving information from agent $i$.
The network $\mathcal{G}$ is considered to be balanced if $\sum\nolimits_{j \in \mathcal{N}_i^{\text{out}}} {{a_{ji}} = } \sum\nolimits_{j \in \mathcal{N}_i^{\text{in}}} {{a_{ij}}} $, $i \in \mathcal{V}$, and unbalanced otherwise. Both $\mathcal{N}_i^{\text{in}}$ and $\mathcal{N}_i^{\text{out}}$ include agent $i$.
\begin{assumption}\label{A1}
The weight matrix $\underline{A} = {\left[ {{a_{ij}}} \right]_{1 \le i,j \le m}} \in {\mathbb{R}^{m \times m}}$ associated with the unbalanced directed network ${\mathcal{G}}$ is primitive and column-stochastic, which means that there exists an integer $K>0$ such that $\underline{A}^K$ is a positive matrix and $1_m^{\top}\underline{A} = 1_m^{\top}$.
\end{assumption}
\begin{remark}\label{R2-0}
A feasible way of designing the weight matrix can be found in \cite[Remark 2]{Pu2021}. Moreover, under Assumption \ref{A1}, it is straightforward from \cite[Theorem 8.5.1]{Johnson2013} that the weight matrix $\underline{A}$ has a unique positive right eigenvector $\pi$ with respect to eigenvalue $1$, with $1_m^{\top}\pi  = 1$ and ${\underline{A}_\infty } := {\lim _{k \to \infty }}{\underline{A}^k} = \pi 1_m^{\top}$.
Let $\bar {\pi }$ and $\underline {\pi }$ denote the maximum element and the minimum element in vector ${\pi }$, respectively. For any vector $x \in  {\mathbb{R}^{n}}$, it can be derived according to the norm equivalence property that ${\left\| {x} \right\|_\pi } \le {\underline {\pi} ^{ - 0.5}}{\left\| {x} \right\|_2}$ and ${\left\| {x} \right\|_2} \le {{\bar \pi }^{0.5}}{\left\| {x} \right\|_\pi }$. Moreover, it can be verified that ${\left\| \underline{A} \right\|_\pi } = {\left\| {{\underline{A}_\infty }} \right\|_\pi } = {\left\| {{I_{m}} - {\underline{A}_\infty }} \right\|_\pi } = 1$ under Assumption \ref{A1}.
\end{remark}
\subsection{Problem Reformulation}\label{Section 2-3}
To resolve problem (\ref{E1-1}) in a decentralized manner, we introduce $z^i$, $i \in \mathcal{V}$, as local copies of decision variable $\tilde z$, and reformulate problem (\ref{E1-1}) as follows:
\begin{equation}\label{E2-1}
\begin{aligned}
\!\!\!\!\!\!&\mathop {\min }\limits_{z \in {\mathbb{R}^{mn}}} f\left( z \right) :=  \frac{1}{m}\sum\limits_{i = 1}^m {{f_i}\left( {{z^i}} \right)}, \quad {f_i}\left( {{z^i}} \right) = \frac{1}{{{q_i}}}\sum\limits_{j = 1}^{{q_i}} {{f_{i,j}}\left( {{z^i}} \right)}, \\
\!\!\!\!\!\!&\;\,{\text{s}}{\text{.t}}{\text{. }}{z^i} = {z^j},\left( {i,j} \right) \in \mathcal{E}.
\end{aligned}
\end{equation}
In the sequel, we use ${q_{\max }} := {\max _{i \in \mathcal{V}}}{q_i}$ and ${q_{\min }} := {\min _{i \in \mathcal{V}}}{q_i}$ to denote respectively the maximum number and the minimum number among local samples. Let ${\mathcal{Q}_i}: = \left\{ {1,2, \ldots ,{q_i}} \right\}$.
\begin{assumption}($\mu$-strongly convex \cite[Assumption 1]{Xin2020f})\label{A2}
For $i \in \mathcal{V}$, each local objective functions ${{f_i}}$ is $\mu$-strongly convex, such that $\forall a,b \in {\mathbb{R}^n}$, we have
\begin{equation}\label{E2-2}
{\mu _i}\left\| {x - y} \right\|_2^2 \le {\left( {\nabla {f_i}\left( x \right) - \nabla {f_i}\left( y \right)} \right)^ \top }\left( {x - y} \right),
\end{equation}
where $\mu>0$.
\end{assumption}
\begin{assumption}($L$-smoothness \cite[Assumption 2]{Xin2020f})\label{A3}
For $i \in {\mathcal V}$, each component function $f_{i,h}$, $h \in {\mathcal{Q}_i}$, has a Lipschitz continuous gradient, such that $\forall a,b \in {\mathbb{R}^n}$, there exists
\begin{equation}\label{E2-3}
\| {\nabla f_{i,h}( a ) - \nabla f_{i,h}( b )} \|_{2} \le {L}\| {a - b} \|_{2},
\end{equation}
 where $L > 0$.
\end{assumption}
\begin{remark}\label{R2-1}
Note that Assumptions \ref{A1}-\ref{A3} are not uncommon in recent literature \cite{Saadatniaki2020,Qureshi2020a,Pu2021}. Furthermore, it can be obtained from \cite[chaper 3]{Bubeck2015c} that $0< \mu \le L$.  Under Assumption \ref{A2}, we know that the globally optimal solution ${z^*} \in {\mathbb{R}^{mn}}$ to problem (\ref{E2-1}) exists uniquely. If we define the category of $\mu$-strongly convex and $L$-smooth functions as ${\mathcal {S}}$, then ${f_i} \in {\mathcal {S}}$, $i \in {\mathcal {V}}$, and thus it can be verified that $\tilde f \in {\mathcal {S}}$.
\end{remark}
\section{Algorithm Development}\label{Section 3}
Consider the noise-based distributed stochastic gradient algorithm \textit{S-ADDOPT} \cite{Qureshi2021}:
\begin{subequations}\label{E3-1-0}
\begin{align}
\label{E3-1-0-1}&x_{k + 1}^i = \sum\limits_{j \in \mathcal{N}_i^{{\text{in}}}} {{a_{ij}}x^j_k}  - \alpha v^i_k,\\
\label{E3-1-0-2}&y^i_{k + 1} =  \sum\limits_{j \in \mathcal{N}_i^{{\text{in}}}} {{a_{ij}}y^j_k},\\
\label{E3-1-0-3}&z^i_{k + 1} =  \frac{{{x^i_{k + 1}}}}{{{y^i_{k + 1}}}},\\
\label{E3-1-0-4}&{v^i_{k + 1}} = \sum\limits_{j \in \mathcal{N}_i^{{\text{in}}}} {{a_{ij}}{v^j_k}}  + \nabla f\left( {{x^i_{k + 1}},{\xi ^i_{k + 1}}} \right) - \nabla f\left( {{x^i_k},{\xi ^i_k}} \right),
\end{align}
\end{subequations}
where $x^i_k$ and $y^i_k$ are auxiliary variables, $z^i_k$ is the decision variable, $v^i_k$ is the gradient tracker, and $\nabla f\left( {{x^i_k},{\xi ^i_k}} \right)$ is the noisy gradient with a bounded random noise ${\xi ^i_k}$. Therein, (\ref{E3-1-0-1})-(\ref{E3-1-0-4}) represent the gradient-descent, (left) Perron eigenvectors estimation, bias correction, and GT steps, respectively. Steps (\ref{E3-1-0-1})-(\ref{E3-1-0-3}) are referred to as "push-sum" updates. Unfortunately, the convergence result of \textit{S-ADDOPT} is affected by the random noise term, leading to either inlinear exact convergence or sub-linear exact convergence. This implies that \textit{S-ADDOPT} fails to achieve a faster exact convergence rate, which is achievable for distributed batch gradient algorithms \cite{Xi2017h,Nedic2017a,Xi2018d,Xi2018e,Hu2021a,Lu2020e,Pu2021}.

\begin{algorithm}[!htp]
	\caption{\textit{Push-LSVRG-UP}: An LSVRG-based distributed stochastic gradient algorithm with uncoordinated triggered probabilities.}
	\label{Algo1}
	\renewcommand{\algorithmicrequire}{\textbf{Input:}}
	\renewcommand{\algorithmicensure}{\textbf{Output:}}
	\begin{algorithmic}[1] 
		\REQUIRE a proper constant step-size $\alpha >0$ and uncoordinated triggered probabilities $p_i \in \left( {0,1} \right]$.
		\STATE {\bf{Initialize:}} the decision variables $z_0^i = x_0^i \in {\mathbb{R}^n}$, $w_0^i = z_0^i$, $v_0^i = g_0^i = \nabla f_i\left( {z_0^i} \right)$, and $y_0^i=1$, together with the weights ${a_{ji}} \ge 0$ such that $\sum\nolimits_{j \in \mathcal{N}_i^{{\text{out}}}} {{a_{ji}}}  = 1$, $ j \in \mathcal{N}_i^{{\text{out}}}$.
		\FORALL {$k = 0, 1, 2, \ldots$}
		\FORALL {nodes $i \in \mathcal{V}$ in parallel}
		\STATE {\bf{Local stochastic gradient estimation:}} Choosing $s^i_k$ uniformly randomly from the local sample set ${\mathcal{Q}_i}$, and then evaluate the stochastic gradient according to $g^i_{k + 1} = \nabla {f_i^{s^i_k}}\left( {z^i_k} \right) - \nabla {f_i^{s^i_k}}\left( {w^i_k} \right) + \nabla {f_i}\left( {w^i_k} \right)$;
		\label{step 1}
		\STATE {\bf{Triggering mechanism:}} Following a heterogeneous probabilistic triggering mechanism: $w^i_{k + 1} = \left\{ \begin{gathered}
		z^i_k,{p_i} \hfill \\
		w^i_k,1 - {p_i} \hfill \\
		\end{gathered}  \right.$;
		\label{step 2}
		\STATE {\bf{Local gradient descent:}} $x^i_{k + 1} = \sum\limits_{j \in \mathcal{N}_i^{{\text{in}}}} {{a_{ij}}x^j_k}  - \alpha v^i_k$;
		\label{step 3}
		\STATE {\bf{Eigenvector vector estimation:}} $y^i_{k + 1} = \sum\limits_{j \in \mathcal{N}_i^{{\text{in}}}} {{a_{ij}}y^j_k} $;
		\label{step 4}
		\STATE {\bf{State transformation:}} $z^i_{k + 1} = \frac{{x^i_{k + 1}}}{{y^i_{k + 1}}}$;
		\label{step 5}
		\STATE {\bf{GT step:}} $v^i_{k + 1} = \sum\limits_{j \in \mathcal{N}_i^{{\text{in}}}} {{a_{ij}}v^j_k}  + g^i_{k + 1} - g^i_k$.
		\label{step 6}
		\ENDFOR
		\ENDFOR
		\ENSURE{all decision variables ${x^i_k}$.}
	\end{algorithmic}
\end{algorithm}

To reach this goal, we need to remove the variance ($\sigma ^2$) in the convergence error as follows:
\begin{equation}\label{E2-3-2}
\mathop {\lim }\limits_{k \to \infty } \sup {e_k} := \alpha \mathcal{O}\left( {\frac{{{\sigma ^2}}}{{m\mu }}} \right) + {\alpha ^2}\mathcal{O}\left( {\frac{{{L^2}{\sigma ^2}}}{{{\mu ^2}{{\left( {1 - \sigma _A^2} \right)}^4}}}} \right),
\end{equation}
where the variance $\sigma ^2$ dominates the asymptotic convergence error ${e_k}$. In view of this, this section aims to develop a distributed \textit{LSVRG}-based stochastic optimization algorithm to remove this variance for the finite-sum optimization problem (\ref{E2-1}). The reason that we adopt \textit{LSVRG}-type variance-reduced strategy here lies in the fact that in contrast to \textit{SAGA}-based methods \textit{DSA} \cite{Ribeiro2015}, \textit{GT-SAGA} \cite{Xin2020f}, and \textit{Push-SAGA} \cite{Qureshi2020a}, it does not require an expensive storage cost $\mathcal{O}\left( {nQ} \right)$ thanks to its probabilistic nature of computing the local batch gradient. Inspired by the "push-sum" method \cite{Nedic2017a,Xi2018d}, the GT technique \cite{Xu2015b} and the decentralized \textit{LSVRG} method \cite{Ye2020}, we design \textit{Push-LSVRG-UP} to resolve large-scale optimization problems over unbalanced directed networks and the execution details of \textit{Push-LSVRG-UP} are presented in Algorithm \ref{A1}.
\begin{remark}\label{R3-1}
Since the inner loop number in \textit{SVRG}-based distributed stochastic optimization algorithms, for example \textit{GT-SVRG} \cite{Xin2020f}, is a global parameter, it is adverse to the distributed implementation. Thus,
different with the stochastic double-loop distributed algorithm \textit{GT-SVRG} \cite{Xin2020f}, \textit{Push-LSVRG-UP} described in Algorithm \ref{Algo1} not only removes the inner loop via applying a probabilistic triggered mechanism, but is available to work in a class of generic unbalanced directed networks. Another advantage of \textit{Push-LSVRG-UP} in contrast to other \textit{LSVRG}-based algorithms \cite{Kovalev2019,Gorbunov2019,Ye2020}, is the introduction of the uncoordinated triggered probabilistic mechanism, which not only improves the independence and flexibility of each agent $i$ in optimization procedures, but is helpful for its distributed execution. However, this incorporated mechanism also incurs some challenges in deriving the explicit linear convergence and an iteration complexity of \textit{Push-LSVRG-UP}, which has been well-addressed in Section \ref{Section 4}.
\end{remark}
\begin{remark}\label{R3-2}
(A lower per-iteration computational complexity).
In the algorithmic of distributed batch gradient algorithms, such as \cite{Nedic2015f,Xi2016,Xi2017h,Nedic2017a,Xi2018d,Xi2018e,Hu2021a,Lu2020e,Pu2021}, each agent $i$, $i \in \mathcal{V}$, suffers from a per-iteration computational complexity of $\mathcal{O}\left( {{q_i}} \right)$ for computing the local batch gradients. Nevertheless, \textit{Push-LSVRG-UP} inherits a merit from the centralized \textit{LSVRG}-based method \cite{Kovalev2019}, which allows each agent $i$ to calculate only two component gradients $\nabla {f_{i,s_{k + 1}^i}}\left( {z_{k + 1}^i} \right)$ and $\nabla {f_{i,s_{k + 1}^i}}\left( {w_k^i} \right)$, and thus the corresponding per-iteration computational complexity of \textit{Push-LSVRG-UP} is $\mathcal{O}\left( 1 \right)$ if the local batch gradient-computation is not triggered. This reduction in the per-iteration computational complexity is more significant when $q_i$ becomes larger and the uncoordinated triggered probabilities satisfy $0<p_i < 1$. %
\end{remark}
For the convenience of the subsequent convergence analysis, some important definitions are given by $\forall k \ge 0$:
\begin{equation*}
\begin{aligned}
{x_k} =& {\left[ {{{\left( {x_k^1} \right)}^ \top }, {{\left( {x_k^2} \right)}^ \top }, \ldots ,{{\left( {x_k^m} \right)}^ \top }} \right]^ \top },\\
{y_k} =& {\left[ {y_k^1, y_k^2, \ldots ,y_k^m} \right]^{\top} },\\
{z_k} =& {\left[ {{{\left( {z_k^1} \right)}^ \top }, {{\left( {z_k^2} \right)}^ \top }, \ldots ,{{\left( {z_k^m} \right)}^ \top }} \right]^ \top },\\
{g_k} =& {\left[ {{{\left( {g_k^1} \right)}^ \top }, {{\left( {g_k^2} \right)}^ \top }, \ldots ,{{\left( {g_k^m} \right)}^ \top }} \right]^ \top },\\
{v_k} =& {\left[ {{{\left( {v_k^1} \right)}^ \top },{{\left( {v_k^2} \right)}^ \top },  \ldots ,{{\left( {v_k^m} \right)}^ \top }} \right]^ \top },\\
{A} =& \underline{A} \otimes {I_n}, \\
{A_\infty } =& {\lim _{k \to \infty }}{A^k} = \left( {\pi 1_m^{\top}} \right) \otimes {I_n}, \\
{Y_k} =& {\text{diag}}\left\{ {{y_k}} \right\} \otimes {I_n},\\
Y =& {\text{su}}{{\text{p}}_k}{\left\| {{Y_k}} \right\|_2} \ge 1, \\
\tilde Y =& {\text{su}}{{\text{p}}_k}{\left\| {Y_k^{ - 1}} \right\|_2} \ge 1.\\
\end{aligned}
\end{equation*}
Based on the above definitions, we now give the vector-matrix form of Algorithm \ref{Algo1} as follows:
\begin{subequations}\label{E3-1}
\begin{align}
\label{E3-1-1}{x_{k + 1}} =& A{x_k} - \alpha {v_k},\\
\label{E3-1-2}{y_{k + 1}} =& \underline{A} {y_k},\\
\label{E3-1-3}{z_{k + 1}} =& Y_k^{ - 1}{x_{k + 1}},\\
\label{E3-1-4}{v_{k + 1}} =& A{v_k} + {g_{k + 1}} - {g_k}.
\end{align}
\end{subequations}
Based on (\ref{E3-1}), the following useful notations are defined
\begin{equation*}
\begin{aligned}
&{\bar x_k} = \frac{1}{m}\left( {1_m^{\top} \otimes {I_n}} \right){x_k},\\
\end{aligned}
\end{equation*}
\begin{equation*}
\begin{aligned}
&{\bar v_k} = \frac{1}{m}\left( {1_m^{\top} \otimes {I_n}} \right){v_k},\\
&\nabla F\left( {{z_k}} \right) \!=\! {\left[ {{{\left( {\nabla {f_1}\left( {z_k^1} \right)} \right)}^{\top}}, {{\left( {\nabla {f_2}\left( {z_k^2} \right)} \right)}^{\top}}, \ldots ,{{\left( {\nabla {f_m}\left( {z_k^m} \right)} \right)}^{\top}}} \right]^{\top}},\\
&{{\bar h}_k} = \frac{1}{m}\left( {1_m^{\top} \otimes {I_n}} \right)\nabla F\left( {{z_k}} \right),\\
&{{\bar g}_k} = \frac{1}{m}\left( {1_m^{\top} \otimes {I_n}} \right){g_k},\\
&{{\tilde p}_k} = \frac{1}{m}\left( {1_m^{\top} \otimes {I_n}} \right)\nabla F\left( {{1_m} \otimes {{\bar x}_k}} \right).
\end{aligned}
\end{equation*}

\section{Convergence Analysis}\label{Section 4}
In this section, the iteration complexity and the linear convergence rate of Algorithm \ref{Algo1} are derived. Moreover, the step-size condition is also provided when \textit{Push-LSVRG-UP} converges linearly to the globally optimal solution.
In what follows, let $\bar p$ and $\underline{p}$ represent respectively the maximum value and minimum value of the uncoordinated triggered probabilities $p_i$, $i \in \mathcal{V}$, where $0<\underline{p}\le \bar p \le 1$.
\subsection{Main Results}\label{Section 4.1}
Before presenting the main results of this paper, we denote the condition number of functions in ${\mathcal {S}}$ as $\mathcal{Q} := L/\mu$ such that $\mathcal{Q} \ge 1$ and define the matrix norm $\sigma_A  := {\left\| {A - {A_\infty }} \right\|_\pi }$ such that $0 < \sigma_A  < 1$ can be guaranteed under Assumption \ref{A1} (see \cite[Lemma 1]{Xin2019c} and \cite[Section IV]{Qureshi2020a} for details).
\begin{theorem}\label{Theorem 1}
Suppose that Assumptions \ref{A1}-\ref{A3} hold. Considering Algorithm \ref{Algo1} and for a directivity constant $\delta  \ge 1$ defined in Lemma \ref{L3}, if the step-size satisfies
\begin{equation}\label{E4-1}
0 < \alpha  \le \min \left\{\frac{{ \left( {1 - \sigma_A } \right)\underline{p}}}{{6\mu }}, {\frac{{ {{\left( {1 - \sigma_A } \right)}^2}\underline{p}}}{{480\delta L\mathcal{Q}\bar p}}} \right\},
\end{equation}
then the sequence ${\left\{ {{z_k}} \right\}_{k \ge 0}}$ generated by Algorithm \ref{Algo1} converges linearly to the optimal solution ${{{\tilde z}^*}}$ at the rate of $\mathcal{O}( {{{\left( {\eta  + \zeta } \right)}^k}} )$, where $0<\eta<1$ is defined in Lemma \ref{L7} and $\zeta$ is an arbitrarily small positive constant such that $0 < \eta  + \zeta  < 1$. This means that \textit{Push-LSVRG-UP} achieves an $\epsilon$-accurate solution, i.e., $ \mathbb{E}\left[ {\left\| {{z_k} - {{1_m} \otimes {{\tilde z}^*}} } \right\|_2^2} \right] \le \epsilon$ in at least
\begin{equation}\label{E4-2}
k \ge \mathcal{O}\left( {\max \left\{ {\frac{1}{{\left( {1 - \sigma_A } \right)\underline{p} }},\frac{{\delta {\mathcal{Q}^2}\bar p}}{{{{\left( {1 - \sigma_A } \right)}^2}\underline{p} }}} \right\}\ln \frac{1}{\epsilon}} \right)
\end{equation}
iterations (component gradient computations) at each agent.
\begin{proof}
The detailed proof of Theorem \ref{Theorem 1} is placed in Section \ref{Appendix 1} of Appendix to enhance coherence of the paper.
\end{proof}
\end{theorem}
\begin{remark}\label{R4-1}
One may be aware that the step-size condition (\ref{E4-1}) contains the network information $\sigma_A$, which is indeed a global information engendered by the conservative convergence analysis. However, the fully distributed running of \textit{Push-LSVRG-UP} can be still available via exerting a notion of sufficiently small but positive step-size, which is not uncommon in literature \cite{Xi2018d,Xi2017h}.
\end{remark}

\begin{table}[!htp]
\centering
\caption{Convergence performance comparison.}
\resizebox{8.8cm}{2.1cm}{
\begin{tabular}{cc}
  \toprule
  \bf{Algorithm} & \bf{Convergence Rate} \\
  \midrule
  \textit{DSGT} \cite{Pu2021a}, $\mathcal{S}$-$\mathcal{A}\mathcal{B}$ \cite{Xin2019c} and \textit{S-ADDOPT} \cite{Qureshi2021}
  & linear but inexact convergence \\
  \textit{DSA} \cite{Ribeiro2015}
  & $\mathcal{O}\left( {\max \{ {{q_{\max }}Q,\frac{{{Q^4}}}{{1 - \sigma_A }},\frac{1}{{{{\left( {1 - \sigma_A } \right)}^2}}}} \}\ln \frac{1}{\epsilon}}\right)$ \\
  \textit{Push-SAGA} \cite{Qureshi2020a}
  & $\mathcal{O}\left( {\max \{ {{q_{\max }},\frac{{{q_{\max }}}}{{{q_{\min }}}}\frac{{{\delta Q^2} }}{{{{\left( {1 - \sigma_A } \right)}^2}}}} \}\ln \frac{1}{\epsilon}} \right)$ \\
  \textit{GT-SAGA} \cite{Xin2020f}
  & $\mathcal{O}\left( {\max \{ {{q_{\max }},\frac{{{q_{\max }}{Q^2}}}{{{q_{\min }}{{\left( {1 - \sigma_A } \right)}^2}}}} \}\ln \frac{1}{\epsilon}} \right)$ \\
  \textit{GT-SVRG} \cite{Xin2020f}
  & \begin{tabular}[c]{@{}l@{}} $\mathcal{O}\left( {( {{q_{\max }} + \frac{{{Q^2}\log Q}}{{{{\left( {1 - {\sigma_A ^2}} \right)}^2}}}} )\ln \frac{1}{\epsilon}} \right)$ \\ with inner-loop number $T = \mathcal{O}\left( {\frac{{{Q^2}\log Q}}{{{{\left( {1 - \sigma_A } \right)}^2}}}} \right)$ \end{tabular} \\
  \textit{ADD-OPT} \cite{Xi2018d}, \textit{Push-Pull} \cite{Pu2021}, \textit{Diffusion-AVRG} \cite{Yuan2019e}, and \cite{Wang2019}
  & linear (no explicit convergence rate) \\
  \textit{PMGT-SAGA}/\textit{PMGT-LSVRG} \cite{Ye2020}
  & $\mathcal{O}\left(\max \left\{ {Q,\frac{1}{p}} \right\}\ln \frac{1}{\epsilon}\right)$ with $p = 1/{q_{\max }}$ \\
  \textit{Push-LSVRG-UP} (this work)
  & $\mathcal{O}\left( {\max \{ {\frac{1}{{\left( {1 - \sigma_A } \right)\underline{p} }},\frac{{\delta {\mathcal{Q}^2}\bar p}}{{{{\left( {1 - \sigma_A } \right)}^2}\underline{p} }}} \}\ln \frac{1}{\epsilon}}\right)$ \\
  \bottomrule
\end{tabular}}
\label{Table 2}
\end{table}

\begin{remark}\label{R4-2}
In a big data framework, each agent maintains a large number of data, which leads to the fact that ${q_{\max }} \approx {q_{\min }} \gg \delta {Q^2}{\left( {1 - \sigma_A } \right)^{ - 2}}$ with $\delta  \ge 1$ and $Q \ge 1$. Therefore, \textit{Push-LSVRG-UP} achieves a network-independent computational complexity of $\mathcal{O}\left( {{q_{\max }}\ln\frac{1}{\epsilon}} \right)$ via setting $\underline{p}  = 1/{q_{\max }}{\left( {1 - \sigma_A } \right)}$ and $\bar p = 1/{q_{\min }}{\left( {1 - \sigma_A } \right)}$, which shares the same convergence results with \textit{PMGT-SAGA}/\textit{PMGT-LSVRG} \cite{Ye2020}, \textit{GT-SAGA} \cite{Xin2020f}, and \textit{Push-SAGA} \cite{Qureshi2020a} under the big data framework. This convergence rate is $m$ times faster than the centralized stochastic algorithms \textit{SAGA} \cite{Defazio2014c}, \textit{SVRG} \cite{Johnson2013c}, and \textit{LSVRG} \cite{Kovalev2019}. Furthermore, this improvement under the big data framework may have the potential to be further advanced via setting a pair of different values for triggering probabilities $p_i$, $i \in \mathcal{V}$, according to the requirements in practice. Another improvement is the theoretical results of \textit{Push-LSVRG-UP} available for more generic unbalanced directed networks, while the convergence rates of some elegant distributed stochastic methods, such as \textit{DSA} \cite{Ribeiro2015}, \textit{Diffusion-AVRG} \cite{Yuan2019e}, \textit{GT-SAGA}/\textit{GT-SVRG} \cite{Xin2020f}, \textit{PMGT-SAGA}/\textit{PMGT-LSVRG} \cite{Ye2020}, and \cite{Wang2019} are only available for undirected networks, and moreover the convergence rate of \textit{PMGT-LSVRG} \cite{Ye2020} is based on a coordinated triggering probability, i.e., $p := {p_1} = {p_2} =  \cdots  = {p_m} = 1/{q_i}$, $i \in \left\{ {1, \ldots ,m} \right\}$, under the assumption of ${q_1} = {q_2} =  \cdots  = {q_m}$. In a general data framework, as shown in Table \ref{Table 2}, an explicit linear convergence rate of \textit{Push-LSVRG-UP} for strongly convex objective functions is still meaningful for distributed stochastic optimization over unbalanced directed networks.
\end{remark}

\subsection{Sketch of The Proof}\label{Section 4.2}
To establish the linear convergence of \textit{Push-LSVRG-UP}, we need to analyze the following four error terms:
\begin{enumerate}
  \item the network agreement error: $\left\| {{x_{k }} - {A_\infty }{x_{k }}} \right\|_\pi ^2$,
  \item the convergence error: $m\left\| {{{\bar x}_{k }} - {{\tilde z}^*}} \right\|_2^2$,
  \item the gradient-learning error: ${\Delta _k} := \sum\nolimits_{i = 1}^m {\left( {1/{q_i}} \right)\sum\nolimits_{j = 1}^{{q_i}} {\left\| {\nabla {f_{i,j}}\left( {w_k^i} \right) - \nabla {f_{i,j}}\left( {{{\tilde z}^*}} \right)} \right\|_2^2} } $,
  \item the GT error: $\left\| {{v_{k }} - {A_\infty }{v_{k }}} \right\|_\pi ^2$.
\end{enumerate}
Specifically, we aim to systematically construct a discrete linear time invariant (DLTI) system associated with the above four error terms, see Proposition \ref{P1}. Then, the linear convergence can be obtained through solving for a specific interval of the constant step-size $\alpha$ to guarantee strictly the spectral radius of the system matrix less than $1$.
\subsection{Auxiliary Results}\label{Section 4.3}
The following lemma describes the contraction property of primitive and column-stochastic weight matrix $A$.
\begin{lemma}\label{L1}
Suppose that Assumption \ref{A1} holds. Then $\forall x \in {\mathbb{R}^{mn}}$, we have
\begin{equation}\label{E4-3}
{\left\| {Ax - {A_\infty }x} \right\|_\pi } \le \sigma_A {\left\| {x - {A_\infty }x} \right\|_\pi },
\end{equation}
\begin{proof}
According to compatibility of matrix norms, we have
\begin{equation}\label{E4-3-1}
\begin{aligned}
{\left\| {Ax - {A_\infty }x} \right\|_\pi } =& {\left\| {\left( {A - {A_\infty }} \right)\left( {x - {A_\infty }x} \right)} \right\|_\pi }\\
\le& {\left\| {A - {A_\infty }} \right\|_\pi }{\left\| {x - {A_\infty }x} \right\|_\pi },
\end{aligned}
\end{equation}
where the proof is ended by setting $\sigma_A  = {\left\| {A - {A_\infty }} \right\|_\pi }$.
\end{proof}
\end{lemma}
The next lemma derives upper bounds on the network agreement error.
\begin{lemma}\label{L2}
Suppose that Assumption \ref{A1} holds. Considering the sequence ${\left\{ {{x_k}} \right\}_{k \ge 0}}$ generated by Algorithm \ref{Algo1}, then $\forall k \ge 0$, there hold
\begin{equation}\label{E4-4}
\begin{aligned}
\mathbb{E}\left[ {\left\| {{x_{k + 1}} - {A_\infty }{x_{k + 1}}} \right\|_\pi ^2} \right]\le& \frac{{1 + {\sigma_A ^2}}}{2}\mathbb{E}\left[ {\left\| {{x_k} - {A_\infty }{x_k}} \right\|_\pi ^2} \right]\\
 &+ \frac{{2{\alpha ^2}}}{{1 - {\sigma_A ^2}}}\mathbb{E}\left[ {\left\| {{v_k} - {A_\infty }{v_k}} \right\|_\pi ^2} \right],
\end{aligned}
\end{equation}
and
\begin{equation}\label{E4-5}
\begin{aligned}
\!\!\!\!\mathbb{E}\left[ {\left\| {{x_{k + 1}} - {A_\infty }{x_{k + 1}}} \right\|_\pi ^2} \right]\le& 2\mathbb{E}\left[ {\left\| {{x_k} - {A_\infty }{x_k}} \right\|_\pi ^2} \right]\\
\!\!\!\!& + 2{\alpha ^2}\mathbb{E}\left[ {\left\| {{v_k} - {A_\infty }{v_k}} \right\|_\pi ^2} \right].
\end{aligned}
\end{equation}
\begin{proof}
According to (\ref{E3-1-1}), we have for $ r >0$,
\begin{equation}\label{E4-6}
\begin{aligned}
&\left\| {{x_{k + 1}} - {A_\infty }{x_{k + 1}}} \right\|_\pi ^2\\
=& \left\| {A{x_k} - {A_\infty }{x_k}} \right\|_\pi ^2 + {\alpha ^2}\left\| {{v_k} - {A_\infty }{v_k}} \right\|_\pi ^2\\
&- 2\alpha {\left\langle {A{x_k} - {A_\infty }{x_k},{v_k} - {A_\infty }{v_k}} \right\rangle _\pi }\\
\le& {\sigma_A ^2}\left\| {{x_k} - {A_\infty }{x_k}} \right\|_\pi ^2 + {\alpha ^2}\left\| {{v_k} - {A_\infty }{v_k}} \right\|_\pi ^2\\
&+ 2\sigma_A \alpha {\left\| {{x_k} - {A_\infty }{x_k}} \right\|_\pi }{\left\| {{v_k} - {A_\infty }{v_k}} \right\|_\pi }\\
\le&\left( {1 + r} \right){\sigma_A ^2}\left\| {{x_k} - {A_\infty }{x_k}} \right\|_\pi ^2 \!+\! \left( {1 + {r^{ - 1}}} \right){\alpha ^2}\left\| {{v_k} - {A_\infty }{v_k}} \right\|_\pi ^2,
\end{aligned}
\end{equation}
where the first inequality follows Lemma \ref{L1} and the last inequality employs the Young’s inequality. The proof is completed by setting $ r = \left( {1 - {\sigma_A ^2}} \right)/2{\sigma_A ^2}$ and $r = 1$, respectively.
\end{proof}
\end{lemma}
The next lemma captures the contraction property of $Y_k$. This is not necessary in the convergence analysis of \textit{DSA} \cite{Ribeiro2015}, \textit{Diffusion-AVRG} \cite{Yuan2019e}, \textit{GT-SAGA}/\textit{GT-SVRG} \cite{Xin2020f}, \textit{PMGT-SAGA}/\textit{PMGT-LSVRG} \cite{Ye2020}, and \cite{Wang2019} since these algorithm are only available to undirected networks and not involved in imbalanced information exchange.
\begin{lemma}\label{L3}
Suppose that Assumption \ref{A1} holds. Recalling the definition of ${Y_k}$ and denoting ${Y_\infty } = {\lim _{k \to \infty }}{Y_k}$, we have $\forall k \ge 0$
\begin{equation}\label{E4-7}
{\left\| {{Y_k} - {Y_\infty }} \right\|_2} \le T{\sigma_A ^k},
\end{equation}
where $T := \sqrt \vartheta  {\left\| {{1_m} - m\pi } \right\|_2}$ with $\vartheta  := \bar \pi /\underline{\pi }  > 1$.
\begin{proof}
To explore the upper bound of ${\left\| {{Y_k} - {Y_\infty }} \right\|_2}$, we need to define ${y_\infty } := {{\underline{A} }_\infty }{y_k}$ according to (\ref{E3-1-2}). Then, $\forall k \ge 0$, there holds
\begin{equation}\label{E4-8}
\begin{aligned}
{\left\| {{Y_k} - {Y_\infty }} \right\|_2}=& {\left\| {{\text{diag}}\left\{ {{y_k}} \right\} \otimes {I_n} - {\text{diag}}\left\{ {{y_\infty }} \right\} \otimes {I_n}} \right\|_2}\\
\le & {\left\| {{y_k} - {y_\infty }} \right\|_2}\\
\le& {{\bar \pi }^{0.5}}\sigma_A {\left\| {{y_{k - 1}} - {y_\infty }} \right\|_\pi }\\
\le& {\left( {\frac{{\bar \pi }}{{\underline{\pi } }}} \right)^{0.5}}{\left\| {{1_m} - m\pi } \right\|_2}{\sigma_A ^k},
\end{aligned}
\end{equation}
where the first inequality applies the definitions of standard 2-norm for vectors and matrices, and the last inequality is owing to the fact that ${y_0} = {1_m}$ and ${y_\infty } = m\pi $. The proof is ended by setting $\vartheta  = \bar \pi /\underline{\pi }  > 1$ and $T = \sqrt \vartheta  {\left\| {{1_m} - m\pi } \right\|_2}$.
\end{proof}
\end{lemma}
The next lemma derives the upper bound on the convergence error. We define a directivity constant $\delta  := Y\left( {1 + T} \right)\vartheta {\tilde Y^2}$ and denote ${{\mathcal{F}_k}}$ as the filter of the history of the dynamical system yielded by the sequence $\left\{ {s_k^i} \right\}_{k \ge 0}^{i = 1, \ldots ,m}$.
\begin{lemma}\label{L4}
Suppose that Assumptions \ref{A1}-\ref{A3} hold. Recalling the definition of ${{\bar x}_{k}} $, if the step-size satisfies $0 < \alpha  \le \min \left\{ {1/{L},{{m\mu }}/{{8{L^2}}}} \right\}$, then $\forall k \ge 0$, we have
\begin{equation}\label{E4-8-1}
\begin{aligned}
&\mathbb{E}\left[ {m\left\| {{{\bar x}_{k + 1}} - {{\tilde z}^*}} \right\|_2^2} \right]\\
 \le& \left( {1 - \frac{{\mu \alpha }}{2}} \right)\mathbb{E}\left[ {m\left\| {{{\bar x}_k} \! - \!{{\tilde z}^*}} \right\|_2^2} \right] \!+\! \frac{{2\delta \bar \pi \alpha {L^2}}}{\mu }\mathbb{E}\left[ {\left\| {{x_k} - {A_\infty }{x_k}} \right\|_\pi ^2} \right]\\
&+\frac{{2{\alpha ^2}}}{m}\mathbb{E}\left[ {{\Delta _k}} \right] + \frac{{2\delta T \alpha {L^2}{\sigma_A ^k}}}{\mu }\mathbb{E}\left[ {\left\| {{x_k}} \right\|_2^2} \right],
\end{aligned}
\end{equation}
and
\begin{equation}\label{E4-8-2}
\begin{aligned}
&\mathbb{E}\left[ {m\left\| {{{\bar x}_{k + 1}} - {{\tilde z}^*}} \right\|_2^2} \right]\\
 \le& 3\mathbb{E}\left[ {m\left\| {{{\bar x}_k} - {{\tilde z}^*}} \right\|_2^2} \right] + 6\delta \bar \pi {L^2}{\alpha ^2}\mathbb{E}\left[ {\left\| {{x_k} - {A_\infty }{x_k}} \right\|_\pi ^2} \right]\\
&+ \frac{{2{\alpha ^2}}}{m}\mathbb{E}\left[ {{\Delta _k}} \right] + 6\delta T{L^2}{\alpha ^2}{\sigma_A ^k}\mathbb{E}\left[ {\left\| {{x_k}} \right\|_2^2} \right].
\end{aligned}
\end{equation}
\begin{proof}
Recall the definition of ${{{\bar x}_k}}$, considering (\ref{E3-1-1}), we have
\begin{equation}\label{E4-9}
\begin{aligned}
&\mathbb{E}\left[ {\left\| {{{\bar x}_{k + 1}} - {{\tilde z}^*}} \right\|_2^2} |{\mathcal{F}_k}\right]\\
=& \mathbb{E}\left[ {\left\| {{{\bar x}_k} - \alpha {{\bar g}_k} - {{\tilde z}^*}} \right\|_2^2} \right]\\
=& \left\| {{{\bar x}_k} - {{\tilde z}^*}} \right\|_2^2 - 2\alpha \mathbb{E}\left[ {\left\langle {{{\bar x}_k} - {{\tilde z}^*},{{\bar g}_k}} \right\rangle |{\mathcal{F}_k}} \right]+ {\alpha^2}\mathbb{E}\left[ {\left\| {{{\bar g}_k}} \right\|_2^2|{\mathcal{F}_k}} \right]\\
=&\left\| {{{\bar x}_k} - {{\tilde z}^*}} \right\|_2^2 - 2\alpha \left\langle {{{\bar x}_k} - {{\tilde z}^*},{{\tilde p}_k}} \right\rangle  + 2\alpha \left\langle {{{\bar x}_k} - {{\tilde z}^*},{{\tilde p}_k} - {{\bar h}_k}} \right\rangle\\
&+ {\alpha ^2}\mathbb{E}\left[ {\left\| {{{\bar g}_k} - {{\bar h}_k}} \right\|_2^2|{\mathcal{F}_k}} \right] + {\alpha ^2}\left\| {{{\bar h}_k} - {{\tilde p}_k} + {{\tilde p}_k}} \right\|_2^2\\
=& \left\| {{{\bar x}_k} - \alpha {{\tilde p}_k} - {{\tilde z}^*}} \right\|_2^2 + {\alpha ^2}\left\| {{{\tilde p}_k} - {{\bar h}_k}} \right\|_2^2\\
&+ 2\alpha \left\langle {{{\bar x}_k} - \alpha {{\tilde p}_k} - {{\tilde z}^*},{{\tilde p}_k} - {{\bar h}_k}} \right\rangle  + {\alpha ^2}\mathbb{E}\left[ {\left\| {{{\bar g}_k} - {{\bar h}_k}} \right\|_2^2|{\mathcal{F}_k}} \right]\\
\le & {\left( {1 \!-\! \alpha \mu } \right)^2}\left\| {{{\bar x}_k} \!-\! {{\tilde z}^*}} \right\|_2^2 \!+\! {\alpha ^2}\left\| {{{\tilde p}_k} \!-\! {{\bar h}_k}} \right\|_2^2 \!+\! {\alpha ^2}\mathbb{E}[ {\| {{{\bar g}_k} \!-\! {{\bar h}_k}} \|_2^2|{\mathcal{F}_k}} ]\\
&+ \alpha \left( {1 - \alpha \mu } \right)\left( {\mu \left\| {{{\bar x}_k} - {{\tilde z}^*}} \right\|_2^2 + \frac{1}{\mu }\left\| {{{\tilde p}_k} - {{\bar h}_k}} \right\|_2^2} \right)\\
=& \left( {1 - \alpha \mu } \right)\left\| {{{\bar x}_k} - {{\tilde z}^*}} \right\|_2^2 \!+\! \frac{\alpha }{{m\mu }}\left\| {\nabla F\left( {{1_m} \otimes {{\bar x}_k}} \right) - \nabla F\left( {{z_k}} \right)} \right\|_2^2\\
& + {\alpha ^2}\mathbb{E}\left[ {\left\| {{{\bar g}_k} - {{\bar h}_k}} \right\|_2^2|{\mathcal{F}_k}} \right]\\
\le& \left( {1 - \alpha \mu } \right)\left\| {{{\bar x}_k} - {{\tilde z}^*}} \right\|_2^2 + \left( {\frac{{\alpha {L^2}}}{{m\mu }}} \right)\left\| {{z_k} - {1_m} \otimes {{\bar x}_k}} \right\|_2^2\\
& + {\alpha ^2}\mathbb{E}\left[ {\left\| {{{\bar g}_k} - {{\bar h}_k}} \right\|_2^2|{\mathcal{F}_k}} \right],
\end{aligned}
\end{equation}
where the first inequality employs both a standard contraction theorem \cite[Lemma 1]{Xin2020f} in convex optimization provided that the step-size satisfies $0 < \alpha  \le 1/L$ for the first term and the Young's Inequality for the last term. The last inequality uses the $L$-smoothness, i.e., the Lipschitz continuity, of the local objective functions.
We proceed to handle $\mathbb{E}\left[ {\left\| {{{\bar g}_k} - {{\bar h}_k}} \right\|_2^2|{\mathcal{F}_k}} \right]$ as follows:
\begin{equation}\label{E4-10}
\begin{aligned}
\mathbb{E}\left[ {\left\| {{{\bar g}_k} - {{\bar h}_k}} \right\|_2^2|{\mathcal{F}_k}} \right] =& \frac{1}{{{m^2}}}\mathbb{E}\left[ {\left\| {\sum\limits_{i = 1}^m {\left( {g_k^i - \nabla {f_i}\left( {z_k^i} \right)} \right)} } \right\|_2^2|{\mathcal{F}_k}} \right]\\
=& \frac{1}{{{m^2}}}\mathbb{E}\left[ {\sum\limits_{i = 1}^n {\left\| {g_k^i - \nabla {f_i}\left( {z_k^i} \right)} \right\|_2^2} |{\mathcal{F}_k}} \right] \\
=& \frac{1}{{{m^2}}}\mathbb{E}\left[ {\left\| {{g_k} - \nabla F\left( {{z_k}} \right)} \right\|_2^2|{\mathcal{F}_k}} \right],
\end{aligned}
\end{equation}
where the first equality is obtained from the definitions of ${{\bar g}_k}$ and ${{\bar h}_k}$, and the second equality is owing to the fact that $\mathbb{E}\left[ {\sum\limits_{i \ne j} {\left\langle {g_k^i - \nabla {f_i}\left( {z_k^i} \right),g_k^j - \nabla {f_j}\left( {z_k^j} \right)} \right\rangle } |{\mathcal{F}_k}} \right]=0$. In the next step, we continue to seek the upper bound on the expectation of variance term $\mathbb{E}\left[ {\left\| {{g_k} - \nabla F\left( {{z_k}} \right)} \right\|_2^2|{\mathcal{F}_k}} \right]$. Before deriving the upper bound, we first define $a := \nabla {f_{i,s_k^i}}\left( {z_k^i} \right) - \nabla {f_{i,s_k^i}}\left( {{{\tilde z}^*}} \right)$ and $b := \nabla {f_{i,s_k^i}}\left( {w_k^i} \right) - \nabla {f_{i,s_k^i}}\left( {{{\tilde z}^*}} \right)$, and then recall the definition of $g_k^i$ such that
\begin{equation}\label{E4-11}
\begin{aligned}
&\mathbb{E}\left[ {\left\| {g_k^i - \nabla {f_i}\left( {z_k^i} \right)} \right\|_2^2|{\mathcal{F}_k}} \right]\\
=& \mathbb{E}\left[ {\left\| {\nabla\! {f_{i,s_k^i}}\left( {z_k^i} \right) \!-\! \nabla {f_{i,s_k^i}}\left( {w_k^i} \right) \!+\! \nabla {f_i}\left( {w_k^i} \right) \!-\! \nabla\! {f_i}\left( {z_k^i} \right)} \! \right\|_2^2|{\mathcal{F}_k}} \right]\\
=&\mathbb{E}\left[ {\left\| {a - \mathbb{E}\left[ {a|{\mathcal{F}_k}} \right] - \left( {b - \mathbb{E}\left[ {b|{\mathcal{F}_k}} \right]} \right)} \right\|_2^2|{\mathcal{F}_k}} \right]\\
\le& 2\mathbb{E}\left[ {\left\| {a - \mathbb{E}\left[ {a|{\mathcal{F}_k}} \right]} \right\|_2^2|{\mathcal{F}_k}} \right] + 2\mathbb{E}\left[ {\left\| {b - \mathbb{E}\left[ {b|{\mathcal{F}_k}} \right]} \right\|_2^2|{\mathcal{F}_k}} \right]\\
=& 2\mathbb{E}\left[ {\left\| {\nabla {f_{i,s_k^i}}\left( {z_k^i} \right)- \nabla {f_{i,s_k^i}}\left( {{{\tilde z}^*}} \right)} \right\|_2^2|{\mathcal{F}_k}} \right]\\
&+ 2\mathbb{E}\left[ {\left\| {\nabla {f_{i,s_k^i}}\left( {w_k^i} \right) - \nabla {f_{i,s_k^i}}\left( {{{\tilde z}^*}} \right)} \right\|_2^2|{\mathcal{F}_k}} \right]\\
& - 2\left\| {\nabla {f_i}\left( {w_k^i} \right) - \nabla {f_i}\left( {{{\tilde z}^*}} \right)} \right\|_2^2 - 2\left\| {\nabla {f_i}\left( {z_k^i} \right) - \nabla {f_i}\left( {{{\tilde z}^*}} \right)} \right\|_2^2\\
\le& 2\mathbb{E}\left[ {\left\| {\nabla {f_{i,s_k^i}}\left( {z_k^i} \right) - \nabla {f_{i,s_k^i}}\left( {{{\tilde z}^*}} \right)} \right\|_2^2|{\mathcal{F}_k}} \right]\\
&+ 2\mathbb{E}\left[ {\left\| {\nabla {f_{i,s_k^i}}\left( {w_k^i} \right) - \nabla {f_{i,s_k^i}}\left( {{{\tilde z}^*}} \right)} \right\|_2^2|{\mathcal{F}_k}} \right]\\
=&\frac{2}{{{q_i}}}\sum\limits_{j = 1}^{{q_i}} {\left\| {\nabla {f_{i,j}}\left( {z_k^i} \right) - \nabla {f_{i,j}}\left( {{{\tilde z}^*}} \right)} \right\|_2^2}\\
&+ \underbrace {\frac{2}{{{q_i}}}\sum\limits_{j = 1}^{{q_i}} {\left\| {\nabla {f_{i,j}}\left( {w_k^i} \right) - \nabla {f_{i,j}}\left( {{{\tilde z}^*}} \right)} \right\|_2^2} }_{:= 2\Delta _k^i}\\
\le& 4{L^2}\left( {\left\| {z_k^i - {{\bar x}_k}} \right\|_2^2 + \left\| {{{\bar x}_k} - {{\tilde z}^*}} \right\|_2^2} \right)+2\Delta _k^i,
\end{aligned}
\end{equation}
where the last inequality uses the $L$-smoothness of $\nabla f_{i,j}$. Summing (\ref{E4-11}) over $i$ and taking the total expectation yield
\begin{equation}\label{E4-12}
\begin{aligned}
&\mathbb{E}\left[ {\left\| {{g_k} - \nabla F\left( {{z_k}} \right)} \right\|_2^2|{\mathcal{F}_k}} \right]\\
\le& 4{L^2}\left( {\left\| {{z_k} - {1_m} \otimes {{\bar x}_k}} \right\|_2^2 + m\left\| {{{\bar x}_k} - {{\tilde z}^*}} \right\|_2^2} \right)+2{\Delta _k}.
\end{aligned}
\end{equation}
Plugging (\ref{E4-12}) into (\ref{E4-10}) obtains
\begin{equation}\label{E4-13}
\begin{aligned}
&\mathbb{E}\left[ {\left\| {{{\bar g}_k} - {{\bar h}_k}} \right\|_2^2|{\mathcal{F}_k}} \right]\\
 \le& \frac{{4{L^2}}}{{{m^2}}}\left\| {{z_k} - {1_m} \otimes {{\bar x}_k}} \right\|_2^2 + \frac{{4{L^2}}}{m}\left\| {{{\bar x}_k} - {{\tilde z}^*}} \right\|_2^2+\frac{2}{{{m^2}}}{\Delta _k}.
 \end{aligned}
\end{equation}
Combining (\ref{E4-9}) with (\ref{E4-13}) reduces to
\begin{equation}\label{E4-14}
\begin{aligned}
\mathbb{E}\left[ {\left\| {{{\bar x}_{k +\! 1}} \!-\! {{\tilde z}^*}} \right\|_2^2}|{\mathcal{F}_k} \right]\!\le&( {1 \!-\! \alpha \mu  \!+\! \frac{{4{L^2}{\alpha ^2}}}{m}} )\left\| {{{\bar x}_k} \!-\! {{\tilde z}^*}} \right\|_2^2 \!+\! \frac{{2{\alpha ^2}}}{{{m^2}}}{\Delta _k}\\
&+ \!( {\frac{\alpha {{L^2}}}{{m\mu }} +\! \frac{{4{L^2}{\alpha ^2}}}{{{m^2}}}} )\left\| {{z_k} - {1_m} \otimes {{\bar x}_k}} \right\|_2^2.
\end{aligned}
\end{equation}
We next handle $\left\| {{z_k} - {1_m} \otimes {{\bar x}_k}} \right\|_2^2$ as follows:
\begin{equation}\label{E4-15}
\begin{aligned}
&\left\| {{z_k} - {1_m} \otimes {{\bar x}_k}} \right\|_2^2\\
=& \!\left\| {Y_k^{ - 1}\left( {{x_k} \!-\! {Y_\infty }\left( {{1_m} \otimes {{\bar x}_k}} \right)} \right) \!+\! \left( {Y_k^{ \!- \!1}{Y_\infty } \!-\! {I_{mn}}} \right)\left( {{1_m} \otimes {{\bar x}_k}} \right)} \right\|_2^2\\
\le& \!\left\| {Y_k^{ \!-\! 1}( {{x_k} \!-\! {Y_\infty }\left( {{1_m} \!\otimes\! {{\bar x}_k}} \right)} )} \right\|_2^2 \! + \!\left\| {( {Y_k^{ \!-\! 1}{Y_\infty } \!-\! {I_{mn}}} )( {{1_m} \!\otimes\! {{\bar x}_k}} )} \right\|_2^2\\
&+\! 2{\left\| {Y_k^{ -\! 1}\left( {{x_k} \!-\! {Y_\infty }\left( {{1_m} \!\otimes\! {{\bar x}_k}} \right)} \right)\left( {Y_k^{ -\! 1}{Y_\infty } \!-\! {I_{mn}}} \right)( {{1_m} \!\otimes\! {{\bar x}_k}} )} \right\|_2}\\
\le& {{\tilde Y}^2}\left\| {{x_k} - {A_\infty }{x_k}} \right\|_2^2 + \left\| {\left( {Y_k^{ - 1}{Y_\infty } - {I_{mn}}} \right)\left( {{1_m} \otimes {{\bar x}_k}} \right)} \right\|_2^2\\
&+ 2\tilde Y{\left\| {{x_k} - {A_\infty }{x_k}} \right\|_2}{\left\| {\left( {Y_k^{ - 1}{Y_\infty } - {I_{mn}}} \right)\left( {{1_m} \otimes {{\bar x}_k}} \right)} \right\|_2}\\
\le& {{\tilde Y}^2}\left\| {{x_k} - {A_\infty }{x_k}} \right\|_2^2 + 2T{\sigma_A ^k}{{\tilde Y}^2}{\left\| {{x_k} - {A_\infty }{x_k}} \right\|_2}{\left\| {{x_k}} \right\|_2}\\
&+{\left( {\tilde YT{\sigma_A ^k}} \right)^2}\left\| {{x_k}} \right\|_2^2\\
\le&\bar \pi \left( {1 + T} \right){{\tilde Y}^2}\left\| {{x_k} - {A_\infty }{x_k}} \right\|_\pi ^2 + T\left( {T + 1} \right){{\tilde Y}^2}{\sigma_A ^k}\left\| {{x_k}} \right\|_2^2,
\end{aligned}
\end{equation}
where the third inequality applies Lemma \ref{L3} and the last inequality uses the fact that ${0 < \sigma_A  < 1}$. Via defining ${d_1} :=  \left( {1 + T} \right){\tilde Y^2}$ and ${d_2} :=  T\left( {T + 1} \right){\tilde Y^2}$, one can attain
\begin{equation}\label{E4-16}
\!\left\| {{z_k} - {1_m} \otimes {{\bar x}_k}} \right\|_2^2 \le \bar \pi {d_1}\left\| {{x_k} - {A_\infty }{x_k}} \right\|_\pi ^2 + {d_2}{\sigma_A ^k}\left\| {{x_k}} \right\|_2^2.
\end{equation}
Recalling the definition of $\delta$, we know $d_1\le \delta$ and $d_2\le \delta T$. Plugging (\ref{E4-16}) into (\ref{E4-14}) yields
\begin{equation}\label{E4-17}
\begin{aligned}
&\mathbb{E}\left[ {\left\| {{{\bar x}_{k + 1}} - {{\tilde z}^*}} \right\|_2^2}|{\mathcal{F}_k} \right]\\
\le& \left( {1 - \alpha \mu  + \frac{{4{L^2}{\alpha ^2}}}{m}} \right)\left\| {{{\bar x}_k} - {{\tilde z}^*}} \right\|_2^2 + \left( {\frac{{\alpha {L^2}}}{{m\mu }} + \frac{{4{L^2}{\alpha ^2}}}{{{m^2}}}} \right){d_1}\bar \pi\\
&\times \!\left\| {{x_k} \!-\! {A_\infty }{x_k}} \right\|_\pi ^2 \!+\! \frac{{2{\alpha ^2}}}{{{m^2}}}{\Delta _k} \!+\! \left( {\frac{{\alpha {L^2}}}{{m\mu }} \!+\! \frac{{4{L^2}{\alpha ^2}}}{{{m^2}}}} \right){d_2}{\sigma_A ^k}\| {{x_k}} \|_2^2.
\end{aligned}
\end{equation}
If one picks $0 < \alpha  \le m/\left( {4\mu } \right)$ on the final term, then there holds
\begin{equation}\label{E4-18}
\begin{aligned}
&\mathbb{E}\left[ {\left\| {{{\bar x}_{k + 1}} - {{\tilde z}^*}} \right\|_2^2}|{\mathcal{F}_k} \right]\\
\le&\left( {1 \!-\! \alpha \mu  \!+\! \frac{{4{L^2}{\alpha ^2}}}{m}} \right)\left\| {{{\bar x}_k} \!-\! {{\tilde z}^*}} \right\|_2^2 + \frac{{2\bar \pi {d_1}\alpha {L^2}}}{{m\mu }}\left\| {{x_k} \!-\! {A_\infty }{x_k}} \right\|_\pi ^2\\
&+\frac{{2{L^2}{\alpha ^2}}}{{{m^2}}}{\Delta _k} + \frac{{2{d_2}\alpha {L^2}{\sigma_A ^k}}}{{m\mu }}\left\| {{x_k}} \right\|_2^2.
\end{aligned}
\end{equation}
If one further chooses $0 < \alpha  \le m\mu /\left( {8{L^2}} \right)$, it holds
\begin{equation}\label{E4-19}
\begin{aligned}
&\mathbb{E}\left[ {\left\| {{{\bar x}_{k + 1}} - {{\tilde z}^*}} \right\|_2^2}|{\mathcal{F}_k} \right]\\
\le& \left( {1 - \frac{{\mu \alpha }}{2}} \right)\left\| {{{\bar x}_k} - {{\tilde z}^*}} \right\|_2^2 + \frac{{2\bar \pi {d_1}\alpha {L^2}}}{{m\mu }}\left\| {{x_k} - {A_\infty }{x_k}} \right\|_\pi ^2\\
&+\frac{{2{\alpha ^2}}}{{{m^2}}}{\Delta _k} + \frac{{2{d_2}\alpha {L^2}{\sigma_A ^k}}}{{m\mu }}\left\| {{x_k}} \right\|_2^2.
\end{aligned}
\end{equation}
Recall that $Y \ge 1$ and $\tilde Y \ge 1$. Then, we have
\begin{equation}\label{E4-20}
\begin{aligned}
&\mathbb{E}\left[ {\left\| {{{\bar x}_{k + 1}} - {{\tilde z}^*}} \right\|_2^2}|{\mathcal{F}_k} \right]\\
\le& \left( {1 - \frac{{\mu \alpha }}{2}} \right)\left\| {{{\bar x}_k} - {{\tilde z}^*}} \right\|_2^2 + \frac{{2\delta \bar \pi \alpha {L^2}}}{{m\mu }}\left\| {{x_k} - {A_\infty }{x_k}} \right\|_\pi ^2\\
&+\frac{{2{\alpha ^2}}}{{{m^2}}}{\Delta _k} + \frac{{2\delta T\alpha {L^2}{\sigma_A ^k}}}{{m\mu }}\left\| {{x_k}} \right\|_2^2.
\end{aligned}
\end{equation}
Multiplying $m$ to the both sides of (\ref{E4-20}) and taking the total expectation on the both sides of the above inequality yields (\ref{E4-8-1}).
In the another technical line, if we modify the first inequality in (\ref{E4-9}) as follows:
\begin{equation}\label{E4-22}
\begin{aligned}
&\mathbb{E}\left[ {\left\| {{{\bar x}_{k + 1}} - {{\tilde z}^*}} \right\|_2^2}|{\mathcal{F}_k} \right]\\
\le& {\left( {1 \!-\! \alpha \mu } \right)^2}\left\| {{{\bar x}_k} \!-\! {{\tilde z}^*}} \right\|_2^2 \!+ {\alpha ^2}\left\| {{{\tilde p}_k} \!-\! {{\bar h}_k}} \right\|_2^2 \!+ {\alpha ^2}\mathbb{E}[ {\| {{{\bar g}_k} \!-\! {{\bar h}_k}} \|_2^2|{\mathcal{F}_k}} ]\\
&+ \left( {1 - \alpha \mu } \right)\left( {\left\| {{{\bar x}_k} - {{\tilde z}^*}} \right\|_2^2 + {\alpha ^2}\left\| {{{\tilde p}_k} - {{\bar h}_k}} \right\|_2^2} \right)\\
\le&2\left\| {{{\bar x}_k} \!-\! {{\tilde z}^*}} \right\|_2^2 \!+\! \frac{{2{\alpha ^2}{L^2}}}{m}\left\| {{z_k} \!-\! {1_m} \!\otimes\! {{\bar x}_k}} \right\|_2^2 \!+\! {\alpha ^2}\mathbb{E}[ {\| {{{\bar g}_k} \!-\! {{\bar h}_k}} \|_2^2|{\mathcal{F}_k}} ]\\
\le& \left( {2 + \frac{{4{L^2}{\alpha ^2}}}{m}} \right)\left\| {{{\bar x}_k} - {{\tilde z}^*}} \right\|_2^2 + \frac{{6{L^2}{\alpha ^2}}}{m}\left\| {{z_k} - {1_m} \otimes {{\bar x}_k}} \right\|_2^2\\
& + \frac{{2{\alpha ^2}}}{{{m^2}}}{\Delta _k}\\
\le& \left( {2 + \frac{{4{L^2}{\alpha ^2}}}{m}} \right)\left\| {{{\bar x}_k} - {{\tilde z}^*}} \right\|_2^2 + \frac{{6\bar \pi {d_1}{L^2}{\alpha ^2}}}{m}\left\| {{x_k} - {A_\infty }{x_k}} \right\|_\pi ^2\\
&+ \frac{{2{\alpha ^2}}}{{{m^2}}}{\Delta _k} + \frac{{6{d_2}{L^2}{\alpha ^2}}}{m}{\sigma_A ^k}\left\| {{x_k}} \right\|_2^2,
\end{aligned}
\end{equation}
where the third inequality uses (\ref{E4-13}) and the last inequality applies (\ref{E4-16}), then
via choosing the step-size ${0 < \alpha  \le \sqrt m /\left( {2L} \right)}$, the above inequality further becomes
\begin{equation}\label{E4-23}
\begin{aligned}
&\mathbb{E}\left[ {\left\| {{{\bar x}_{k + 1}} - {{\tilde z}^*}} \right\|_2^2}|{\mathcal{F}_k} \right]\\
\!\!\!\!\le&3\left\| {{{\bar x}_k} - {{\tilde z}^*}} \right\|_2^2 + \frac{{6\bar \pi {d_1}{L^2}{\alpha ^2}}}{m}\left\| {{x_k} - {A_\infty }{x_k}} \right\|_\pi ^2+ \frac{{2{\alpha ^2}}}{{{m^2}}}{\Delta _k}\\
\!\!\!\!&+ \frac{{6{d_2}{L^2}{\alpha ^2}}}{m}{\sigma_A ^k}\left\| {{x_k}} \right\|_2^2.
\end{aligned}
\end{equation}
Finally, multiplying $m$ to the both sides of (\ref{E4-23}) and taking the total expectation yield another inequality (\ref{E4-8-2}).
\end{proof}
\end{lemma}
In the next lemma, we seek an upper bound on the gradient-learning error of \textit{Push-LSVRG-UP}, which is an important result for distributed stochastic optimization over unbalanced directed netowrks distinguished from the existing distributed stochastic gradient work, such as \textit{DSA} \cite{Ribeiro2015}, \textit{GT-SAGA}/\textit{GT-SVRG} \cite{Xin2020f}, \textit{DSGT} \cite{Pu2021a}, \textit{S-ADDOPT} \cite{Qureshi2021}, \textit{Push-SAGA} \cite{Qureshi2020a}, \textit{PMGT-SAGA}/\textit{PMGT-LSVRG} \cite{Ye2020}, and \cite{Wang2019}.
\begin{lemma}\label{L5}
Suppose that Assumptions \ref{A1} and \ref{A3} hold. Recalling the definition of ${\Delta _k}$, we have
\begin{equation}\label{E4-25}
\begin{aligned}
&\mathbb{E}\left[ {{\Delta _{k + 1}}} \right]\\
\le& \left( {1 - \underline{p} } \right)\mathbb{E}\left[ {{\Delta _k}} \right] + 2\bar p\bar \pi {d_1}{L^2}\mathbb{E}\left[ {\left\| {{x_k} - {A_\infty }{x_k}} \right\|_\pi ^2} \right]\\
&+ 2m\bar p{L^2}\mathbb{E}\left[ {\left\| {{{\bar x}_k} - {{\tilde z}^*}} \right\|_2^2} \right] + 2\bar p{d_2}{L^2}{\sigma_A ^k}\mathbb{E}\left[ {\left\| {{x_k}} \right\|_2^2} \right].
\end{aligned}
\end{equation}
\begin{proof}
Recalling the definition of ${\Delta _{k + 1}^i}$, we have
\begin{equation*}
\begin{aligned}
&\mathbb{E}\left[ {\Delta _{k + 1}^i|{\mathcal{F}_k}} \right]\\
=& \frac{1}{{{q_i}}}\mathbb{E}\left[ {\sum\limits_{j = 1}^{{q_i}} {\left\| {\nabla {f_{i,j}}\left( {w_{k + 1}^i} \right) - \nabla {f_{i,j}}\left( {{{\tilde z}^*}} \right)} \right\|_2^2|{\mathcal{F}_k}} } \right]\\
=& \frac{{1 - {p_i}}}{{{q_i}}}\sum\limits_{j = 1}^{{q_i}} {\left\| {\nabla {f_{i,j}}\left( {w_k^i} \right) - \nabla {f_{i,j}}\left( {{{\tilde z}^*}} \right)} \right\|_2^2}   \\
&+\frac{{{p_i}}}{{{q_i}}}\sum\limits_{j = 1}^{{q_i}} {\left\| {\nabla {f_{i,j}}\left( {z_k^i} \right) - \nabla {f_{i,j}}\left( {{{\tilde z}^*}} \right)} \right\|_2^2}\\
\end{aligned}
\end{equation*}
\begin{equation}\label{E4-26}
\begin{aligned}
=& \left( {1 - {p_i}} \right)\Delta _k^i + \frac{{{p_i}}}{{{q_i}}}\sum\limits_{j = 1}^{{q_i}} {\left\| {\nabla {f_{i,j}}\left( {z_k^i} \right) - \nabla {f_{i,j}}\left( {{{\tilde z}^*}} \right)} \right\|_2^2}  \\
\le& \left( {1 - \underline{p} } \right)\Delta _k^i + \bar p{L^2}\left\| {z_k^i - {{\tilde z}^*}} \right\|_2^2,
\end{aligned}
\end{equation}
where the last inequality follows the $L$-smoothness of $\nabla f_{i, j}$.  Summing (\ref{E4-26}) over $i$ yields
\begin{equation*}
\begin{aligned}
&\mathbb{E}\left[ {{\Delta _{k + 1}}|{\mathcal{F}_k}} \right]\\
\le & \left( {1 - \underline{p} } \right)\sum\limits_{i = 1}^m {\Delta _k^i}  + \bar p{L^2}\sum\limits_{i = 1}^m {\left\| {z_k^i - {{\tilde z}^*}} \right\|_2^2} \\
\le & \left( {1 - \underline{p} } \right){\Delta _k} + 2\bar p{L^2}\left\| {{z_k} - {1_m} \otimes {{\bar x}_k}} \right\|_2^2 + 2m\bar p{L^2}\left\| {{{\bar x}_k} - {{\tilde z}^*}} \right\|_2^2\\
\end{aligned}
\end{equation*}
\begin{equation}\label{E4-27}
\begin{aligned}
\le & \left( {1 - \underline{p} } \right){\Delta _k} \!+\! 2\bar p\bar \pi {d_1}{L^2}\left\| {{x_k} \!-\! {A_\infty }{x_k}} \right\|_\pi ^2 \!+\! 2m\bar p{L^2}\left\| {{{\bar x}_k} - {{\tilde z}^*}} \right\|_2^2\\
&+2\bar p{d_2}{L^2}{\sigma_A ^k}\left\| {{x_k}} \right\|_2^2,
\end{aligned}
\end{equation}
where the last inequality uses the result (\ref{E4-16}).
The proof is completed by taking the total expectation on the both sides of (\ref{E4-27}).
\end{proof}
\end{lemma}
The upper bound on the GT error is sought in the next lemma.
\begin{lemma}\label{L6}
Suppose that Assumptions \ref{A1} and \ref{A3} hold. Considering the sequence ${\left\{ {{v_k}} \right\}_{k \ge 0}}$ generated by Algorithm \ref{Algo1}, if the step-size satisfies $0 < \alpha  \le \min \left\{ {\frac{1}{{2L\sqrt 6 }},\frac{1}{{2LY\sqrt {3\left( {3{{\tilde Y}^2} + 16{d_2}} \right)} }},\frac{{1 - {\sigma_A ^2}}}{{2L\sqrt {2\vartheta \left( {9{{\tilde Y}^2} + 16{d_1} + 48{d_2}} \right)} }}} \right\}$, then $\forall k \ge 0$, we have
\begin{equation}\label{E4-29}
\begin{aligned}
&\mathbb{E}\left[ {\left\| {{v_{k + 1}} - {A_\infty }{v_{k + 1}}} \right\|_\pi ^2} \right]\\
\le& \left( {\frac{{3 \!+\!{\sigma_A ^2}}}{4}} \right)\mathbb{E}\left[ {\left\| {{v_k} \!-\! {A_\infty }{v_k}} \right\|_\pi ^2} \right]\!+\!\frac{{194\delta {L^2}}}{{1 \!-\! {\sigma_A ^2}}}\mathbb{E}\left[ {\left\| {{x_k} \!-\! {A_\infty }{x_k}} \right\|_\pi ^2} \right]\\
&+\frac{{169{{\underline{\pi } }^{ - 1}}{L^2}}}{{1 - {\sigma_A ^2}}}\mathbb{E}\left[ {m\left\| {{{\bar x}_k} - {{\tilde z}^*}} \right\|_2^2} \right] + \frac{{110{{\underline{\pi } }^{ - 1}}}}{{3\left( {1 - {\sigma_A ^2}} \right)}}\mathbb{E}\left[ {{\Delta _k}} \right]\\
&+\frac{{194{{\underline{\pi } }^{ - 1}}T{\delta ^2}{L^2}{\sigma_A ^k}}}{{1 - {\sigma_A ^2}}}\mathbb{E}\left[ {\left\| {{x_k}} \right\|_2^2} \right].
\end{aligned}
\end{equation}
\begin{proof}
According to (\ref{E3-1-4}), we have
\begin{equation}\label{E4-30}
\begin{aligned}
&\mathbb{E}\left[ {\left\| {{v_{k + 1}} - {A_\infty }{v_{k + 1}}} \right\|_\pi ^2|{\mathcal{F}_k}} \right]\\
=& \mathbb{E}\left[ {\left\| {A{v_k} - {A_\infty }{v_k} + \left( {{I_{mn}} - {A_\infty }} \right)\left( {{g_{k + 1}} - {g_k}} \right)} \right\|_\pi ^2|{\mathcal{F}_k}} \right]\\
\le& \! \frac{{1\!+\!{\sigma_A ^2}}}{2}\mathbb{E}\left[ {\left\| {{v_k} \!-\! {A_\infty }{v_k}} \right\|_\pi ^2|{\mathcal{F}_k}} \right]\!+\!\frac{2}{{1 \!-\! {\sigma_A ^2}}}\mathbb{E}\left[ {\left\| {{g_{k +\! 1}} \!-\! {g_k}} \right\|_\pi ^2|{\mathcal{F}_k}} \right],
\end{aligned}
\end{equation}
where the inequality applies Lemma \ref{L1} and the Young's inequality. We next handle $\mathbb{E}\left[ {\left\| {{g_{k + 1}} - {g_k}} \right\|_\pi ^2|{\mathcal{F}_k}} \right]$ as follows:
\begin{equation}\label{E4-31}
\begin{aligned}
&\mathbb{E}\left[ {\left\| {{g_{k + 1}} - {g_k}} \right\|_\pi ^2|{\mathcal{F}_k}} \right]\\
\le& 2\mathbb{E}\left[ {\left\| {{g_{k + 1}} - {g_k} - \nabla F\left( {{z_{k + 1}}} \right) + \nabla F\left( {{z_k}} \right)} \right\|_\pi ^2|{\mathcal{F}_k}} \right]\\
&+ 2\left\| {\nabla F\left( {{z_{k + 1}}} \right) - \nabla F\left( {{z_k}} \right)} \right\|_\pi ^2\\
\le& 2{{\underline{\pi } }^{ - 1}}{L^2}\left\| {{z_{k + 1}} - {z_k}} \right\|_2^2+ 4{{\underline{\pi } }^{ - 1}}\mathbb{E}\left[ {\left\| {{g_k} - \nabla F\left( {{z_k}} \right)} \right\|_\pi ^2|{\mathcal{F}_k}} \right]\\
&+ 4{{\underline{\pi } }^{ - 1}}\mathbb{E}\left[ {\mathbb{E}\left[ {\left\| {{g_{k + 1}} - \nabla F\left( {{z_{k + 1}}} \right)} \right\|_\pi ^2|{\mathcal{F}_{k + 1}}} \right]|{\mathcal{F}_k}} \right].
\end{aligned}
\end{equation}
It follows from (\ref{E4-12}) that
\begin{equation}\label{E4-32}
\begin{aligned}
&\mathbb{E}\left[ {\left\| {{g_k} - \nabla F\left( {{z_k}} \right)} \right\|_\pi ^2|{\mathcal{F}_k}} \right]\\
\le& 4\bar \pi {d_1}{L^2}\left\| {{x_k} - {A^\infty }{x_k}} \right\|_\pi ^2 + 4m{L^2}\left\| {{{\bar x}_k} - {{\tilde z}^*}} \right\|_2^2+2{\Delta _k}\\
&+ 4{d_2}{L^2}{\sigma_A ^k}\left\| {{x_k}} \right\|_2^2,
\end{aligned}
\end{equation}
where the inequality applies (\ref{E4-16}). Similarly,
\begin{equation}\label{E4-33}
\begin{aligned}
&\mathbb{E}\left[ {\mathbb{E}\left[ {\left\| {{g_{k + 1}} - \nabla F\left( {{z_{k + 1}}} \right)} \right\|_\pi ^2|{\mathcal{F}_{k + 1}}} \right]|{\mathcal{F}_k}} \right]\\
\le& 4{L^2}\left\| {{z_{k + 1}} - {1_m} \otimes {{\bar x}_{k + 1}}} \right\|_2^2 + 4m{L^2}\left\| {{{\bar x}_{k + 1}} - {{\tilde z}^*}} \right\|_2^2\\
&+2\mathbb{E}\left[ {{\Delta _{k + 1}}|{\mathcal{F}_k}} \right]\\
\le& 4{d_1}\bar \pi {L^2}\left\| {{x_{k + 1}} - {A_\infty }{x_{k + 1}}} \right\|_\pi ^2 + 4{d_2}{L^2}{\sigma_A ^{k + 1}}\left\| {{x_{k + 1}}} \right\|_2^2\\
&+ 4m{L^2}\left\| {{{\bar x}_{k + 1}} - {{\tilde z}^*}} \right\|_2^2+2\mathbb{E}\left[ {{\Delta _{k + 1}}|{\mathcal{F}_k}} \right]\\
\le&8\bar \pi {d_1}{L^2}{\alpha ^2}\mathbb{E}\left[ {\left\| {{v_k} \!-\! {A_\infty }{v_k}} \right\|_\pi ^2|{\mathcal{F}_k}} \right] \!+\! 4m{L^2}\left( {3 \!+\! \bar p} \right)\left\| {{{\bar x}_k} \!-\! {{\tilde z}^*}} \right\|_2^2\\
&+4\bar \pi {d_1}{L^2}\left( {6{L^2}{\alpha ^2} + \bar p + 2} \right)\left\| {{x_k} - {A_\infty }{x_k}} \right\|_\pi ^2\\
&+ 2\left( {\frac{{4{L^2}{\alpha ^2}}}{m} + \left( {1 - \underline{p}} \right)} \right){\Delta _k}+8{d_2}{L^2}{\alpha ^2}{\sigma_A ^k}\mathbb{E}\left[ {\left\| {{v_k}} \right\|_2^2|{\mathcal{F}_k}} \right]\\
&+ 4{d_2}{L^2}\left( {6{L^2}{\alpha ^2} + \bar p + 2} \right){\sigma_A ^k}\left\| {{x_k}} \right\|_2^2,
\end{aligned}
\end{equation}
where the last inequality uses the results from Lemmas \ref{L2}, \ref{L4}-\ref{L5}. Considering $0<\underline{p}\le \bar p \le 1$ and picking the step-size as $0 < \alpha  \le 1/\left( {2L\sqrt 6 } \right)$, one can further obtain
\begin{equation}\label{E4-34}
\begin{aligned}
&\mathbb{E}\left[ {\mathbb{E}\left[ {\left\| {{g_{k + 1}} - \nabla F\left( {{z_{k + 1}}} \right)} \right\|_\pi ^2|{\mathcal{F}_{k + 1}}} \right]|{\mathcal{F}_k}} \right]\\
\le& 13\bar \pi {d_1}{L^2}\left\| {{x_k} \!-\! {A_\infty }{x_k}} \right\|_\pi ^2 \!+\! 16m{L^2}\left\| {{{\bar x}_k} - {{\tilde z}^*}} \right\|_2^2 \!+\! \frac{{6m \!+\! 1}}{{3m}}{\Delta _k}\\
&+ 8\bar \pi {d_1}{L^2}{\alpha ^2}\left\| {{v_k} - {A_\infty }{v_k}} \right\|_\pi ^2 + 13{d_2}{L^2}{\sigma_A ^k}\left\| {{x_k}} \right\|_2^2\\
&+8{d_2}{L^2}{\alpha ^2}{\sigma_A ^k}\mathbb{E}\left[ {\left\| {{v_k}} \right\|_2^2|{\mathcal{F}_k}} \right].
\end{aligned}
\end{equation}
We continue to handle $\left\| {{z_{k + 1}} - {z_k}} \right\|_2^2$ in (\ref{E4-31}) as follows:
\begin{equation}\label{E4-35}
\begin{aligned}
&\left\| {{z_{k + 1}} - {z_k}} \right\|_2^2\\
=& \left\| {Y_{k + 1}^{ - 1}\left( {A - {I_{mn}}} \right){x_k} - \alpha Y_{k + 1}^{ - 1}{v_k} + \left( {Y_{k + 1}^{ - 1} - Y_k^{ - 1}} \right){x_k}} \right\|_2^2\\
\le&\left\| {Y_{k + 1}^{ - 1}\left( {A - {I_{mn}}} \right)\left( {{x_k} - {A_\infty }{x_k}} \right)} \right\|_2^2+{{\tilde Y}^2}{\alpha ^2}\left\| {{v_k}} \right\|_2^2\\
&+ \left\| {Y_k^{ - 1}\left( {{Y_k} - {Y_{k + 1}}} \right)Y_{k + 1}^{ - 1}} \right\|_2^2\left\| {{x_k}} \right\|_2^2\\
&+2{\left\| {Y_{k + 1}^{ - 1}\left( {A - {I_{mn}}} \right){x_k}} \right\|_2}{\left\| {\alpha Y_{k + 1}^{ - 1}{v_k}}
\right\|_2}\\
&+2{\left\| {\alpha Y_{k + 1}^{ - 1}{v_k}} \right\|_2}{\left\| {Y_{k + 1}^{ - 1} - Y_k^{ - 1}} \right\|_2}{\left\| {{x_k}} \right\|_2}\\
&+ 2{\left\| {Y_{k + 1}^{ - 1}\left( {A - {I_{mn}}} \right){x_k}} \right\|_2}{\left\| {Y_{k + 1}^{ - 1} - Y_k^{ - 1}} \right\|_2}{\left\| {{x_k}} \right\|_2}\\
\le&12\bar \pi {{\tilde Y}^2}\left\| {{x_k} - {A_\infty }{x_k}} \right\|_\pi ^2 + 3{{\tilde Y}^2}{\alpha ^2}\left\| {{v_k}} \right\|_2^2\\
&+ 12{T^2}{{\tilde Y}^4}{\sigma_A ^{2k}}\left\| {{x_k}} \right\|_2^2,
\end{aligned}
\end{equation}
where the last inequality applies Lemma \ref{L2} and the results from \cite[Lemma 8]{Xi2018d}. Then, taking the total expectation on the both sides of (\ref{E4-35}) yields
\begin{equation}\label{E4-36}
\begin{aligned}
&\mathbb{E}\left[ {\left\| {{z_{k + 1}} - {z_k}} \right\|_2^2} \right]\\
\le& 12\bar \pi {\tilde Y^2}\mathbb{E}\left[ {\left\| {{x_k} - {A_\infty }{x_k}} \right\|_\pi ^2} \right] + 3{\tilde Y^2}{\alpha ^2}\mathbb{E}\left[ {\left\| {{v_k}} \right\|_2^2|{\mathcal{F}_k}} \right]\\
&+ 12{T^2}{\tilde Y^4}{\sigma_A ^{2k}}\mathbb{E}\left[ {\left\| {{x_k}} \right\|_2^2} \right].
\end{aligned}
\end{equation}
We next handle $\mathbb{E}\left[ {\left\| {{v_k}} \right\|_2^2|{\mathcal{F}_k}} \right]$ in both (\ref{E4-34}) and (\ref{E4-36}) as follows:
\begin{equation}\label{E4-37}
\begin{aligned}
&\mathbb{E}\left[ {\left\| {{v_k}} \right\|_2^2|{\mathcal{F}_k}} \right]\\
\le& 3\mathbb{E}\left[ {\left\| {{v_k} - {Y_\infty }\left( {{1_m} \otimes {{\bar g}_k}} \right)} \right\|_2^2|{\mathcal{F}_k}} \right] + 3\left\| {{Y_\infty }\left( {{1_m} \otimes {{\tilde p}_k}} \right)} \right\|_2^2\\
&+ 3\mathbb{E}\left[ {\left\| {{Y_\infty }\left( {{1_m} \otimes {{\bar g}_k}} \right) - {Y_\infty }\left( {{1_m} \otimes {{\tilde p}_k}} \right)} \right\|_2^2|{\mathcal{F}_k}} \right]\\
\le& 3\bar \pi \mathbb{E}\left[ {\left\| {{v_k} - {A_\infty }{v_k}} \right\|_\pi ^2|{\mathcal{F}_k}} \right] + 6m{Y^2}\mathbb{E}\left[ {\left\| {{{\bar g}_k} - {{\bar h}_k}} \right\|_2^2|{\mathcal{F}_k}} \right]\\
&+ 6m{Y^2}\left\| {{{\bar h}_k} - {{\tilde p}_k}} \right\|_2^2 + 3m{Y^2}\left\| {{{\tilde p}_k} - \nabla \tilde f\left( {{\tilde z^*}} \right)} \right\|_2^2\\
\le& 3\bar \pi \mathbb{E}\left[ {\left\| {{v_k} - {A_\infty }{v_k}} \right\|_\pi ^2|{\mathcal{F}_k}} \right] + 6m{Y^2}\mathbb{E}\left[ {\left\| {{{\bar g}_k} - {{\bar h}_k}} \right\|_2^2|{\mathcal{F}_k}} \right]\\
&+ 6{L^2}{Y^2}\left\| {{z_k} - {1_m} \otimes {{\bar x}_k}} \right\|_2^2 + 3{L^2}{Y^2}\left\| {{{\bar x}_k} - {{\tilde z}^*}} \right\|_2^2\\
\le&3\bar \pi \mathbb{E}\left[ {\left\| {{v_k} - {A_\infty }{v_k}} \right\|_\pi ^2|{\mathcal{F}_k}} \right] + 30\bar \pi {d_1}{L^2}{Y^2}\left\| {{x_k} - {A_\infty }{x_k}} \right\|_\pi ^2\\
&+ 27{L^2}{Y^2}\left\| {{{\bar x}_k} - {{\tilde z}^*}} \right\|_2^2 + \frac{{12{Y^2}}}{m}{\Delta _k} + 30{d_2}{L^2}{Y^2}{\sigma_A ^k}\left\| {{x_k}} \right\|_2^2,
\end{aligned}
\end{equation}
where the third inequality uses the $L$-smoothness of the local objective functions, and the last inequality applies (\ref{E4-13}) and (\ref{E4-16}). Then, combining (\ref{E4-37}) and (\ref{E4-34}) obtains
\begin{equation}\label{E4-38}
\begin{aligned}
&\mathbb{E}\left[ {\mathbb{E}\left[ {\left\| {{g_{k + 1}} - \nabla F\left( {{z_{k + 1}}} \right)} \right\|_\pi ^2|{\mathcal{F}_{k + 1}}} \right]|{\mathcal{F}_k}} \right]\\
\le&\bar \pi {d_1}{L^2}\left( {13 + 240{d_2}{Y^2}{\alpha ^2}{L^2}{\sigma_A ^k}} \right)\left\| {{x_k} - {A_\infty }{x_k}} \right\|_\pi ^2\\
&+ 8{L^2}\left( {2m + 27{d_2}{Y^2}{\alpha ^2}{L^2}{\sigma_A ^k}} \right)\left\| {{{\bar x}_k} - {{\tilde z}^*}} \right\|_2^2\\
&+ \frac{1}{m}\left( {\frac{{6m + 1}}{3} + 96{d_2}{Y^2}{L^2}{\alpha ^2}{\sigma_A ^k}} \right){\Delta _k}\\
&+8\bar \pi {L^2}{\alpha ^2}\left( {{d_1} + 3{d_2}{\sigma_A ^k}} \right)\left\| {{v_k} - {A_\infty }{v_k}} \right\|_\pi ^2\\
& + {d_2}{L^2}\left( {13 + 240{d_2}{Y^2}{\alpha ^2}{L^2}{\sigma_A ^k}} \right){\sigma_A ^k}\left\| {{x_k}} \right\|_2^2,
\end{aligned}
\end{equation}
and plugging (\ref{E4-37}) into (\ref{E4-35}) gives
\begin{equation}\label{E4-39}
\begin{aligned}
&\mathbb{E}\left[ {\left\| {{z_{k + 1}} - {z_k}} \right\|_2^2|{\mathcal{F}_k}} \right]\\
\le&6\bar \pi {{\tilde Y}^2}\left( {2 + 15{d_1}{L^2}{Y^2}{\alpha ^2}} \right)\left\| {{x_k} - {A_\infty }{x_k}} \right\|_\pi ^2 \!+\! \frac{{36{Y^2}{{\tilde Y}^2}{\alpha ^2}}}{m}{\Delta _k}\\
&+\! 81{L^2}{Y^2}{{\tilde Y}^2}{\alpha ^2}\left\| {{{\bar x}_k} \!-\! {{\tilde z}^*}} \right\|_2^2  \!+\! 9\bar \pi {{\tilde Y}^2}{\alpha ^2}\mathbb{E}[ {\| {{v_k} \!-\! {A_\infty }{v_k}} \|_\pi ^2|{\mathcal{F}_k}} ]\\
&+ 6{{\tilde Y}^2}\left( {15{d_2}{L^2}{Y^2}{\alpha ^2} + 2{T^2}{{\tilde Y}^2}{\sigma_A ^k}} \right){\sigma_A ^k}\left\| {{x_k}} \right\|_2^2.
\end{aligned}
\end{equation}
Then, plugging (\ref{E4-32}), (\ref{E4-38}), and (\ref{E4-39}) into (\ref{E4-31}) reduces to
\begin{equation}\label{E4-40}
\begin{aligned}
&\mathbb{E}\left[ {\left\| {{g_{k + 1}} - {g_k}} \right\|_\pi ^2|{\mathcal{F}_k}} \right]\\
\le& \vartheta {L^2}\left( {24{{\tilde Y}^2} + 68{d_1} + 60{d_1}{L^2}{Y^2}\left( {3{{\tilde Y}^2} + 16{d_2}} \right){\alpha ^2}} \right)\\
&\times\left\| {{x^k} - {A_\infty }{x^k}} \right\|_\pi ^2\\
&+ {{\underline{\pi } }^{ - 1}}{L^2}\left( {80m + 54{L^2}{Y^2}\left( {3{{\tilde Y}^2} \!+\! 16{d_2}} \right){\alpha ^2}} \right)\left\| {{{\bar x}_k} - {{\tilde z}^*}} \right\|_2^2\\
&+ \frac{{4{{\underline{\pi } }^{ - 1}}}}{m}\left( {\frac{{12m + 1}}{3} + \left( {18{{\tilde Y}^2} + 96{d_2}} \right){L^2}{Y^2}{\alpha ^2}} \right){\Delta ^k}\\
&+\vartheta \left( {9{{\tilde Y}^2} + 16{d_1} + 48{d_2}} \right){L^2}{\alpha ^2}\mathbb{E}\left[ {\left\| {{v_k} - {A_\infty }{v_k}} \right\|_\pi ^2|{\mathcal{F}_k}} \right]\\
&+\left( {24{T^2}{{\tilde Y}^4}{\sigma_A ^k} + 68{d_2} + 60{d_2}{L^2}{Y^2}\left( {3{{\tilde Y}^2} + 16{d_2}{\sigma_A ^k}} \right){\alpha ^2}} \right)\\
&\times {{\underline{\pi } }^{ - 1}}{L^2}{\sigma_A ^k}\left\| {{x_k}} \right\|_2^2,
\end{aligned}
\end{equation}
where the inequality uses the fact that ${0 < \sigma_A  < 1}$. Then, combing (\ref{E4-30}) with (\ref{E4-40}) obtains
\begin{equation}\label{E4-41}
\begin{aligned}
&\mathbb{E}\left[ {\left\| {{v_{k + 1}} - {A_\infty }{v_{k + 1}}} \right\|_\pi ^2|{\mathcal{F}_k}} \right]\\
\le& \left( {\frac{{1+{\sigma_A ^2}}}{2} + \frac{{2\vartheta \left( {9{{\tilde Y}^2} + 16{d_1} + 48{d_2}} \right){L^2}{\alpha ^2}}}{{1 - {\sigma_A ^2}}}} \right)\\
&\times\mathbb{E}\left[ {\left\| {{v_k} - {A_\infty }{v_k}} \right\|_\pi ^2|{\mathcal{F}_k}} \right]\\
&+\left( {\frac{{48{{\tilde Y}^2} + 136{d_1}}}{{1 - {\sigma_A ^2}}} + \frac{{120{d_1}\left( {3{{\tilde Y}^2} + 16{d_2}} \right){L^2}{Y^2}{\alpha ^2}}}{{1 - {\sigma_A ^2}}}} \right)\\
&\times \vartheta {L^2}\left\| {{x_k} - {A_\infty }{x_k}} \right\|_\pi ^2\\
&+\!\left( {\frac{{160m}}{{1 \!-\! {\sigma_A ^2}}} \!+\! \frac{{108( {3{{\tilde Y}^2} \!+\! 16{d_2}} ){L^2}{Y^2}{\alpha ^2}}}{{1 - {\sigma_A ^2}}}} \right){{\underline{\pi } }^{ \!- 1}}{L^2}\left\| {{{\bar x}_k}\! - \! {{\tilde z}^*}} \right\|_2^2\\
&+ \left( {\frac{{96m + 8}}{{3m\left( {1 - {\sigma_A ^2}} \right)}} + \frac{{24\left( {3{{\tilde Y}^2} + 16{d_2}} \right){L^2}{Y^2}{\alpha ^2}}}{{m\left( {1 - {\sigma_A ^2}} \right)}}} \right){{\underline{\pi } }^{ - 1}}{\Delta _k}\\
&+\left( {\frac{{48{T^2}{{\tilde Y}^4} + 136{d_2}}}{{1 - {\sigma_A ^2}}} + \frac{{120{d_2}\left( {3{{\tilde Y}^2} + 16{d_2}} \right){L^2}{Y^2}{\alpha ^2}}}{{1 - {\sigma_A ^2}}}} \right)\\
&\times{{\underline{\pi } }^{ - 1}}{L^2}{\sigma_A ^k}\left\| {{x_k}} \right\|_2^2.
\end{aligned}
\end{equation}
Via choosing $0 < \alpha  \le 1/\left( {2LY\sqrt {3\left( {3{{\tilde Y}^2} + 16{d_2}} \right)} } \right)$ for the first term and $0 < \alpha  \le \left( {1 - {\sigma_A ^2}} \right)/\left( {2L\sqrt {2\vartheta \left( {9{{\tilde Y}^2} + 16{d_1} + 48{d_2}} \right)} } \right)$ for the rest terms in the left-hand-side of (\ref{E4-41}), one can obtain
\begin{equation}\label{E4-42}
\begin{aligned}
&\mathbb{E}\left[ {\left\| {{v_{k + 1}} - {A_\infty }{v_{k + 1}}} \right\|_\pi ^2|{\mathcal{F}_k}} \right]\\
\le& \left( {\frac{{3+{\sigma_A ^2}}}{4}} \right)\mathbb{E}\left[ {\left\| {{v_k} - {A_\infty }{v_k}} \right\|_\pi ^2|{\mathcal{F}_k}} \right]+\left( {\frac{{48{{\tilde Y}^2} + 146{d_1}}}{{1 - {\sigma_A ^2}}}} \right)\\
&\times \vartheta {L^2}\| {{x_k} - {A_\infty }{x_k}} \|_\pi ^2 + \left( {\frac{{160m + 9}}{{1 \!-\! {\sigma_A ^2}}}} \right){{\underline{\pi } }^{ - 1}}{L^2}\| {{{\bar x}_k} - {{\tilde z}^*}}\|_2^2\\
&+\left( {\frac{{48{T^2}{{\tilde Y}^4} + 146{d_2}}}{{1 - {\sigma_A ^2}}}} \right){{\underline{\pi } }^{ - 1}}{L^2}{\sigma_A ^k}\left\| {{x_k}} \right\|_2^2.
\end{aligned}
\end{equation}
Recalling the definitions of ${d_1}$, $d_2$, and $\delta $, it holds that $\vartheta {\tilde Y^2} \le \vartheta {d_1} \le \delta $. Therefore, one can rewrite (\ref{E4-42}) as follows:
\begin{equation}\label{E4-43}
\begin{aligned}
&\mathbb{E}\left[ {\left\| {{v_{k + 1}} - {A_\infty }{v_{k + 1}}} \right\|_\pi ^2|{\mathcal{F}_k}} \right]\\
\le& \left( {\frac{{3\!+\!{\sigma_A ^2}}}{4}} \right)\mathbb{E}\left[ {\left\| {{v_k} \!-\! {A_\infty }{v_k}} \right\|_\pi ^2|{\mathcal{F}_k}} \right]\!+\!\frac{{194\delta {L^2}}}{{1 \!-\! {\sigma_A ^2}}}\left\| {{x_k} \!-\! {A_\infty }{x_k}} \right\|_\pi ^2\\
&+\frac{{169m{{\underline{\pi } }^{ - 1}}{L^2}}}{{1 - {\sigma_A ^2}}}\left\| {{{\bar x}_k} - {{\tilde z}^*}} \right\|_2^2 + \frac{{110{{\underline{\pi } }^{ - 1}}}}{{3\left( {1 - {\sigma_A ^2}} \right)}}{\Delta _k}\\
&+\frac{{194{{\underline{\pi } }^{ - 1}}T{\delta ^2}{L^2}{\sigma_A ^k}}}{{1 - {\sigma_A ^2}}}\left\| {{x_k}} \right\|_2^2.
\end{aligned}
\end{equation}
Taking the total expectation on the both sides of (\ref{E4-43}) completes the proof.
\end{proof}
\end{lemma}
Based on Lemmas \ref{L2}, \ref{L4}-\ref{L6}, it is straightforward to build a DLTI system in the following proposition.
\begin{proposition}\label{P1}
Suppose that Assumptions \ref{A1}-\ref{A3} hold.
If the step-size satisfies $0 < \alpha  \le \left( {1 - {\sigma_A ^2}} \right)\sqrt {\underline{p} } /\left( {28L\mathcal{Q}\delta \sqrt {\bar p} } \right)$, then $\forall k \ge 0$, the following DLTI system inequality holds
\begin{equation}\label{E4-44}
{t_{k + 1}} \leq {H_\alpha }{t_k} + {G_k}{\tau _k},
\end{equation}
where the inequality is taken element-wise and the vectors are defined as:
\[{t_k} :=  \left[ {\begin{array}{*{20}{c}}
  {\mathbb{E}\left[ {\left\| {{x_k} - {A_\infty }{x_k}} \right\|_\pi ^2} \right]} \\
  {\mathbb{E}\left[ {m\left\| {{{\bar x}_k} - {{\tilde z}^*}} \right\|_2^2} \right]} \\
  {\mathbb{E}\left[ {{\Delta _k}} \right]} \\
  {\mathbb{E}\left[ {\left\| {{v_k} - {A_\infty }{v_k}} \right\|_\pi ^2} \right]}
\end{array}} \right],{\tau _k} := \left[ {\begin{array}{*{20}{c}}
  {\mathbb{E}\left[ {\left\| {{x_k}} \right\|_2^2} \right]} \\
  0 \\
  0 \\
  0
\end{array}} \right],\]
and the matrices are indicated by
\[{H_\alpha } := \left[ {\begin{array}{*{20}{c}}
  {\frac{{1 + {\sigma_A ^2}}}{2}}&0&0&{\frac{{2{\alpha ^2}}}{{1 - {\sigma_A ^2}}}} \\
  {\frac{{2\bar \pi \delta \alpha {L^2}}}{\mu }}&{1 - \frac{{\mu \alpha }}{2}}&{\frac{{2{\alpha ^2}}}{m}}&0 \\
  {2\bar p\bar \pi {d_1}{L^2}}&{2m\bar p{L^2}}&{1 - \underline{p} }&0 \\
  {\frac{{194\delta {L^2}}}{{1 - {\sigma_A ^2}}}}&{\frac{{169{{\underline{\pi } }^{ - 1}}{L^2}}}{{1 - {\sigma_A ^2}}}}&{\frac{{110{{\underline{\pi } }^{ - 1}}}}{{3\left( {1 - {\sigma_A ^2}} \right)}}}&{\frac{{3 +{\sigma_A ^2}}}{4}}
\end{array}} \right],\]
\[{G_k} := \left[ {\begin{array}{*{20}{c}}
  0&0&0&0 \\
  {\frac{{2\delta \alpha {L^2}}}{\mu }}&0&0&0 \\
  {2\bar p{d_2}{L^2}}&0&0&0 \\
  {\frac{{194{{\underline{\pi } }^{ - 1}}{L^2}{\delta ^2}}}{{1 - {\sigma_A ^2}}}}&0&0&0
\end{array}} \right]T{\sigma_A ^k}.\]
\end{proposition}
\noindent With the help of Proposition \ref{P1}, we aim at solving for a specific range of the constant step-size $\alpha$ to guarantee that $\rho \left( {{H_\alpha }} \right) < 1$ (see Appendix), which is a necessary condition to establish the linear convergence of Algorithm \ref{Algo1}. Note that $G_k$ decays linearly at the rate of $\sigma_A$.
\section{Experimental Results}\label{Section 5}
To manifest the effectiveness and  practicability of \textit{Push-LSVRG-UP}, we provide two case studies to compare them with existing state-of-art distributed algorithms over both unbalanced directed networks and undirected networks. The total number of training samples is denoted as $N$, which is randomly and evenly allocated among $m$ agents. Then, each agent $i$, $i \in \mathcal{V}$, maintains ${q_i} = N/m$ local samples. We denote ${c_{ij}} \in {\mathbb{R}^n}$ as the $j$-th training sample, and ${b_{ij}} \in \left\{ { + 1, - 1} \right\}$ is the corresponding label accessed only by agent $i$, $i\in\mathcal{V}$.
In the following simulations, the optimal gap is indicated by residual: $\left( {1/m} \right){\sum\nolimits_{i = 1}^m {\left\| {z_k^i - {{\tilde z}^*}} \right\|} _2}$, and each epoch indicates an effective pass of the local samples.
To investigate extensively the convergence performance of \textit{Push-LSVRG-UP}, we also consider the special case of \textit{Push-LSVRG-UP}, i.e., fixing the uncoordinated triggering probabilities $p_i$ as a coordinated one with $p_i = p$, $i \in \mathcal{V}$. This special case of \textit{Push-LSVRG-UP} is named Push-LSVRG in the sequel. Throughout the simulations,
The coordinated triggering probabilities of Push-LSVRG is fixed as $p=1/Q$ and the uncoordinated triggering probability of \textit{Push-LSVRG-UP} is randomly selected in an interval: $1/Q \le p_i \le m/Q$, $i \in \mathcal{V}$.
All simulations are carried out in MATLAB on a Dell PowerEdge R740 with 2.10 GHz, 26 Cores, 52 Threads, Intel Xeon Gold 6230R processor and 256GB memory.
\begin{figure}[!htp]
  \centering
  \includegraphics[width=8cm]{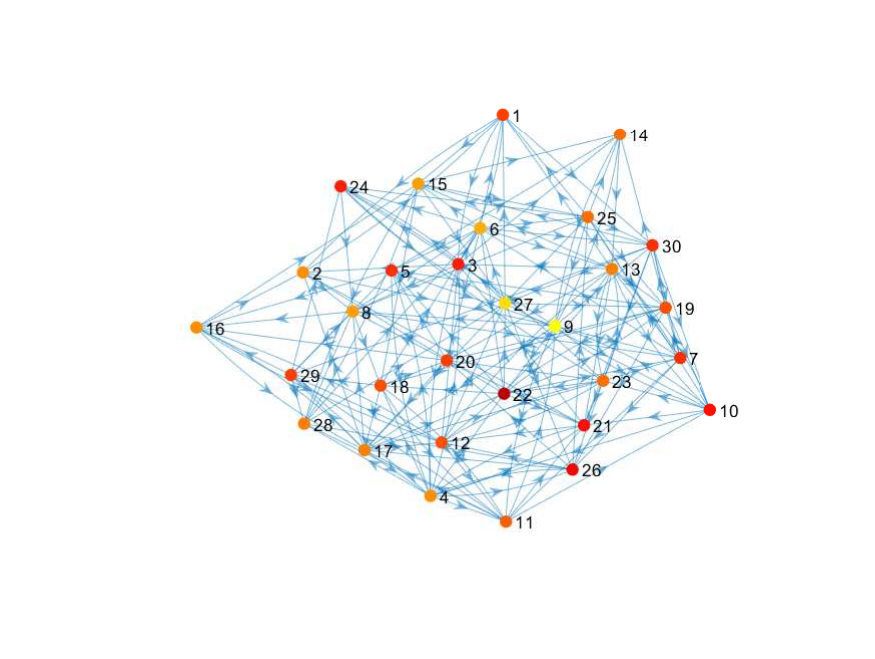}
  \caption{An unbalanced directed network with $m=30$.}\label{fig. 1}
\end{figure}
\subsection{Case Study One: Distributed Logistic Regression}\label{Section 5-1}
In the first case study, Push-LSVRG, \textit{Push-LSVRG-UP} and the other tested distributed algorithms are utilized to identify whether a mushroom is poisonous or not according to its different features, such as ``gill-color", ``stalk-root", ``veil-type",``cap-shape", ``habitat" and so on. Each feature may contains several options, for example, the options of  ``cap-shape" are varied from ``bell", ``conical", ``convex", ``flat", ``knobbed", and ``sunken". Mushroom dataset provided by UCI Machine Learning Repository \cite{Dua2013a} contains a total number of 8124 samples and each sample has $n=112$ dimensions that indicate different features. We randomly choose $N=6000$ samples from the total samples to train the discriminator and the rest of samples are used for testing.
\begin{figure}[!htp]
\centering
\subfigure[Training performance.]{\includegraphics[width=1.62in,height=1.30in]{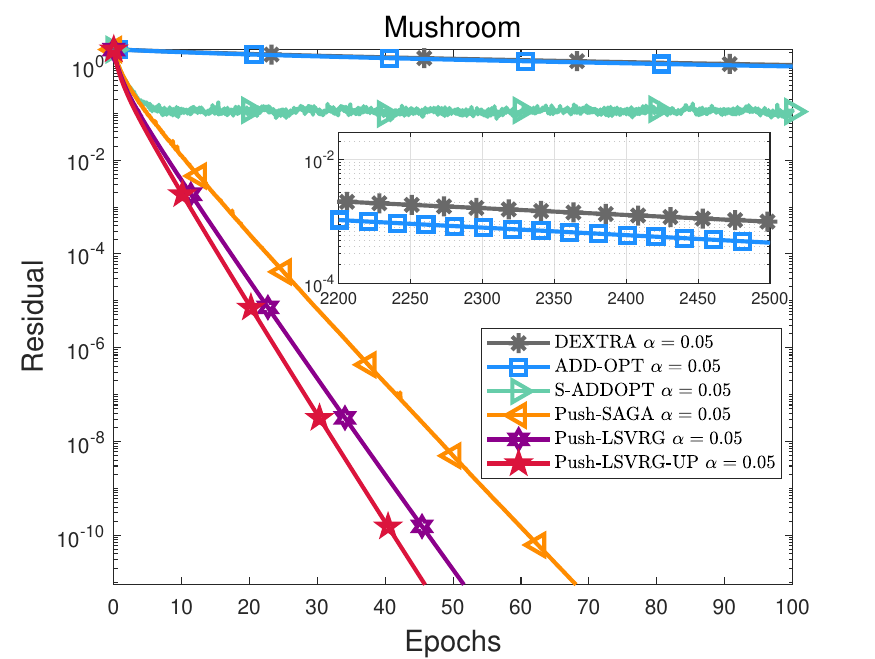}
\label{fig. 2}}
\hfil
\subfigure[Testing performance.]{\includegraphics[width=1.62in,height=1.30in]{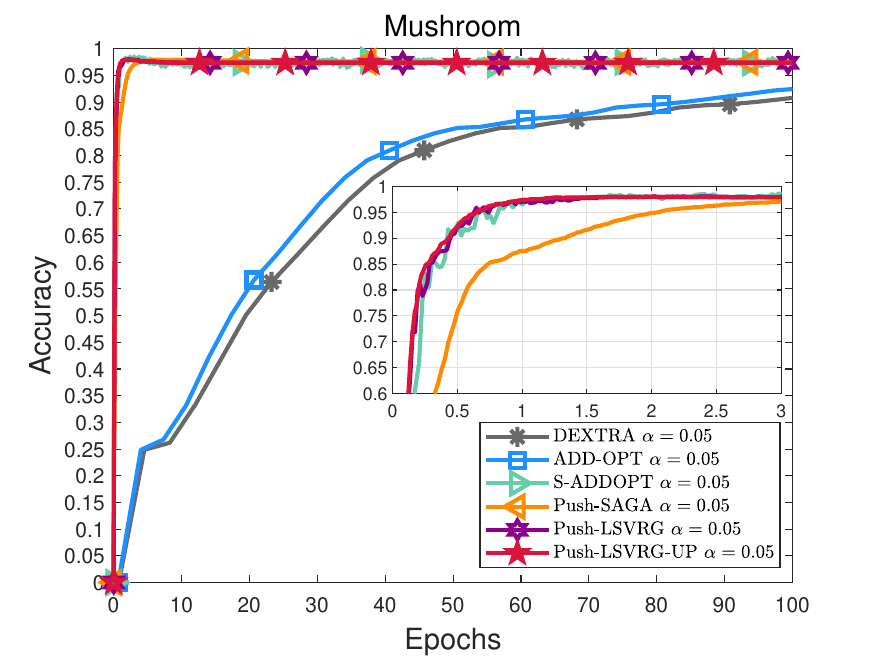}
\label{fig. 3}}
\caption{Performance comparison over epochs.}
\label{ConvPerfEpo}
\end{figure}
Specifically, a network of $m=30$ agents train cooperatively a regularized logistic regression model for binary classification as follows:
\begin{equation}\label{E5-1}
\mathop {\min }\limits_{\tilde z \in {\mathbb{R}^n}} \tilde {f}\left( {\tilde z} \right)\! := \! \frac{\beta }{2}\left\| {\tilde z} \right\|_2^2 + \frac{1}{m}\sum\limits_{i = 1}^m {\frac{1}{{{q_i}}}} \sum\limits_{j = 1}^{{q_i}} {\log \left( {1 \!+\! \exp \left( { - {b_{ij}}c_{ij}^{\top} \tilde z} \right)} \right)},
\end{equation}
where $\beta$ is the regularized constant and we set $\beta=5$ in this case study.
Note that ${b_{ij}} =  + 1$ when the sample ${c_{ij}}$ is poisonous while ${b_{ij}} =  - 1$ when the sample ${c_{ij}}$ is edible. To solve (\ref{E5-1}) in a distributed framework,
we conduct the simulation in a multi-agent system with strongly-connected unbalanced directed networks and each agent has $6$ out-neighbors as shown in Fig. \ref{fig. 1}.
Through utilizing a powerful "centrality" function in Matlab, Fig. \ref{fig. 1} maps different "authorities" of agents to various color.
Furthermore, the number of the total training samples is allocated equally to all agents.
Therefore, each agent $i$ maintains $q_i = N/m = 200$ training samples.
Figs. \ref{fig. 2}-\ref{02271422Training_time} and Figs. \ref{fig. 3}-\ref{02271422Testing_time} show respectively the training performance and testing performance of all the tested algorithms. From Figs. \ref{ConvPerfEpo}-\ref{ConvPerfTime}, one can apparently see that \textit{Push-LSVRG} and \textit{Push-LSVRG-UP} achieve high convergence and testing accuracy faster than all the other tested algorithms in terms of both epochs and CPU running time over an unbalanced directed network as shown in Fig. \ref{fig. 1}. In Fig. \ref{ConvPerfTime}, even though \textit{S-ADDOPT} \cite{Qureshi2021} also reaches this highest accuracy in a fast way, it is relatively unstable when achieving the highest testing accuracy. This phenomenon is understandable since \textit{S-ADDOPT} is an inexact distributed stochastic algorithm.
\begin{figure}[!htp]
\centering
\subfigure[Training performance.]{\includegraphics[width=1.62in,height=1.30in]{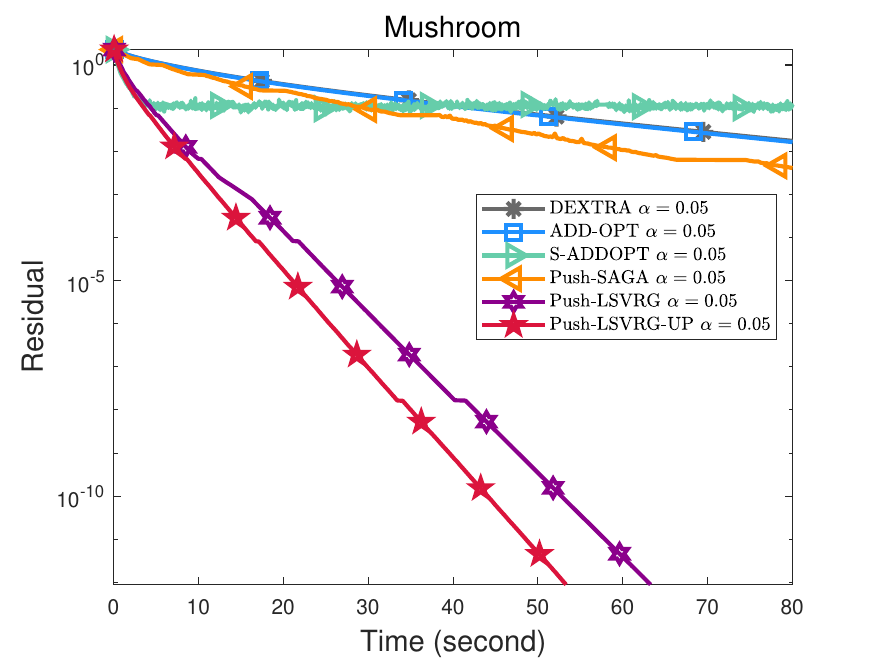}
\label{02271422Training_time}}
\hfil
\subfigure[Testing performance.]{\includegraphics[width=1.62in,height=1.30in]{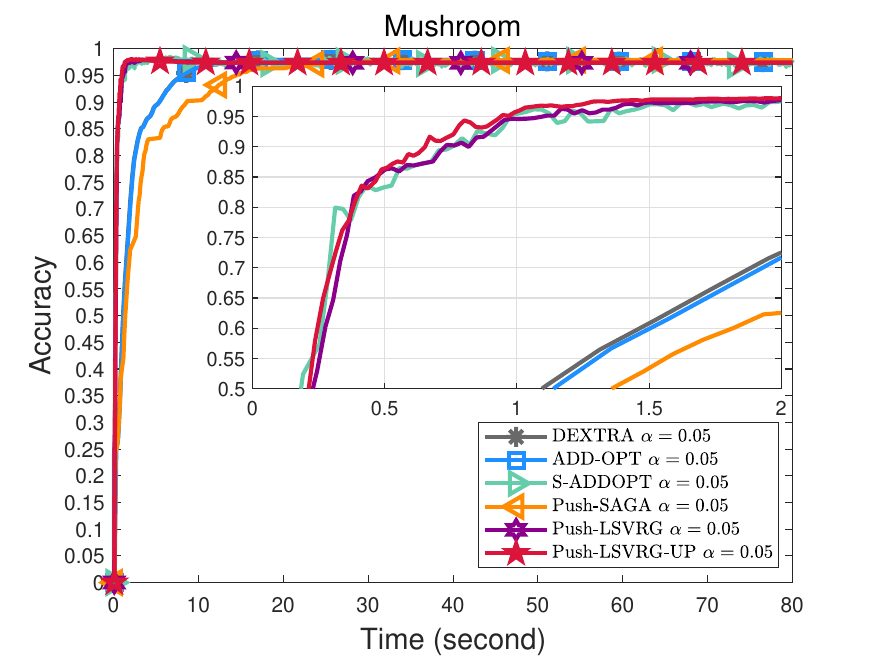}
\label{02271422Testing_time}}
\caption{Performance comparison over CPU running time.}
\label{ConvPerfTime}
\end{figure}
To further explore the impact of network sizes towards the convergence performance, we compare all the tested algorithms over three exponential networks \cite{Yuan2021} as shown in Figs. \ref{ExGra4}-\ref{ExGra16}. These exponential networks only differ on the number of agents and they share the same network structure.
\begin{figure}[!htp]
\centering
\subfigure[$m=4$.]{\includegraphics[width=1.10in,height=1.10in]{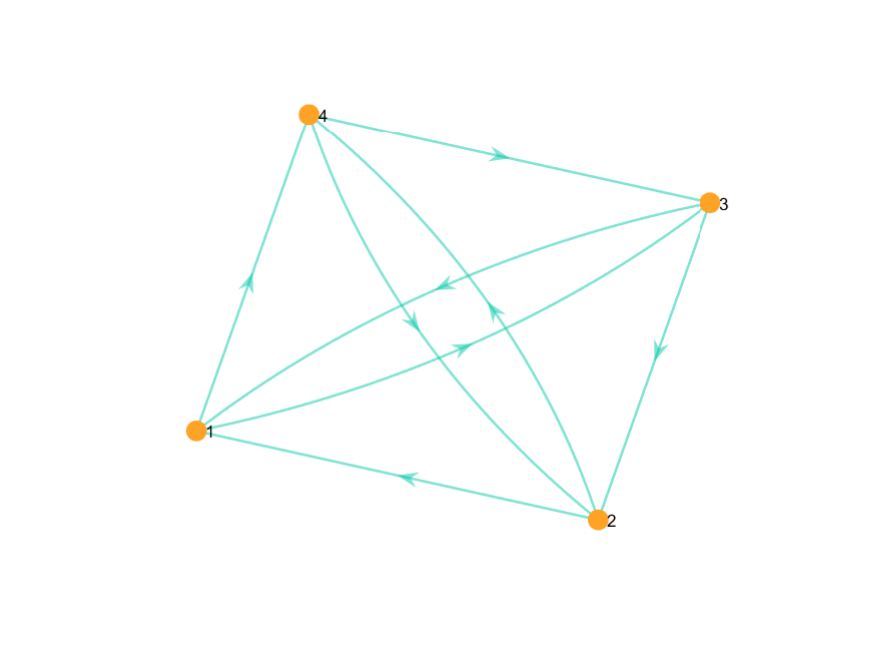}
\label{ExGra4}}
\hfil
\subfigure[$m=8$.]{\includegraphics[width=1.10in,height=1.10in]{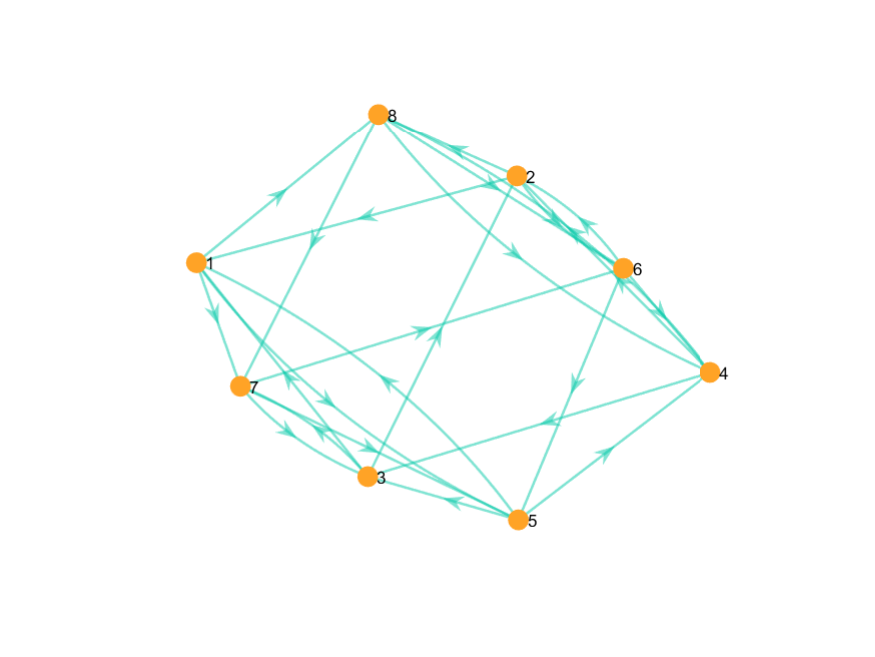}
\label{ExGra8}}
\hfil
\subfigure[$m=16$.]{\includegraphics[width=1.10in,height=1.10in]{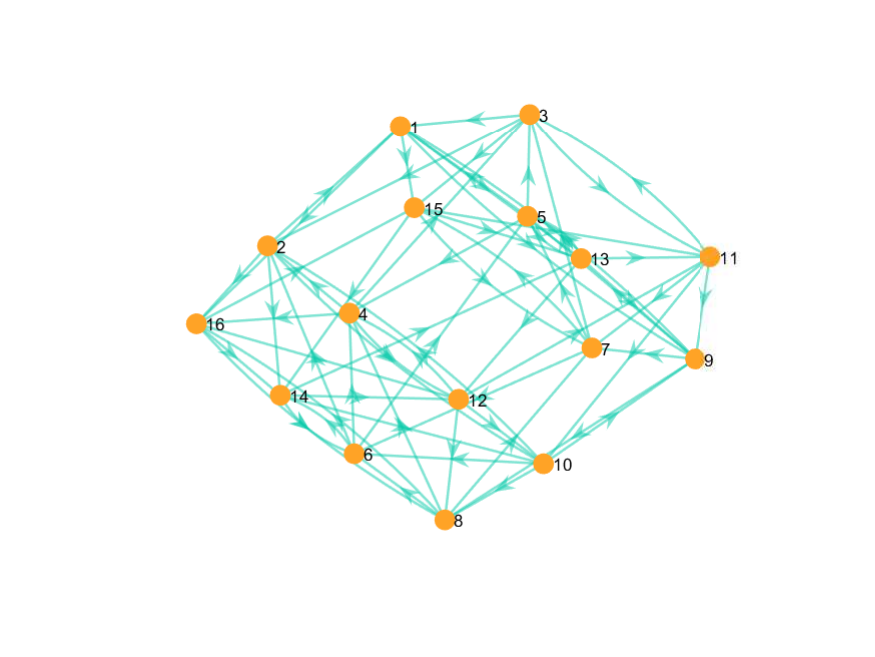}
\label{ExGra16}}
\caption{Exponential unbalanced directed networks with different agents.}
\label{fig. ExGra}
\end{figure}
\begin{figure}[!htp]
\centering
\subfigure[Performance comparison over Fig. \ref{ExGra4}.]{\includegraphics[width=1.10in,height=1.00in]{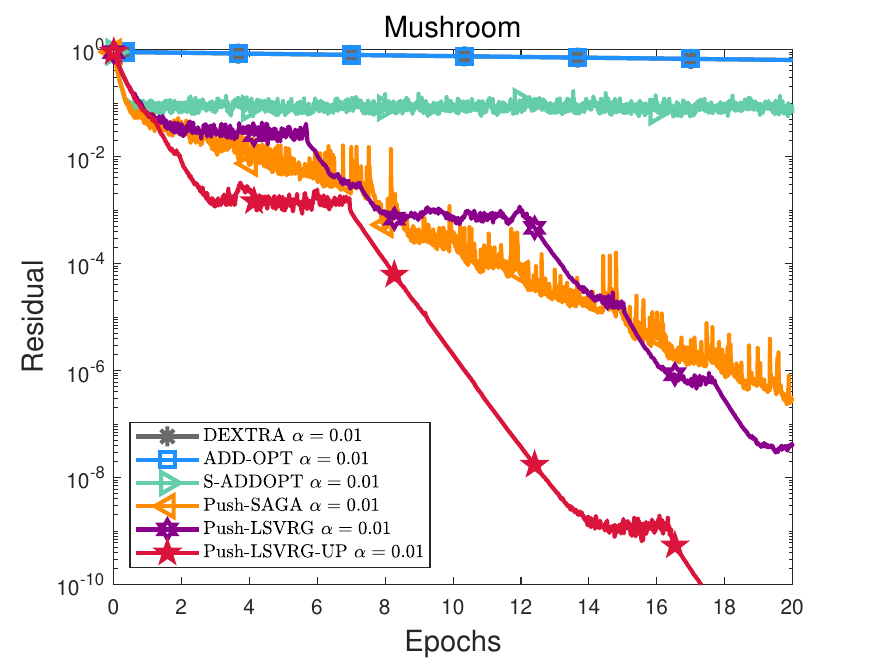}
\label{Training_epochs_ExGra4}}
\hfil
\subfigure[Performance comparison over Fig. \ref{ExGra8}.]{\includegraphics[width=1.10in,height=1.00in]{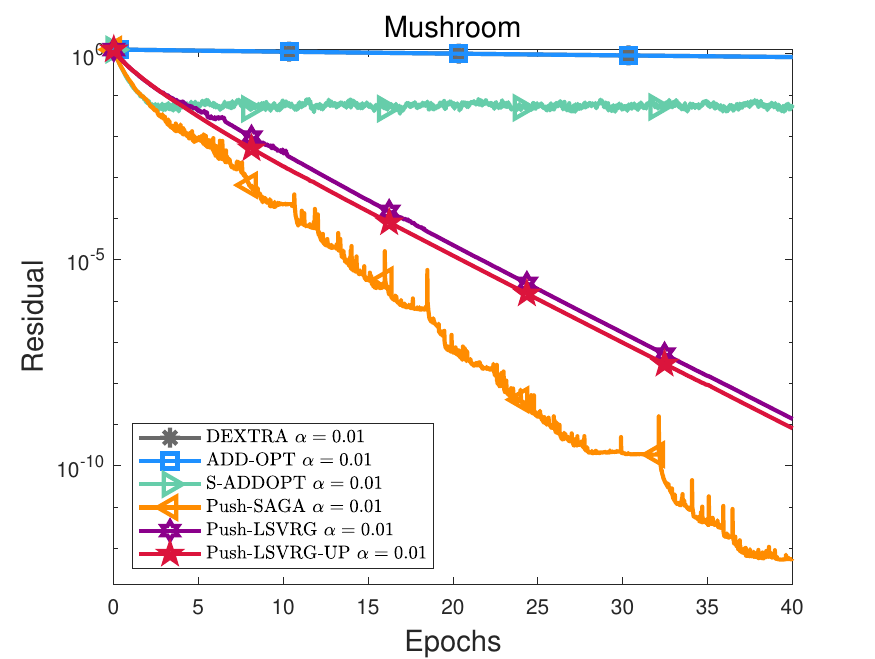}
\label{Training_epochs_ExGra8}}
\hfil
\subfigure[Performance comparison over Fig. \ref{ExGra16}.]{\includegraphics[width=1.10in,height=1.00in]{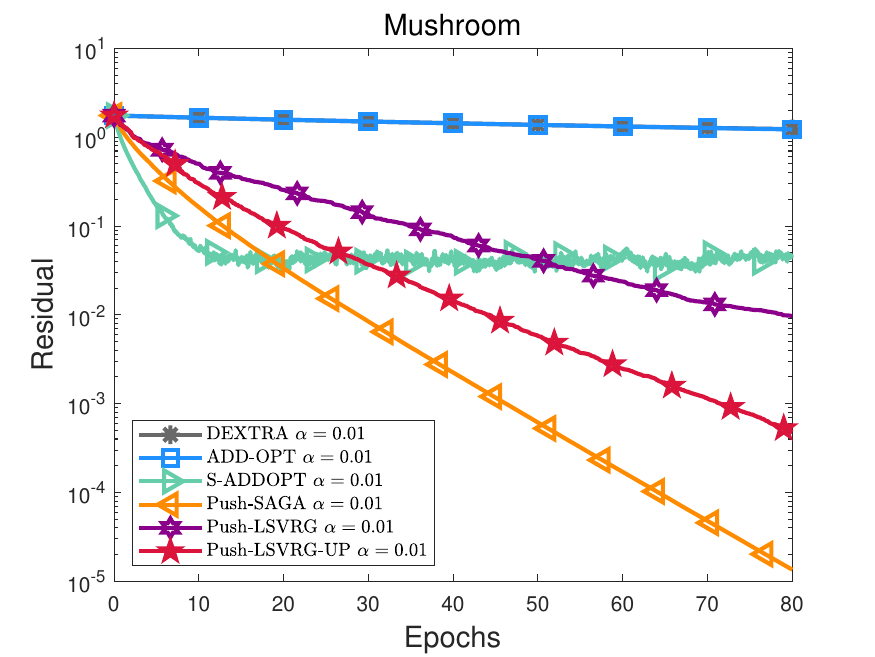}
\label{Training_epochs_ExGra16}}
\caption{Convergence performance comparison over epochs with different network sizes.}
\label{Training_epochs}
\end{figure}
\begin{figure}[!htp]
\centering
\subfigure[Performance comparison over Fig. \ref{ExGra4}.]{\includegraphics[width=1.10in,height=1.00in]{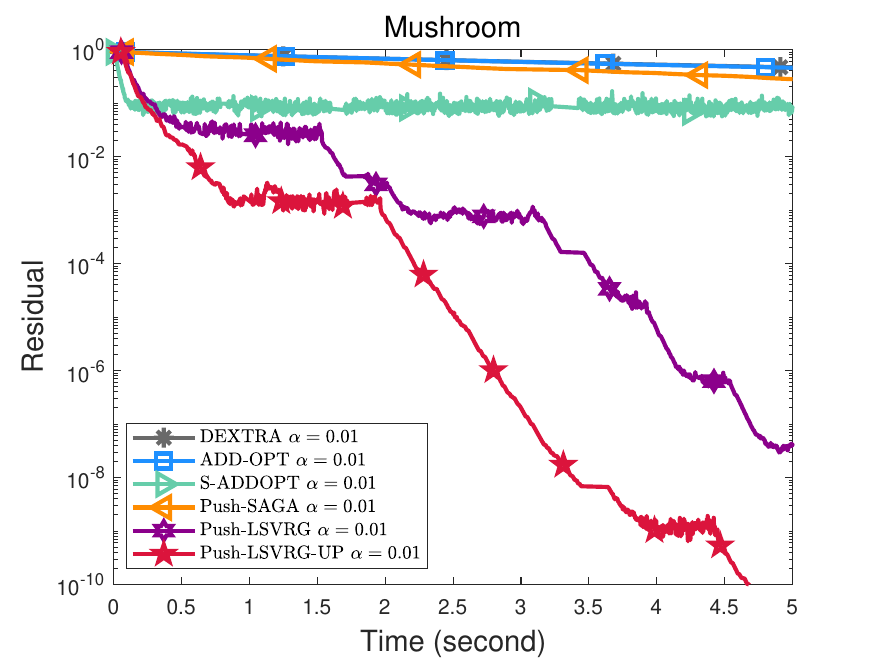}
\label{Training_time_ExGra4}}
\hfil
\subfigure[Performance comparison over Fig. \ref{ExGra8}.]{\includegraphics[width=1.10in,height=1.00in]{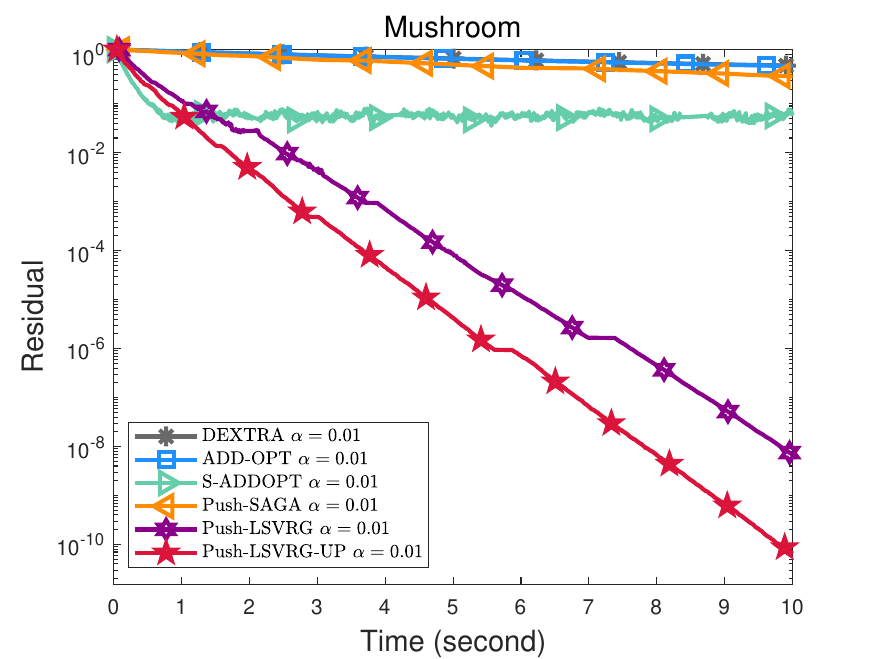}
\label{Training_time_ExGra8}}
\hfil
\subfigure[Performance comparison over Fig. \ref{ExGra16}.]{\includegraphics[width=1.10in,height=1.00in]{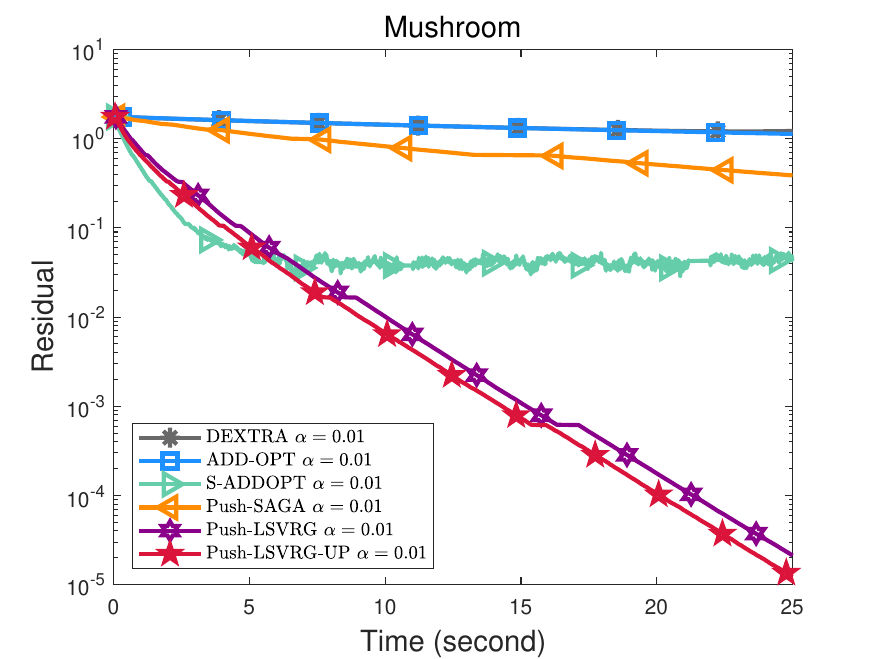}
\label{Training_time_ExGra16}}
\caption{Convergence performance comparison over CPU running time with different network sizes.}
\label{Training_time}
\end{figure}
It is worth mentioning that when the network size $m$ changes, the studied optimization problem also changes according to (\ref{E1-1}), which means that the global optimal solutions of all the tested algorithms over Figs \ref{ExGra4}-\ref{ExGra16} are different. Therefore, we cannot simply compare the convergence performance of all the tested algorithms over different network sizes. Nevertheless, some interest results can still be obtained from Figs. \ref{Training_epochs}-\ref{Training_time}. Specifically, when the network size is increased from Fig. \ref{ExGra4} to Fig. \ref{ExGra16}, Figs. \ref{Training_time_ExGra4}-\ref{Training_time_ExGra16} show that \textit{Push-LSVRG-UP} has certain acceleration than the other tested algorithms regarding CPU running time with different network sizes.
We note that \textit{Push-LSVRG-UP} does not require any additional storage, which is an important advantage to \textit{SAGA}-based algorithms, for instance \textit{Push-SAGA} \cite{Qureshi2020a}, since these algorithms requires an expensive storage cost of $\mathcal{O}\left( {nQ} \right)$ under the same structured large-scale optimization problem.
However, even though \textit{Push-LSVRG-UP} shows its priority in both CPU running time and a less requirement of storage to \textit{Push-SAGA}, Figs. \ref{Training_epochs_ExGra4}-\ref{Training_epochs_ExGra16} demonstrate that the convergence performance comparison in terms of epochs between \textit{Push-LSVRG-UP} and \textit{Push-SAGA} is on a case-by-case basis. In some cases, \textit{Push-SAGA} can achieve certain accelerated convergence in terms of epochs than \textit{Push-LSVRG-UP} at the expense of expensive storage costs.

\subsection{Case Study Two: Distributed Support Vector Machine with Smoothed Hinge Loss}\label{Section 5-2}
In the second case study, we show the accelerated convergence of \textit{Push-LSVRG-UP} in contrast to the existing notable distributed optimization algorithms in a large-scale multi-agent system over undirected networks.
\begin{figure}[!htp]
\centering
\begin{minipage}{4cm}
    \includegraphics[width=4cm, height=3cm]{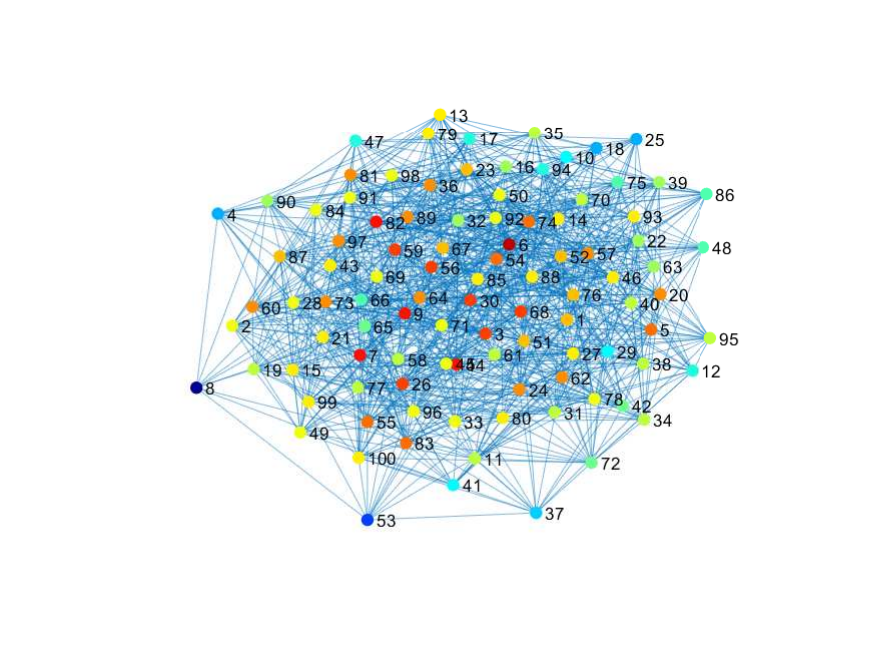}
\caption{An undirected network with $m=100$.}
\label{fig. 4}
  \end{minipage}
  \centering
\begin{minipage}{4cm}
    \includegraphics[width=4cm, height=3cm]{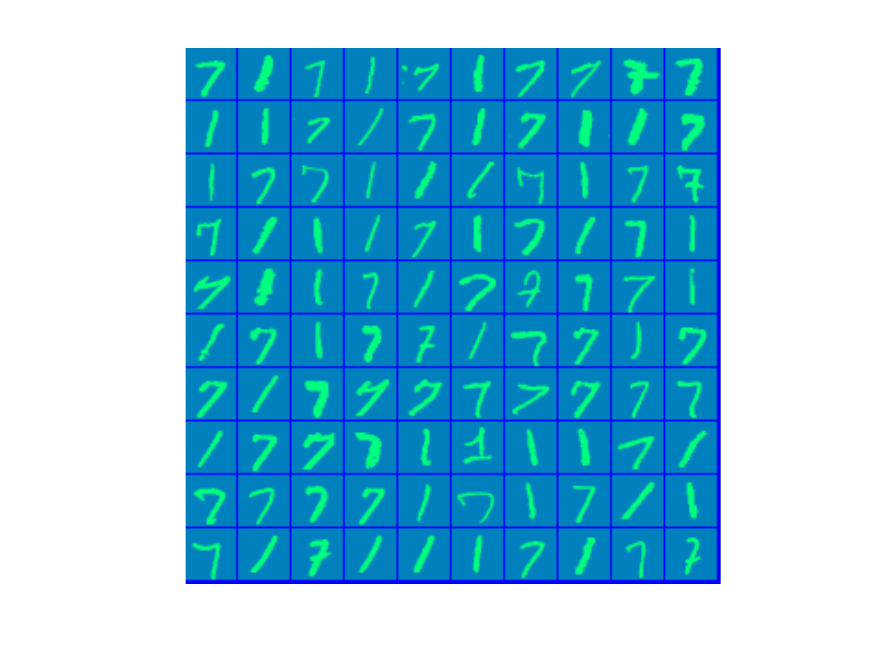}
\caption{100 random samples from MNIST dataset.}
\label{fig. 5}
  \end{minipage}
\end{figure}

Specifically, a network of $m=100$ agents cooperatively solve for a support vector machine problem to train a separating hyperplane via optimizing the following nonlinear cost function
\begin{equation}\label{E5-3}
\begin{aligned}
\mathop {\min }\limits_{\omega  \in {\mathbb{R}^n},\upsilon  \in \mathbb{R}} \tilde {f}\left( {\omega ,\upsilon } \right) :=&  \frac{1}{m}\sum\limits_{i = 1}^m {\frac{\lambda}{{{q_i}}}\sum\limits_{j = 1}^{{q_i}} {h\left( {{b_{ij}}\left( {c_{ij}^{\top} \omega  + \upsilon } \right)} \right)} } \\
&+\frac{1}{2}\left( {{\left\| \omega  \right\|_2^2} + {\upsilon ^2}} \right),
\end{aligned}
\end{equation}
where $\lambda $ is a penalty parameter in this simulation; the hinge loss function $h_i\left( u \right)$ is initially introduced in \cite{Rennie2005} as follows:
\begin{equation}\label{E5-4}
h_i\left( u \right) := \left\{ \begin{gathered}
   - 0.5 - u,\quad{\text{     if  }}u < 0, \hfill \\
  0.5{\left( {1 - u} \right)^2},{\text{ if  0}} \le u < 1, \hfill \\
  0,\qquad\qquad\;\,{\text{                if 1}} \le u. \hfill \\
\end{gathered}  \right.
\end{equation}
To solve (\ref{E5-4}) in a distributed manner, we need to define $\tilde z = {\left[ {{\omega ^{\top}},\upsilon } \right]^{\top} } \in {\mathbb{R}^{n + 1}}$ and ${{\tilde c}_{ij}} = {\left[ {c_{ij}^{\top} ,1} \right]^{\top} } \in {\mathbb{R}^{n + 1}}$. The simulation is based on a multi-agent system over a randomly-generated undirected network as depicted in Fig. \ref{fig. 4} with network connectivity ratio $0.2$.
\begin{figure}[!htp]
\centering
\subfigure[Training performance.]{\includegraphics[width=1.62in,height=1.30in]{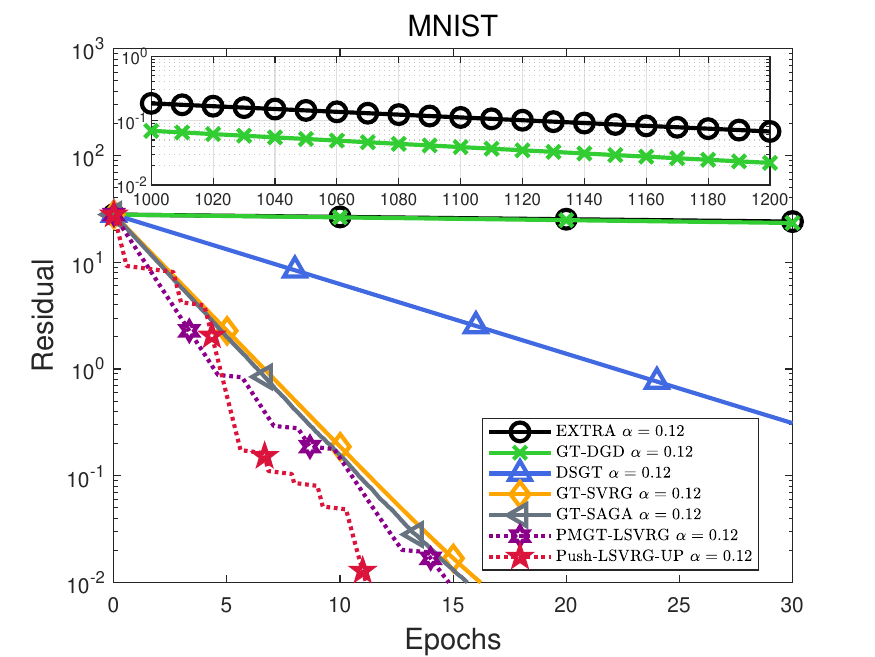}
\label{fig. 6}}
\hfil
\subfigure[Testing performance.]{\includegraphics[width=1.62in,height=1.30in]{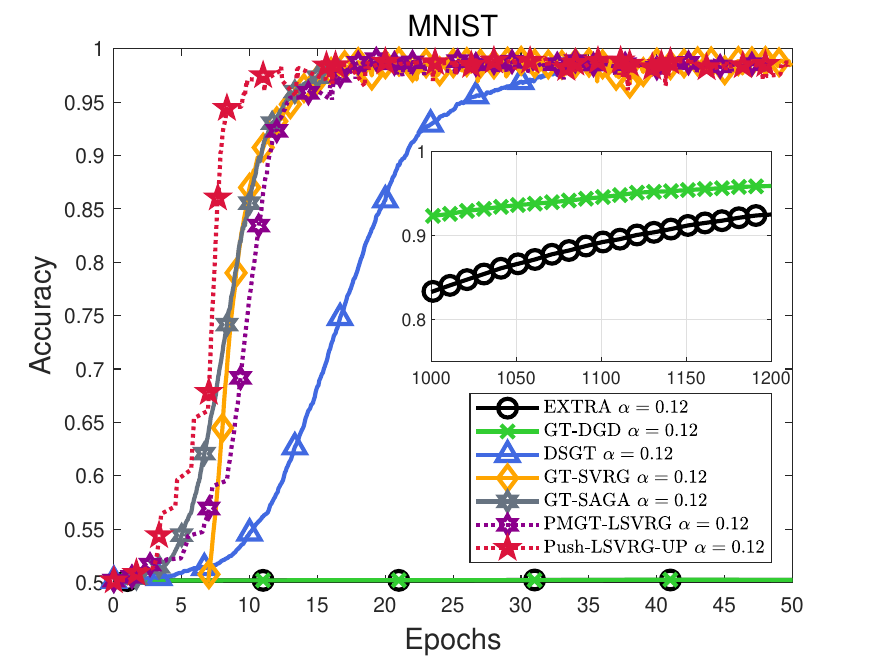}
\label{fig. 7}}
\caption{Performance comparison over epochs.}
\label{Epochs_SVM}
\end{figure}
\begin{figure}[!htp]
\centering
\subfigure[Training performance.]{\includegraphics[width=1.62in,height=1.30in]{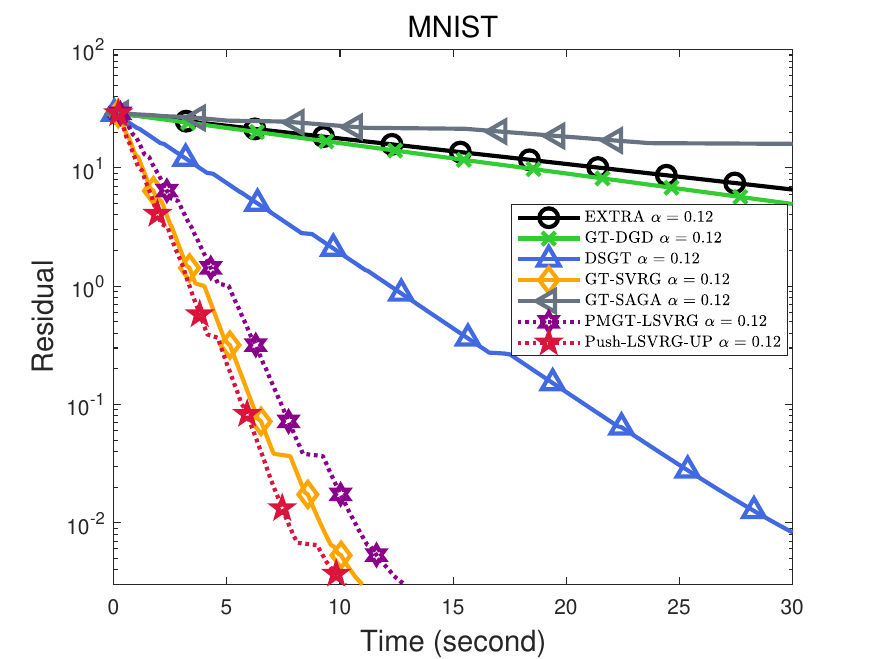}
\label{Training_time_SVM}}
\hfil
\subfigure[Testing performance.]{\includegraphics[width=1.62in,height=1.30in]{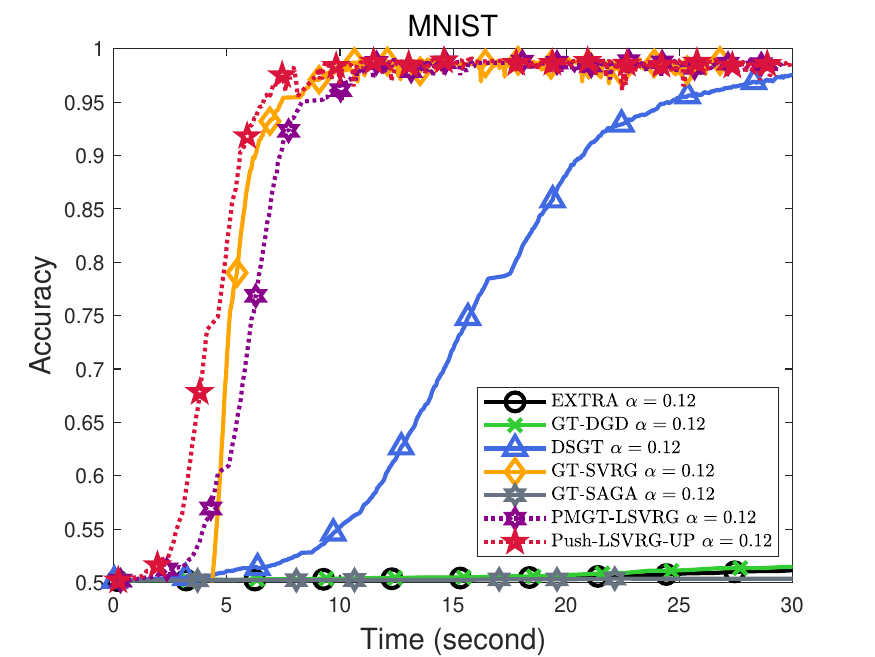}
\label{Testing_time_SVM}}
\caption{Performance comparison over CPU running time.}
\label{Time_SVM}
\end{figure}
\begin{figure}[!htp]
\centering
\subfigure[A ring network.]{\includegraphics[width=1.62in,height=1.30in]{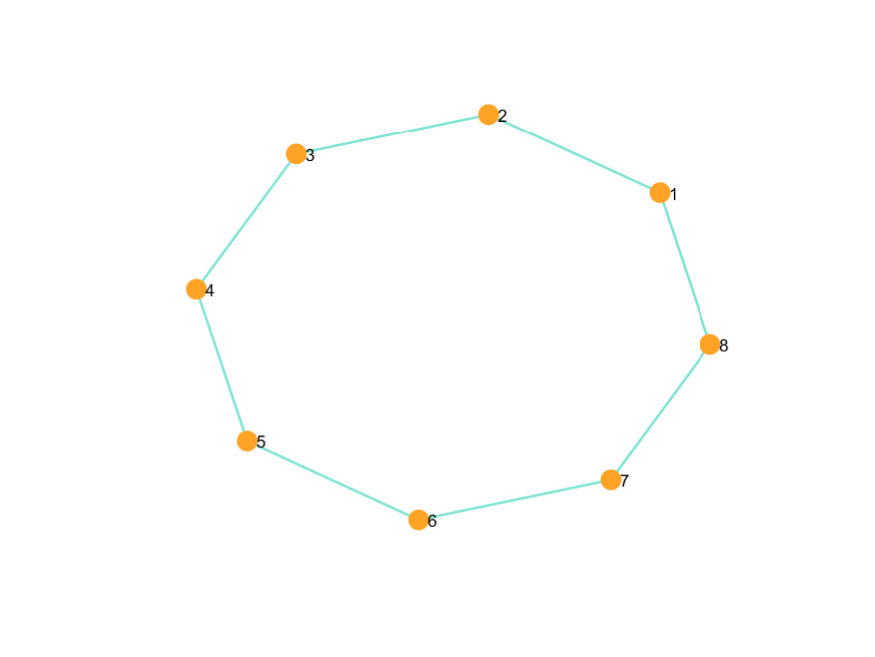}
\label{Ring_network}}
\hfil
\subfigure[A mesh network.]{\includegraphics[width=1.62in,height=1.30in]{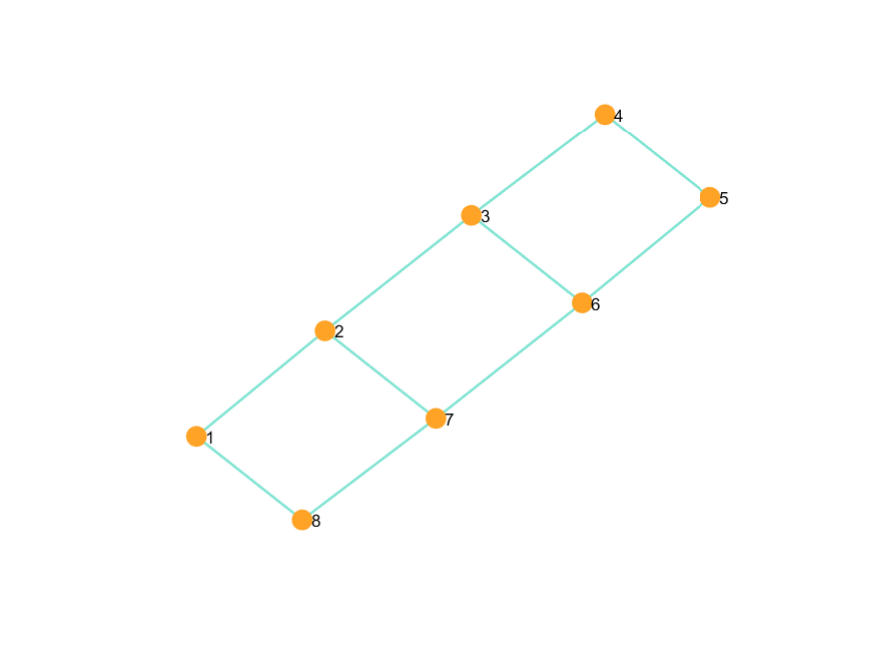}
\label{Mesh_network}}
\subfigure[A symmetric exponential network.]{\includegraphics[width=1.62in,height=1.30in]{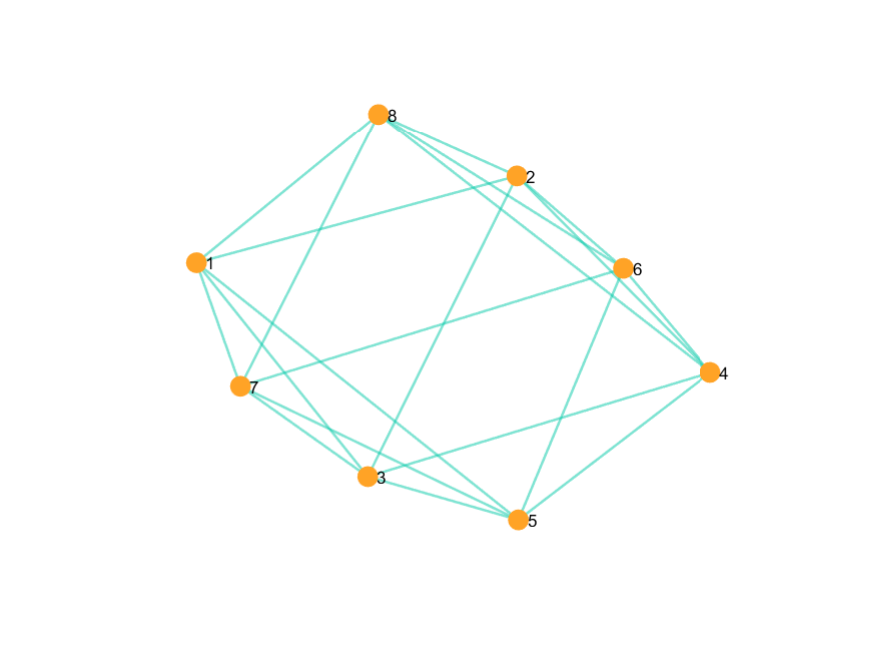}
\label{Symmetric_exponential_network}}
\hfil
\subfigure[A full-connected network.]{\includegraphics[width=1.62in,height=1.30in]{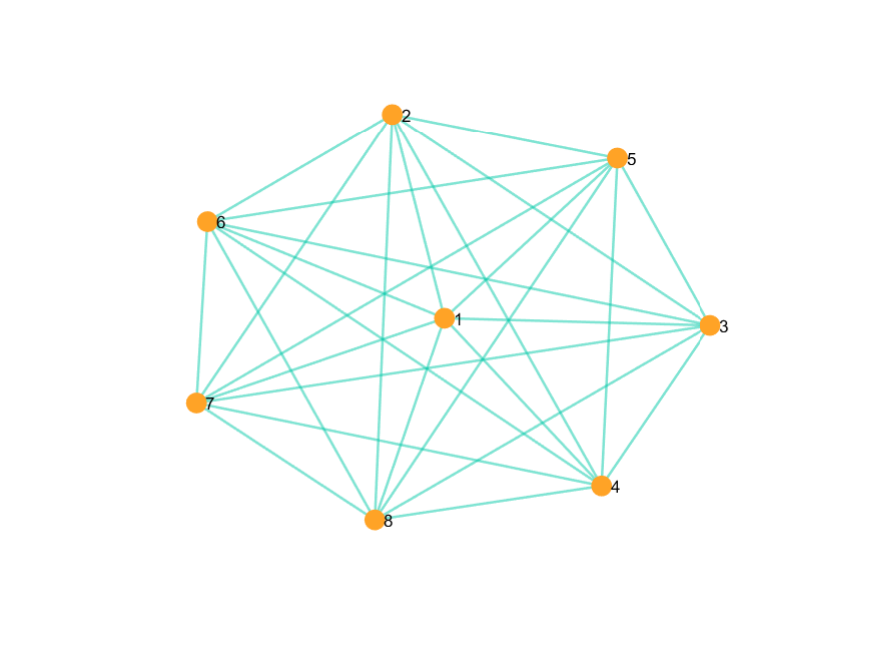}
\label{Full_connected network}}
\caption{Different network topologies with $m=8$ agents.}
\label{different networks}
\end{figure}
A total number of $12000$ samples of number $1$ and $7$ is randomly chosen from the MNIST dataset \cite{LeCun2010}, from which we randomly select $8000$ samples from the total samples for training the separating hyperplane, and the rest samples are used for testing.
Note that ${b_{ij}} =  + 1$ when the sample ${c_{ij}}$ is number $7$ while ${b_{ij}} =  - 1$ when the sample ${c_{ij}}$ is number $1$. Fig. \ref{fig. 5} visualizes 100 samples randomly selected from the MNIST dataset and each sample is quantified as a $n = 784$-dimensional vector.
We set the penalty parameter as $\lambda = 0.01$ and the total training samples are allocated evenly to each agent in Fig. \ref{fig. 4}. Therefore, each agent $i$ maintains $q_i = N/m = 80$ training samples.
Via employing a powerful "centrality" function in Matlab, Fig. \ref{fig. 4} maps different "closeness" of each agent with the other agents to various color.
From Figs. \ref{Epochs_SVM}-\ref{Time_SVM}, one can see that \textit{Push-LSVRG-UP} achieves accelerated convergence than the other tested algorithms in terms of both epochs and CPU running time over an undirected network as shown in Fig. \ref{fig. 4}. Especially when compared with the batch gradient algorithms \textit{EXTRA} \cite{Shi2015a} and \textit{GT-DGD} \cite{Xu2015b,Nedic2017a}, the accelerated performance is more obvious. Note that the inner loop number of \textit{GT-SVRG} \cite{Xin2020f} is set as $Q$.
\begin{figure}[!htp]
\centering
\subfigure[Convergence performance over epochs.]{\includegraphics[width=1.62in,height=1.30in]{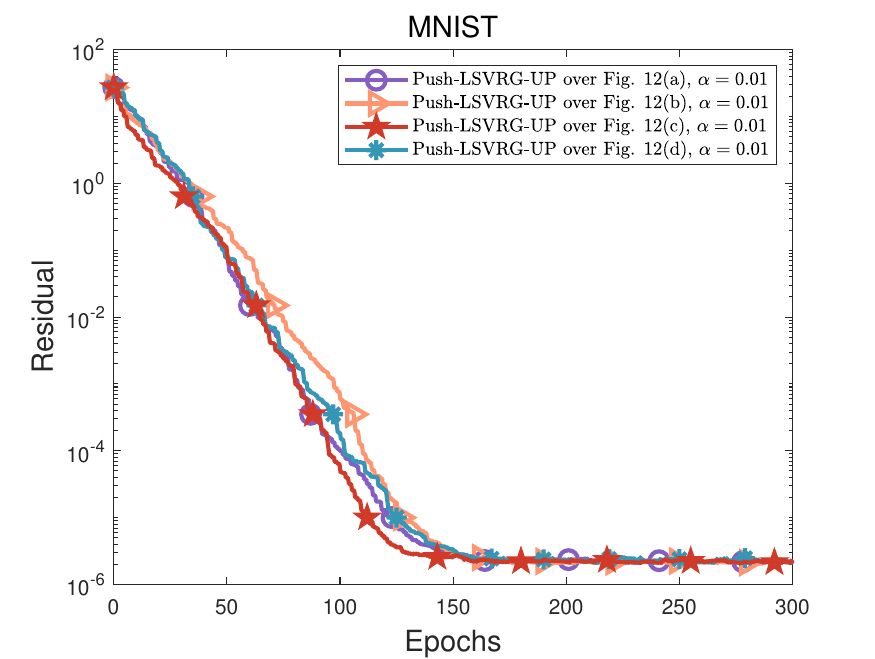}
\label{Training_epochs_topologies}}
\hfil
\subfigure[Convergence performance over CPU running time.]{\includegraphics[width=1.62in,height=1.30in]{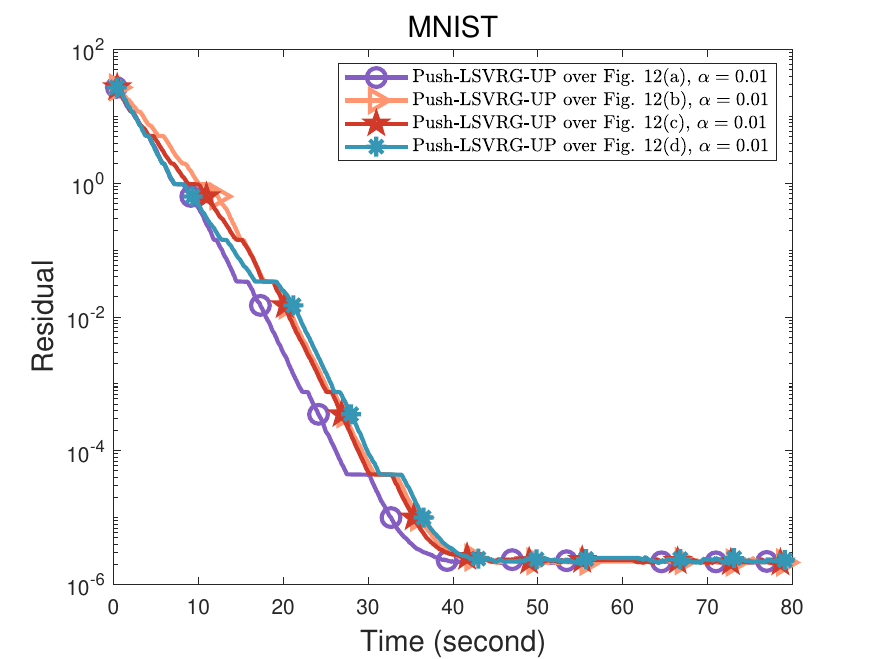}
\label{Training_time_topologies}}
\caption{Convergence performance comparison under different network topologies.}
\label{Training_topologies}
\end{figure}
\begin{figure}[!htp]
\centering
\subfigure[Convergence performance over epochs.]{\includegraphics[width=1.62in,height=1.30in]{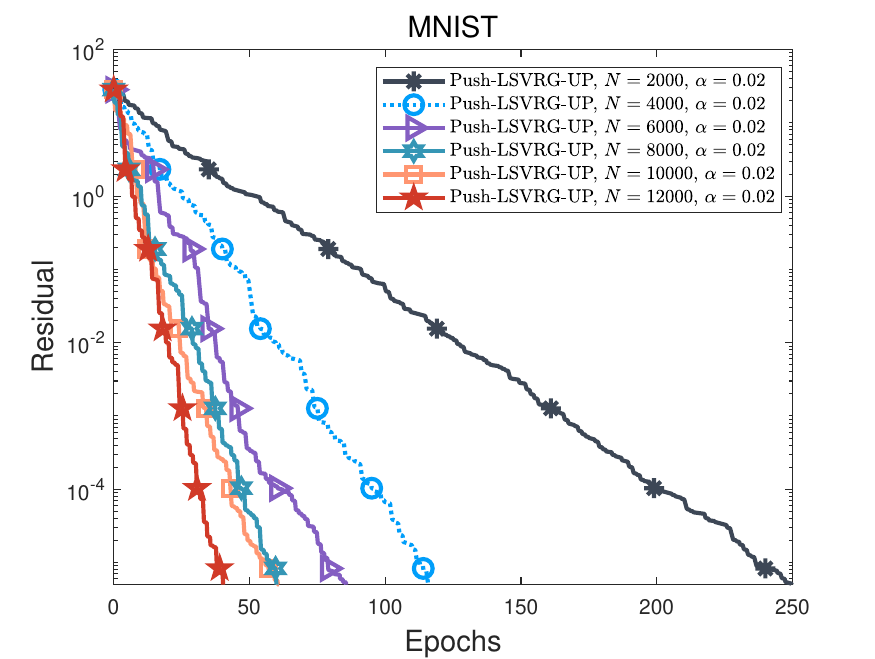}
\label{Training_epochs_samples}}
\hfil
\subfigure[Convergence performance over CPU running time.]{\includegraphics[width=1.62in,height=1.30in]{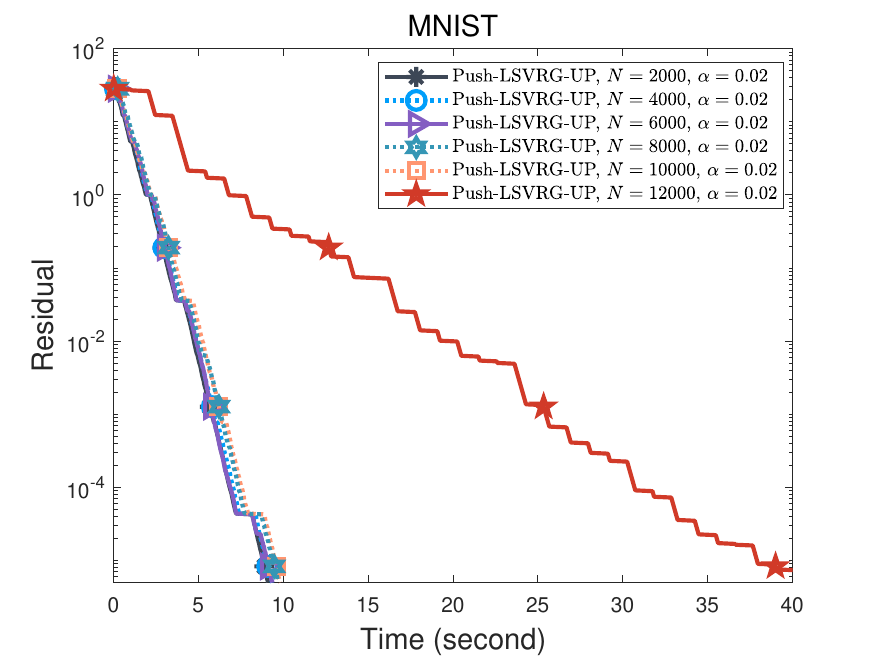}
\label{Training_time_samples}}
\caption{Convergence performance comparison under different numbers of samples.}
\label{Training_samples}
\end{figure}
To further investigate the impact of network topologies on the convergence performance of \textit{Push-LSVRG-UP}, four different network topologies are constructed as shown in Figs. \ref{Ring_network}-\ref{Full_connected network} and these network topologies are not uncommon in practice. Figs. \ref{Training_epochs_topologies}-\ref{Training_time_topologies} show the convergence of \textit{Push-LSVRG-UP} over these four different network topologies in terms of both epochs and CPU running time. From these two figures, one can see that \textit{Push-LSVRG-UP} converges faster over symmetric exponential network Fig. \ref{Symmetric_exponential_network} in terms of epochs and ring network Fig. \ref{Ring_network} in terms of CPU running time, respectively. However, the performance difference of \textit{Push-LSVRG-UP} over these four network topologies is not obvious.
To study the convergence performance of \textit{Push-LSVRG-UP} under different numbers of samples, we utilize a set of $N=2000, 4000, 6000, 8000, 10000$ samples to compare the convergence performance of \textit{Push-LSVRG-UP} over symmetric exponential network Fig. \ref{Symmetric_exponential_network}. Fig. \ref{Training_samples} shows that when the number of samples becomes larger, the convergence of \textit{Push-LSVRG-UP} regarding epochs becomes faster while the convergence of \textit{Push-LSVRG-UP} regarding CPU running time becomes slower. These results are understandable since when the number of samples becomes larger and the other parameters of \textit{Push-LSVRG-UP} remain unchanged, each epoch indicates more gradient computation and each agent needs to spend more time on computing the local batch gradients than before.

\section{Conclusions and Future Work}\label{Section 6}
In this paper, we first proposed a distributed stochastic optimization algorithm named \textit{Push-LSVRG-UP} to resolve large-scale optimization problems over unbalanced directed networks in a consensus manner. In theoretical aspects, a linear exact convergence rate, the iteration complexity and an explicit feasible step-size interval are derived, which are the first results of the \textit{LSVRG}-type method in multi-agent systems over generic unbalanced directed networks. In simulations, we provided two case studies to manifest the effectiveness and practicability of \textit{Push-LSVRG-UP}. The results of these two case studies also demonstrated the improved performance of \textit{Push-LSVRG-UP} over both undirected networks and unbalanced directed networks. However, \textit{Push-LSVRG-UP} is not perfect and the robustness of \textit{Push-LSVRG-UP} has potential to be enhanced as \textit{Push-LSVRG-UP} depends on certain network synchrony at each iteration and is also not immune to possible failures of communication links.


\appendix[Proof of Theorem \ref{Theorem 1}]\label{Appendix 1}
Before establishing the linear convergence of Algorithm \ref{Algo1}, we first show stability of the DLTI system matrix $H_{\alpha}$ in the sequel lemma.
\begin{lemma}\label{L7}
Suppose that Assumptions \ref{A1}-\ref{A3} hold. If the step-size satisfies $0 < \alpha  \le {\left( {1 - \sigma_A } \right)^2 \underline{p}}\min \left\{ {1/6\mu, 1/480\delta \mu {\mathcal{Q}^2\bar p}} \right\}$, then we have
\begin{equation}\label{E7-1}
0<\rho \left( {{H_\alpha }} \right) \le\left\| {{H_\alpha }} \right\|_\infty ^\theta  \le \eta  := 1 -  \frac{{\mu \alpha }}{4}< 1,
\end{equation}
where $\left\| H_\alpha \right\|_\infty ^\theta $ is a matrix norm induced by its corresponding max-vector norm and $\theta  = {\left[ {1,9\bar \pi \delta {\mathcal{Q}^2},60\bar \pi \delta {L^2}{\mathcal{Q}^2}\bar p/\underline{p}, 20165\bar p\vartheta \delta {L^2}{\mathcal{Q}^2}/\left( {\underline{p}}{{{\left( {1 - {\sigma_A ^2}} \right)}^2}} \right)} \right]^ \top }$ is one of the feasible choices.
\begin{proof}
To begin with, similar with \cite{Hu2021a}, we aim at solving for a proper interval of the step-size $\alpha$ and a positive vector $\theta  = {\left[ {{\theta _1},{\theta _2},{\theta _3},{\theta _4}} \right]^{\top }}$ such that ${H_\alpha }\theta  \le \eta \theta $, then $\rho \left( {{H_\alpha }} \right) \le \left\| {{H_\alpha }} \right\|_\infty ^\theta  \le \eta $ with $\eta  = 1 - {\mu \alpha }/{4}$, which is equivalent to
\begin{subequations}\label{E7-2}
\begin{align}
\label{E7-2-1}&\frac{{1 + {\sigma_A ^2}}}{2}{\theta _1} + \frac{{2{\alpha ^2}}}{{1 - {\sigma_A ^2}}}{\theta _4} \!\le\! \left( {1 - \frac{{\mu \alpha }}{4}} \right){\theta _1},\;\\
\label{E7-2-2}&\frac{{2\bar \pi \delta \alpha {L^2}}}{\mu }{\theta _1} + \left( {1 - \frac{{\mu \alpha }}{2}} \right){\theta _2} + \frac{{2{\alpha ^2}}}{m}{\theta _3} \!\le\! \left( {1 - \frac{{\mu \alpha }}{4}} \right){\theta _2},\;\\
\label{E7-2-3}&2\bar p\bar \pi {d_1}{L^2}{\theta _1} + 2m\bar p{L^2}{\theta _2} + \left( {1 - \underline{p}} \right){\theta _3} \!\le\! \left( {1 - \frac{{\mu \alpha }}{4}} \right){\theta _3},\;\\
\label{E7-2-4}&\frac{{194\delta {L^2}{\theta _1}}}{{1 - {\sigma_A ^2}}}\!\!+\!\! \frac{{169{{\underline{\pi } }^{ \!-\! 1}}{L^2}{\theta _2}}}{{1 - {\sigma_A ^2}}} \!\!+\! \frac{{110{{\underline{\pi } }^{ \!-\! 1}}{\theta _3}}}{{3( {1 \!-\! {\sigma_A ^2}} )}} \!+\!\! \frac{{3 \!+\! {\sigma_A ^2}}}{4}{\theta _4} \! \le \!\! ( {1 \!-\! \frac{{\mu \alpha }}{4}} ){\theta _4}.\;
\end{align}
\end{subequations}
To proceed, we rearrange (\ref{E7-2}) as follows:
\begin{subequations}\label{E7-3}
\begin{align}
\label{E7-3-1}&\frac{{2{\alpha ^2}}}{{1 - {\sigma_A ^2}}}{\theta _4} + \frac{{\mu \alpha }}{4}{\theta _1} \le \frac{{1 - {\sigma_A ^2}}}{2}{\theta _1},\\
\label{E7-3-2}&\frac{{2\alpha }}{m}{\theta _3} \le \frac{\mu }{4}{\theta _2} - \frac{{2\bar \pi \delta {L^2}}}{\mu }{\theta _1},\\
\label{E7-3-3}&2\bar p\bar \pi {d_1}{L^2}{\theta _1} + \frac{{\mu \alpha }}{4}{\theta _3} \le \underline{p}{\theta _3} - 2m\bar p{L^2}{\theta _2},\\
\label{E7-3-4}&\frac{{\mu \alpha }}{4}{\theta _4} \!\le\! \frac{{1 \!-\! {\sigma_A ^2}}}{4}{\theta _4} -\! \frac{{110{{\underline{\pi } }^{ - 1}}}}{{3\left( {1 \!- {\sigma_A ^2}} \right)}}{\theta _3} \!-\! \frac{{169{{\underline{\pi } }^{ - 1}}{L^2}}}{{1 -\! {\sigma_A ^2}}}{\theta _2} \!-\! \frac{{194\delta {L^2}}}{{1 \!-\! {\sigma_A ^2}}}{\theta _1}.
\end{align}
\end{subequations}
Notice that if the right-hand-side of (\ref{E7-3}) is positive, then one can always find some feasible range of step-size to satisfy the relationships. We first need to determine the positive vector $\theta$ as follows:
\begin{subequations}\label{E7-4}
\begin{align}
\label{E7-4-1}&{\theta _2} > 8\bar \pi \delta {\mathcal{Q}^2}{\theta _1}>0, \\
\label{E7-4-2}&{\theta _3} > 6{L^2}{\theta _2}\frac{{\bar p}}{{\underline{p} }}>0, \\
\label{E7-4-3}&{\theta _4} > \frac{4}{{1 \!-\! {\sigma_A ^2}}}\left( {\frac{{110{{\underline{\pi } }^{ - 1}}}}{{3\left( {1 \!-\! {\sigma_A ^2}} \right)}}{\theta _3} \!+\! \frac{{169{{\underline{\pi } }^{ - 1}}{L^2}}}{{1 \!-\! {\sigma_A ^2}}}{\theta _2} \!+\!  \frac{{194\delta {L^2}}}{{1 \!-\! {\sigma_A ^2}}}{\theta _1}} \right)>0.
\end{align}
\end{subequations}
Although there are many feasible vectors $\theta$, here we can pick a feasible one as follows: ${\theta _1} = 1$, ${\theta _2} = 9\bar \pi \delta {\mathcal{Q}^2}$, ${\theta _3} = 22m\bar \pi \delta  {L^2}{\mathcal{Q}^2}{\bar p}/{\underline{p} }$, and ${\theta _4} = {25300}\bar p\vartheta \delta {L^2}{\mathcal{Q}^2}/\left( {\underline{p} {{\left( {1 - {\sigma_A ^2}} \right)}^2}} \right)$. Then, according to (\ref{E7-3-1})-(\ref{E7-3-4}), one can find respectively
\begin{subequations}\label{E7-5}
\begin{align}
\label{E7-5-1}&0 < \alpha  \le \frac{{m\underline{p} }}{{480\bar p\mu {\mathcal{Q}^2}}},\\
\label{E7-5-2}&0 < \alpha  \le \frac{{\underline{p} \left( {1 - {\sigma_A ^2}} \right)}}{{5\mu \bar p}},\\
\label{E7-5-3}&0 < \alpha  \le \frac{{\underline{p} {{\left( {1 - {\sigma_A ^2}} \right)}^2}}}{{300\bar pL\mathcal{Q}}}\sqrt {\frac{1}{{\vartheta \delta }}},\\
\label{E7-5-4}&0 < \alpha  \le \frac{{\underline{p} }}{{6\mu }}.
\end{align}
\end{subequations}
Since $1 - \sigma_A  < 1 - {\sigma_A ^2}$, $0 < \underline{p}  \le \bar p \le 1$, and $m \ge 1$,  combining (\ref{E7-5}) with the step-size condition given in Proposition \ref{P1}, then one can choose a more tight feasible range as follows:
\begin{equation}\label{E7-6}
0 < \alpha  \le  {\left( {1 - \sigma_A } \right)^2}\underline{p}\min \left\{\frac{1}{{6\mu }}, {\frac{1}{{480\delta \mu {\mathcal{Q}^2}\bar p}}} \right\},
\end{equation}
which determines $0<\eta  = 1 -  \mu \alpha/4<1$.
\end{proof}
\end{lemma}
Base on Lemma \ref{L7}, the following lemma establishes the linear convergence of \textit{Push-LSVRG-UP} and the iteration complexity with respect to $\epsilon$-accurate solution.
\begin{lemma}\label{L8}
Suppose that Assumptions \ref{A1}-\ref{A3} hold. If one sets $0 < \alpha  \le \left( {1 - \sigma_A } \right)\underline{p} \min \left\{1/6\mu, {\left( {1 - \sigma_A } \right)/480\delta L\mathcal{Q}\bar p} \right\}$, then ${\left\| {{t_k}} \right\|_2}$ converges linearly to zero at the rate of $\mathcal{O}\left( {{{\left( {\eta  + \zeta } \right)}^k}} \right)$, where $0 < \eta  + \zeta  < 1$ and $\zeta $ is arbitrary small.
\begin{proof}
Via writing recursively (\ref{E4-44}), one can attain
\begin{equation}\label{E7-7}
{t_k} \le H_\alpha ^k{t_0} + \sum\limits_{l = 0}^{k - 1} {H_\alpha ^{k - l - 1}{G_l}{\tau _l}}.
\end{equation}
Taking the norm on the both sides of (\ref{E7-7}) yields
\begin{equation}\label{E7-8}
{\left\| {{t_k}} \right\|_2}\le {\left\| {{t_0}} \right\|_2}{\left\| {H_\alpha ^k} \right\|_2} + \sum\limits_{l = 0}^{k - 1} {{{\left\| {H_\alpha ^{k - l - 1}{G_l}} \right\|}_2}{{\left\| {{\tau _l}} \right\|}_2}}.
\end{equation}
There exist some constants ${\tilde \gamma _1} > 0$ and ${\gamma _2} > 0$ such that ${\left\| {H_\alpha ^k} \right\|_2} \le {\tilde \gamma _1}{\eta ^k}$ and ${\left\| {{G_k}} \right\|_2} = {\gamma _2}{\sigma_A ^k}$. Then, different from \cite{Qureshi2020a,Qureshi2021,Xi2018d}, here we let $0 < \alpha  \le 4\left( {1 - \sigma_A } \right)/\mu $ such that $0 < \sigma_A  \le \eta $ and further define
$\gamma  := {\tilde \gamma _1}{\gamma _2}/\eta  $ such that for $0 \le l \le k - 1$, it holds
\begin{equation}\label{E7-9}
\begin{aligned}
{\left\| {{t_k}} \right\|_2}\le& \left( {{{\left\| {{t_0}} \right\|}_2}{{\tilde \gamma }_1} + \gamma \sum\limits_{l = 0}^{k - 1} {{{\left\| {{\tau _l}} \right\|}_2}} } \right){\eta  ^k}\\
=&\left( {{\gamma _1} + \gamma \sum\limits_{l = 0}^{k - 1} {{{\left\| {{\tau _l}} \right\|}_2}} } \right){\eta  ^k},
\end{aligned}
\end{equation}
where ${\gamma _1} := {\left\| {{t_0}} \right\|_2}{{\tilde \gamma }_1}$. Moreover, it can be verified that
\begin{equation}\label{E7-10}
\begin{aligned}
\left\| {{x_l}} \right\|_2^2 =& \left\| {{x_l} - {A_\infty }{x_l} + {A_\infty }{x_l} - {Y_\infty }{z^*} + {Y_\infty }{z^*}} \right\|_2^2\\
\le& 3\bar \pi \left\| {{x_l} \!-\! {A_\infty }{x_l}} \right\|_\pi ^2 \!+\! 3m{Y^2}\left\| {{{\bar x}_l} \!- {{\tilde z}^*}} \right\|_2^2 \!+\! 3m{Y^2}\left\| {{{\tilde z}^*}} \right\|_2^2,
\end{aligned}
\end{equation}
which gives
\begin{equation}\label{E7-11}
\mathbb{E}\left[ {\left\| {{x_l}} \right\|_2^2} \right] \le 3\left( {\bar \pi  + {Y^2}} \right)\mathbb{E}\left[ {{{\left\| {{t_l}} \right\|}_2}} \right] + 3m{Y^2}\mathbb{E}\left[ {\left\| {{{\tilde z}^*}} \right\|_2^2} \right].
\end{equation}
Then, via setting $b := 3\gamma \left( {\bar \pi  + {Y^2}} \right)$ and $c := 3\gamma m{Y^2}\mathbb{E}\left[ {\left\| {{{\tilde z}^*}} \right\|_2^2} \right]$, we have
\begin{equation}\label{E7-12}
{\left\| {{t_k}} \right\|_2}{\text{ }} \le \left( {{\gamma _1} + kc + b\sum\limits_{l = 0}^{k - 1} {{{\left\| {{t_l}} \right\|}_2}} } \right){\eta ^k}.
\end{equation}
Let ${u_k} = \sum\nolimits_{l = 0}^{k - 1} {{{\left\| {{t_l}} \right\|}_2}} $, ${c_k} := \left( {{\gamma _1} + kc} \right){\eta ^k}$, and ${b_k} := b{\eta ^k}$, and then it can be verified that
\begin{equation}\label{E7-13}
{\left\| {{t_k}} \right\|_2} = {u_{k + 1}} - {u_k} \le \left( {{\gamma _1} + kc + b{u_k}} \right){\eta ^k},
\end{equation}
which gives
\begin{equation}\label{E7-14}
{u_{k + 1}} \le \left( {1 + {b_k}} \right){u_k} + {c_k}.
\end{equation}
Since ${\left\{ {{u_k}} \right\}_{k \ge 0}}$, ${\left\{ {{b_k}} \right\}_{k \ge 0}}$, and ${\left\{ {{c_k}} \right\}_{k \ge 0}}$ are nonnegative sequences, together with $\sum\nolimits_{k = 0}^\infty  {{b_k}}  < \infty $ and $\sum\nolimits_{k = 0}^\infty  {{c_k}}  < \infty $, it follows from \cite[Lemma 7]{Nedic2015f} that the sequence ${\left\{ {{u_k}} \right\}_{k \ge 0}}$ converges and is thus bounded. Therefore, via choosing $\forall \varpi  \in \left( {\eta ,1} \right)$, one can obtain
\begin{equation}\label{E7-15}
\mathop {\lim }\limits_{k \to \infty } \frac{{{{\left\| {{t_k}} \right\|}_2}}}{{{\varpi ^k}}} \le \mathop {\lim }\limits_{k \to \infty } \frac{{\left( {{\gamma _1} + kc + b{u_k}} \right){\eta ^k}}}{{{\varpi ^k}}} = 0.
\end{equation}
That is to say ${\left\| {{t_k}} \right\|_2} = \mathcal{O}\left( {{\varpi ^k}} \right)$, which means that there exists some constant $\chi  > 0$ such that $\forall k \ge 0$,
\begin{equation}\label{E7-16}
{\left\| {{t_k}} \right\|_2} \le \chi {\left( {\eta  + \zeta } \right)^k},
\end{equation}
where $\zeta $ is an arbitrarily small positive constant satisfying $0 < \eta  + \zeta  < 1$. One can utilize the above results to derive the iteration complexity of \textit{Push-LSVRG-UP} in the following. To begin with, recalling (\ref{E3-1}), we have
\begin{equation}\label{E7-17}
\begin{aligned}
&\mathbb{E}\left[ {\left\| {{z_k} - {{1_m} \otimes {{\tilde z}^*}} } \right\|_2^2} \right]\\
\le&3{{\tilde Y}^2}\mathbb{E}\left[ {\left\| {{x_k} - {A_\infty }{x_k}} \right\|_2^2} \right] + 3{{\tilde Y}^2}{Y^2}\mathbb{E}\left[ {m\left\| {{{\bar x}_k} - {{\tilde z}^*}} \right\|_2^2} \right]\\
&+ 3m\left\| {Y_k^{ - 1}{Y_\infty } - {I_{mn}}} \right\|_2^2\mathbb{E}\left[ {\left\| {{{\tilde z}^*}} \right\|_2^2} \right]\\
\le& 3\bar \pi {{\tilde Y}^2}{\left\| {{t_k}} \right\|_2} + 3{{\tilde Y}^2}{Y^2}{\left\| {{t_k}} \right\|_2} + 3m{Y^2}{T^2}{\sigma_A ^{2k}}\mathbb{E}\left[ {\left\| {{{\tilde z}^*}} \right\|_2^2} \right]\\
\le& 3{{\tilde Y}^2}\left( {\bar \pi  + {Y^2}} \right)\chi {\left( {\eta  + \zeta } \right)^k} + 3m{Y^2}{T^2}{\left( {\eta  + \zeta } \right)^k}\mathbb{E}\left[ {\left\| {{{\tilde z}^*}} \right\|_2^2} \right],
\end{aligned}
\end{equation}
where the last inequality uses the fact that $0 < \sigma_A  \le \eta  < 1$. Via defining $\varphi := {3{{\tilde Y}^2}\left( {\bar \pi  + {Y^2}} \right)\chi  + 3m{Y^2}{T^2}\mathbb{E}\left[ {\left\| {{{\tilde z}^*}} \right\|_2^2} \right]} $, we have
\begin{equation}\label{E7-18}
\mathbb{E}\left[ {\left\| {{z_k} - {{1_m} \otimes {{\tilde z}^*}} } \right\|_2^2} \right] \le \varphi {\left( {\eta  + \zeta } \right)^k}.
\end{equation}
Then, to attain an $\epsilon$-accurate solution, i.e., $\mathbb{E}\left[ {\left\| {{z_k} - {{1_m} \otimes {{\tilde z}^*}} } \right\|_2^2} \right] \le \epsilon$, one needs
\begin{equation}\label{E7-19}
\begin{aligned}
\mathbb{E}\left[ {\left\| {{z_k} - {{1_m} \otimes {{\tilde z}^*}} } \right\|_2^2} \right] \le \varphi {\left( {1 - \left( {1 - \left( {\eta  + \zeta } \right)} \right)} \right)^k}\mathop  \le \epsilon.
\end{aligned}
\end{equation}
For $\phi  \in \left( {0,1} \right)$, we always have $1 - \phi  \le {e^{ - \phi }}$. Therefore, one can solve for a sufficient condition as follows:
\begin{equation}\label{E7-30}
 {\left( {1 - \left( {1 - \left( {\eta  + \zeta } \right)} \right)} \right)^k} \le {e^{ - \left( {1 - \left( {\eta  + \zeta } \right)} \right)k}} \le \frac{\epsilon}{\varphi },
\end{equation}
which leads to (\ref{E4-2}) according to (\ref{E7-1}).
\end{proof}
\end{lemma}

\ifCLASSOPTIONcompsoc
  \section*{Acknowledgments}
\else
  \section*{Acknowledgment}
\fi
The authors would like to thank Prof. Lijun Zhu and Miss Xi Zhang for the invaluable assistance and effort to help us complete this paper. The authors would also like to appreciate the UCI Machine Learning Repository \cite{Dua2013a} and \cite{LeCun2010} for the wonderful experimental dataset supports.

\ifCLASSOPTIONcaptionsoff
  \newpage
\fi



%
%
%
\bibliographystyle{IEEEtran}
\bibliography{Push-LSVRG-UP}%

\begin{thebibliography}{10}
\providecommand{\url}[1]{#1}
\csname url@samestyle\endcsname
\providecommand{\newblock}{\relax}
\providecommand{\bibinfo}[2]{#2}
\providecommand{\BIBentrySTDinterwordspacing}{\spaceskip=0pt\relax}
\providecommand{\BIBentryALTinterwordstretchfactor}{4}
\providecommand{\BIBentryALTinterwordspacing}{\spaceskip=\fontdimen2\font plus
\BIBentryALTinterwordstretchfactor\fontdimen3\font minus
  \fontdimen4\font\relax}
\providecommand{\BIBforeignlanguage}[2]{{%
\expandafter\ifx\csname l@#1\endcsname\relax
\typeout{** WARNING: IEEEtran.bst: No hyphenation pattern has been}%
\typeout{** loaded for the language `#1'. Using the pattern for}%
\typeout{** the default language instead.}%
\else
\language=\csname l@#1\endcsname
\fi
#2}}
\providecommand{\BIBdecl}{\relax}
\BIBdecl

\bibitem{Boyd2010b}
S.~Boyd, ``{Distributed optimization and statistical learning via the
  alternating direction method of multipliers},'' \emph{Foundations and
  Trends{\textregistered} in Machine Learning}, vol.~3, no.~1, pp. 1--122,
  2010.

\bibitem{Nedic2020}
N.~Angelia, ``{Distributed gradient methods for convex machine learning
  problems in networks: Distributed optimization},'' \emph{IEEE Signal
  Processing Magazine}, vol.~37, no.~3, pp. 92--101, 2020.

\bibitem{Assran2019}
M.~Assran, N.~Loizou, N.~Ballas, and M.~Rabbat, ``{Stochastic gradient push for
  distributed deep learning},'' in \emph{36th International Conference on
  Machine Learning}, 2019, pp. 344--353.

\bibitem{Li2021}
H.~Li, J.~Hu, L.~Ran, Z.~Wang, Q.~L{\"{u}}, Z.~Du, and T.~Huang,
  ``{Decentralized dual proximal gradient algorithms for non-smooth constrained
  composite optimization problems},'' \emph{IEEE Transactions on Parallel and
  Distributed Systems}, vol.~32, no.~10, pp. 2594--2605, 2021.

\bibitem{Yang2019}
T.~Yang, D.~Wu, H.~Fang, W.~Ren, H.~Wang, Y.~Hong, and K.~H. Johansson,
  ``{Distributed energy resource coordination over time-varying directed
  communication networks},'' \emph{IEEE Transactions on Control of Network
  Systems}, vol.~6, no.~3, pp. 1124--1134, 2019.

\bibitem{Li2019i}
Z.~Li, W.~Shi, and M.~Yan, ``{A decentralized proximal-gradient method with
  network independent step-sizes and separated convergence rates},'' \emph{IEEE
  Transactions on Signal Processing}, vol.~67, no.~17, pp. 4494--4506, 2019.

\bibitem{Lu2020e}
Q.~L{\"{u}}, X.~Liao, H.~Li, and T.~Huang, ``{Achieving acceleration for
  distributed economic dispatch in smart grids over directed networks},''
  \emph{IEEE Transactions on Network Science and Engineering}, vol.~7, no.~3,
  pp. 1988--1999, 2020.

\bibitem{Shi2018d}
X.~Shi, J.~Cao, and W.~Huang, ``{Distributed parametric consensus optimization
  with an application to model predictive consensus problem},'' \emph{IEEE
  Transactions on Cybernetics}, vol.~48, no.~7, pp. 2024--2035, 2018.

\bibitem{Camponogara2011}
E.~Camponogara and H.~F. Scherer, ``{Distributed optimization for model
  predictive control of linear dynamic networks with control-input and output
  constraints},'' \emph{IEEE Transactions on Automation Science and
  Engineering}, vol.~8, no.~1, pp. 233--242, 2011.

\bibitem{Xin2020g}
R.~Xin, S.~Kar, and U.~A. Khan, ``{Decentralized stochastic optimization and
  machine learning: A unified variance-reduction framework for robust
  performance and fast convergence},'' \emph{IEEE Signal Processing Magazine},
  vol.~37, no.~3, pp. 102--113, 2020.

\bibitem{Ye2020}
H.~Ye, W.~Xiong, and T.~Zhang, ``{PMGT-VR: A decentralized proximal-gradient
  algorithmic framework with variance reduction},'' \emph{arXiv preprint
  arXiv:2012.15010}, pp. 1--29, 2020.

\bibitem{Xin2020f}
R.~Xin, U.~A. Khan, and S.~Kar, ``{Variance-reduced decentralized stochastic
  optimization with accelerated convergence},'' \emph{IEEE Transactions on
  Signal Processing}, vol.~68, pp. 6255--6271, 2020.

\bibitem{Pu2021a}
S.~Pu and A.~Nedi{\'{c}}, ``{A distributed stochastic gradient tracking
  method},'' in \emph{2018 IEEE Conference on Decision and Control (CDC)},
  2018, pp. 963--968.

\bibitem{Nedicr2009}
N.~Angelia and O.~Asuman, ``{Distributed subgradient methods for multi-agent
  optimization},'' \emph{IEEE Transactions on Automatic Control}, vol.~54,
  no.~1, pp. 48--61, 2009.

\bibitem{J.C.DuchiA.Agarwal}
J.~C. Duchi, A.~Agarwal, and M.~J. Wainwright, ``{Dual averaging for
  distributed optimization: Convergence analysis and network scaling},''
  \emph{IEEE Transactions on Automatic Control}, vol.~57, no.~3, pp. 592--606,
  2012.

\bibitem{Shi2015a}
W.~Shi, Q.~Ling, G.~Wu, and W.~Yin, ``{EXTRA: An exact first-order algorithm
  for decentralized consensus optimization},'' \emph{SIAM Journal on
  Optimization}, vol.~25, no.~2, pp. 944--966, 2015.

\bibitem{Nedic2017a}
N.~Angelia, O.~Alex, and W.~Shi, ``{Achieving geometric convergence for
  distributed optimization over time-varying graphs},'' \emph{SIAM Journal on
  Optimization}, vol.~27, no.~4, pp. 2597--2633, 2017.

\bibitem{Xu2015b}
J.~Xu, S.~Zhu, Y.~C. Soh, and L.~Xie, ``{Augmented distributed gradient methods
  for multi-agent optimization under uncoordinated constant stepsizes},'' in
  \emph{Proceedings of the IEEE Conference on Decision and Control}, 2015, pp.
  2055--2060.

\bibitem{Jakovetic2019a}
D.~Jakoveti{\'{c}}, ``{A unification and generalization of exact distributed
  first-order methods},'' \emph{IEEE Transactions on Signal and Information
  Processing over Networks}, vol.~5, no.~1, pp. 31--46, 2019.

\bibitem{Xi2018e}
C.~Xi, V.~S. Mai, R.~Xin, E.~H. Abed, and U.~A. Khan, ``{Linear convergence in
  optimization over directed graphs with row-stochastic matrices},'' \emph{IEEE
  Transactions on Automatic Control}, vol.~63, no.~10, pp. 3558--3565, 2018.

\bibitem{Xi2016}
C.~Xi and U.~A. Khan, ``{Distributed subgradient projection algorithm over
  directed graphs},'' \emph{IEEE Transactions on Automatic Control}, vol.~62,
  no.~8, pp. 3986--3992, 2016.

\bibitem{Xi2017h}
C.~Xi and U.~A. Khan, ``{DEXTRA: A fast algorithm for optimization over directed graphs},''
  \emph{IEEE Transactions on Automatic Control}, vol.~62, no.~10, pp.
  4980--4993, 2017.

\bibitem{Nedic2015f}
N.~Angelia and O.~Alex, ``{Distributed optimization over time-varying directed
  graphs},'' \emph{IEEE Transactions on Automatic Control}, vol.~60, no.~3, pp.
  601--615, 2015.

\bibitem{Xi2018d}
C.~Xi, R.~Xin, and U.~A. Khan, ``{ADD-OPT: Accelerated distributed directed
  optimization},'' \emph{IEEE Transactions on Automatic Control}, vol.~63,
  no.~5, pp. 1329--1339, 2018.

\bibitem{Saadatniaki2020}
F.~Saadatniaki, R.~Xin, and U.~A. Khan, ``{Decentralized optimization over
  time-varying directed graphs with row and column-stochastic matrices},''
  \emph{IEEE Transactions on Automatic Control}, vol.~65, no.~11, pp.
  4769--4780, 2020.

\bibitem{Hu2021c}
J.~Hu, Y.~Yan, H.~Li, Z.~Wang, D.~Xia, and J.~Guo, ``{Convergence of an
  accelerated distributed optimisation algorithm over time-varying directed
  networks},'' \emph{IET Control Theory {\&} Applications}, vol.~15, no.~1, pp.
  24--39, 2021.

\bibitem{Pu2021}
S.~Pu, W.~Shi, J.~Xu, and A.~Nedic, ``{Push-Pull gradient methods for
  distributed optimization in networks},'' \emph{IEEE Transactions on Automatic
  Control}, vol.~66, no.~1, pp. 1--16, 2021.

\bibitem{Defazio2014c}
A.~Defazio, F.~Bach, and S.~Lacoste-Julien, ``{SAGA: A fast incremental
  gradient method with support for non-strongly convex composite objectives},''
  in \emph{Advances in Neural Information Processing Systems}, 2014, pp.
  1646--1654.

\bibitem{Ribeiro2015}
A.~Mokhtari and A.~Ribeiro, ``{DSA: Decentralized double stochastic averaging
  gradient algorithm},'' \emph{Journal of Machine Learning Research}, vol.~17,
  no.~1, pp. 2165--2199, 2016.

\bibitem{Wang2019}
Z.~Wang and H.~Li, ``{Edge-based stochastic gradient algorithm for distributed
  optimization},'' \emph{IEEE Transactions on Network Science and Engineering},
  vol.~7, no.~3, pp. 1421--1430, 2019.

\bibitem{Li2022}
H.~Li, L.~Zheng, Z.~Wang, Y.~Yan, L.~Feng, and J.~Guo, ``{S-DIGing : A
  stochastic gradient tracking algorithm for distributed optimization},''
  vol.~6, no.~1, pp. 53--65, 2022.

\bibitem{Ye2020a}
\BIBentryALTinterwordspacing
H.~Ye, L.~Luo, Z.~Zhou, and T.~Zhang, ``{Multi-consensus decentralized
  accelerated gradient descent},'' pp. 1--31, 2020. [Online]. Available:
  \url{http://arxiv.org/abs/2005.00797}
\BIBentrySTDinterwordspacing

\bibitem{Johnson2013c}
R.~Johnson and T.~Zhang, ``{Accelerating stochastic gradient descent using
  predictive variance reduction},'' in \emph{Advances in Neural Information
  Processing Systems}, 2013, pp. 315--323.

\bibitem{Kovalev2019}
D.~Kovalev, S.~Horvath, and P.~Richtarik, ``{Don't jump through hoops and
  remove those loops: SVRG and Katyusha are better without the outer loop},''
  in \emph{Algorithmic Learning Theory}, 2020, pp. 451--467.

\bibitem{Yuan2019e}
K.~Yuan, B.~Ying, J.~Liu, and A.~H. Sayed, ``{Variance-reduced stochastic
  learning by networked agents under random reshuffling},'' vol.~67, no.~2, pp.
  351--366, 2019.

\bibitem{Qureshi2020a}
M.~I. Qureshi, R.~Xin, S.~Kar, and U.~A. Khan, ``{Push-SAGA: A decentralized
  stochastic algorithm with variance reduction over directed graphs},''
  \emph{IEEE Control Systems Letters}, vol.~6, pp. 1202--1207, 2021.

\bibitem{Hu2021a}
J.~Hu, X.~Chen, L.~Zheng, L.~Zhang, and H.~Li, ``{The Barzilai-Borwein Method
  for distributed optimization over unbalanced directed networks},''
  \emph{Engineering Applications of Artificial Intelligence}, vol.~99, p.
  104151, 2021.

\bibitem{Nedic2016}
N.~Angelia and O.~Alex, ``{Stochastic gradient-push for strongly convex
  functions on time-varying directed graphs},'' \emph{IEEE Transactions on
  Automatic Control}, vol.~61, no.~12, pp. 3936--3947, 2016.

\bibitem{Xin2019c}
R.~Xin, A.~K. Sahu, U.~A. Khan, and S.~Kar, ``{Distributed stochastic
  optimization with gradient tracking over strongly-connected networks},'' in
  \emph{Proceedings of the IEEE Conference on Decision and Control}, 2019, pp.
  8353--8358.

\bibitem{Qureshi2021}
M.~I. Qureshi, R.~Xin, S.~Kar, and U.~A. Khan, ``{S-ADDOPT: Decentralized
  stochastic first-order optimization over directed graphs},'' \emph{IEEE
  Control Systems Letters}, vol.~5, no.~3, pp. 953--958, 2021.

\bibitem{Johnson2013}
R.~Horn and C.~R. Johnson, \emph{{Matrix analysis, 2nd ed.}}\hskip 1em plus
  0.5em minus 0.4em\relax Cambridge University Press, 2012.

\bibitem{Bubeck2015c}
S.~Bubeck, ``{Convex optimization: Algorithms and complexity},''
  \emph{Foundations and Trends{\textregistered} in Machine Learning}, vol.~8,
  no. 3-4, pp. 231--357, 2015.

\bibitem{Gorbunov2019}
E.~Gorbunov, F.~Hanzely, and P.~Richt{\'{a}}rik, ``{A unified theory of SGD:
  Variance reduction, sampling, quantization and coordinate descent},'' in
  \emph{International Conference on Artificial Intelligence and Statistics},
  2020, pp. 680--690.

\bibitem{Dua2013a}
D.~Dua and Graff, \emph{{UCI Machine Learning Repository
  [http://archive.ics.uci.edu/ml]. Irvine, CA: University of California, School
  of Information and Computer Science.}}, 2019.

\bibitem{Yuan2021}
K.~Yuan, Y.~Chen, X.~Huang, Y.~Zhang, P.~Pan, Y.~Xu, W.~Yin, and A.~Group,
  ``{DecentLaM: Decentralized momentum SGD for large-batch deep training},'' in
  \emph{Proceedings of the IEEE/CVF International Conference on Computer
  Vision}, 2021, pp. 3029--3039.

\bibitem{Rennie2005}
J.~Rennie and N.~Srebro, ``{Loss functions for preference levels: Regression
  with discrete ordered labels},'' in \emph{Workshop on Advances in
  Preference}, 2005, pp. 180--186.

\bibitem{LeCun2010}
Y.~LeCun, C.~Cortes, and C.~Burges, ``{MNIST handwritten digit database.
  [Online]. Available: http://yann. lecun.com/exdb/m},'' in \emph{AT{\&}T Labs,
  Florham Park, NJ, USA.}, 2020.

\end{thebibliography}

\end{document}